\documentclass[11pt,leqno]{article}

\usepackage{amsmath}
\usepackage{amssymb}
\usepackage{dsfont}
\usepackage{epsfig}

\usepackage{mathrsfs}

\usepackage{float}

\numberwithin{equation}{section}

\newtheorem{Theorem}{Theorem}
\newtheorem{Corollary}[Theorem]{Corollary}
\newtheorem{Lemma}[Theorem]{Lemma}
\newtheorem{Proposition}[Theorem]{Proposition}
\newtheorem{Definition}[Theorem]{Definition}

\numberwithin{Theorem}{section}

\bibstyle{plain}

\newcommand{\al}{\alpha}
\newcommand{\R}{{\mathbf R}}
\newcommand{\ds}{\displaystyle}
\newcommand{\e}{\varepsilon}
\newcommand{\Om}{\Omega}
\newcommand{\lra}{\longrightarrow}
\newcommand{\ra}{\rightarrow}
\newcommand{\p}{\partial}
\newcommand{\la}{\lambda}
\newcommand{\g}{\gamma}
\newcommand{\ov}{\overline}
\newcommand{\C}{{\mathbf C}}

\newcommand{\N}{{\mathbf N}}
\newcommand{\Z}{{\mathbf Z}}

\newcommand{\1}{{\mathds 1}}

\newcommand{\si}{\sigma}
\newcommand{\te}{\theta}
\newcommand{\de}{\delta}

\newcommand{\ph}{\varphi}

\newcommand{\Xo}{{\mathcal X}}
\newcommand{\Yo}{{\mathcal Y}}
\newcommand{\So}{{\mathcal S}}

\newcommand{\Lo}{{\mathcal L}}
\newcommand{\Do}{{\mathcal D}}
\newcommand{\Co}{{\mathcal C}}

\newcommand{\DR}{{\mathcal D}^{1,2}(\R^N)}
\newcommand{\D}{{\mathcal D}^{1,2}}
\newcommand{\wh}{\widehat}

\newcommand{\I}{\mbox{\rm Im}}
\newcommand{\Real}{\mbox{\rm Re}}

\newcommand{\beq}{\begin{equation}}
\newcommand{\eeq}{\end{equation}}

\newcommand{\references}[1]{\theinstitutions 
\begin{thebibliography}{#1}}

\newcommand{\Endrefs}{\end{thebibliography}}

\title{ Traveling waves for nonlinear Schr\"odinger equations with 
nonzero conditions at infinity}

\date{  }

\author{ Mihai MARI\c S \\
          { \small     } \\
         \it Institut de Math\'ematiques de Toulouse UMR 5219 \\
         \it Universit\'e Paul Sabatier \\
         \it 118, Route de Narbonne,   31062 Toulouse, France \\
         \it  e-mail: mihai.maris@math.univ-toulouse.fr
          }

\date{}

\begin{document}

\maketitle

\begin{center}
\it 
Dedicated to Jean-Claude Saut, 

who gave me water to cross the desert.

\end{center}
\rm

\begin{abstract}
For a large class of nonlinear Schr\"odinger equations with 
nonzero conditions at infinity 
and for any speed $c$ less than the sound velocity, 
we prove the existence of nontrivial 
finite energy traveling waves 
moving with speed $c$ in any space dimension $N\geq 3$.  
Our results are valid  as well for the Gross-Pitaevskii equation and 
for NLS with cubic-quintic nonlinearity.


\smallskip


\end{abstract}

\section{Introduction}

We consider the 
nonlinear Schr\"odinger equation
\beq
\label{1.1}
i \frac{ \p \Phi}{\p t} + \Delta \Phi + F( | \Phi |^2) \Phi = 0
\qquad \mbox{ in } \R^N,
\eeq
where $ \Phi : \R^N \times \R \lra \C$ satisfies the "boundary condition" 
$|\Phi | \lra r_0 $ as $ |x| \lra \infty $, 
$r_0 >0$ and $F$ is a real-valued function on $\R_+$ satisfying $ F(r_0^2) = 0$. 

Equations of the form (\ref{1.1}),  
with the considered non-zero conditions at infinity,  arise in a 
large variety of physical problems.   They have been used as  models for  superconductivity,
superfluid Helium II and for Bose-Einstein condensation
(\cite{barashenkov2}, \cite{barashenkov1}, \cite{berloff}, \cite{coste}, \cite{GR}, 
\cite{gross}, \cite{IS}, \cite{JR}, \cite{JPR}).
In nonlinear optics, they appear in the context of dark solitons, which are localized nonlinear waves
(also called "holes") moving on a stable continuous background (see \cite{KL}, \cite{KPS}). 
The boundary condition $|\Phi | \lra r_0 $ at infinity is precisely due to the nonzero background. 

Two important particular cases of (1.1) have been extensively studied 
by physicists and by mathematicians:
the Gross-Pitaevskii (GP) equation (where $F(s)=  1-s$)
and the so-called "cubic-quintic"  Schr\"odinger equation (where
$F(s) = - \al _1 + \al _3 s - \al _5 s ^2 $, $\; \al_1, \, \al _3, \, \al _5$
are positive and  $F$ has two positive roots). 
In both cases we have $ F'( r_0 ^2) < 0$, which ensures that (\ref{1.1}) is defocusing.

The boundary condition $|\Phi | \lra r_0 >0$ at infinity makes the structure of solutions of (\ref{1.1})
much more complicated than in the usual case of zero boundary conditions  
(when the associated dynamics   is essentially governed by dispersion and scattering).
This fact is confirmed by the existence of a remarkable variety of special solutions, such as 
traveling waves or vortex solutions,  
and regimes, like the long wave or the transonic limit.


Using the Madelung transformation
$\Phi(x,t) = \sqrt{\rho(x,t)} e^{i\theta (x,t )}$
(which is well-defined whenever $ \Phi \neq 0$), 
equation (\ref{1.1}) is equivalent to a system of Euler's equations 
for a compressible inviscid fluid of density $ \rho $ and velocity 
$ 2 \nabla \theta$. 
In this context it has been shown that, if $F$ is $C^1$ near $ r_0^2$ 
and $ F'(r_0^2) < 0$, the sound velocity at infinity associated to (\ref{1.1}) is 
$ v_s = r_0 \sqrt{ - 2 F'(r_0^2) } $ (see the introduction of \cite{M8}).

In the defocusing case  $ F'( r_0 ^2) < 0$, 
we perform a simple scaling ($ \Phi (x, t) = r_0 \tilde{\Phi}(\tilde{x}, \tilde{t})$, where
 $ \tilde{x} = r_0 \sqrt{ - F'( r_0 ^2)} x$, $\tilde{t} = - r_0 ^2 F'( r_0 ^2) t$,  and 
$ \tilde{F}(s) = \frac{-1}{r_0 ^2 F'( r_0 ^2)} F( r_0 ^2 s)$) 
and we  assume from now on  that $ r_0 = 1$ and $ F'(r_0^2) = -1$. 
The sound velocity at infinity becomes then $ v_s = \sqrt{2}$.

Equation (\ref{1.1}) is Hamiltonian: 
denoting $ V( s) = \ds \int_s^{1} F( \tau ) \, d \tau$, 
it is easy to see that, at least formally, the "energy"
\beq
\label{1.2}
E(\Phi ) = \int_{\R^N} | \nabla \Phi |^2 \, dx
 + \int_{\R^N} V( |\Phi |^2) \, dx
\eeq
is a conserved quantity.

A second conserved quantity for (\ref{1.1}) is the momentum $P(\Phi ) = ( P_1( \Phi), \dots, P_N(\Phi))$, 
 which describes the evolution of the center of mass of $\Phi$. 
Assuming that $ \Phi \lra 1$ at infinity in a suitable sense 
and denoting by $ \langle \cdot, \cdot \rangle$ the scalar product in $ \C$, the momentum is formally given by 
\beq
\label{mom}
P_k(\Phi) = \int_{\R^N} \langle i \frac{ \p \Phi}{\p x_k}, \Phi -1 \rangle \, dx. 
\eeq

{\bf Traveling waves and the Roberts programme. } 
In a series of papers (see, e.g., \cite{barashenkov2},  \cite{barashenkov1}, 
\cite{GR}, \cite{IS}, \cite{JR}, \cite{JPR}), 
particular attention has been paid to a special class of solutions of (\ref{1.1}), 
namely the traveling  waves. 
These are solutions of the form $ \Phi (x, t) = \psi (x - cty)$, where 
$ y \in S^{N-1}$ is the direction of propagation and $ c >0$ is 
the speed of the traveling wave. 
We say that $\psi $ has finite energy 
if $ \nabla \psi \in L^2(\R^N)$ and $ V(|\psi |^2) \in L^1(  \R^N)$. 
These solutions are supposed to play an important role in the dynamics of (\ref{1.1}).
In view of formal computations and numerical experiments, a list of conjectures,
often referred to as {\it the Roberts programme,} 
has been formulated about the existence, the qualitative properties,   the stability 
of traveling waves and more generally their role in the dynamics of (\ref{1.1}). 

Let $ \psi $ be a finite energy traveling-wave of (\ref{1.1}) moving with  speed $c$. 
Without loss of generality we may assume that $ y = (1, 0, \dots, 0)$.
If $ N \geq 3$, it follows that $ \psi - z_0 \in L^{2^*} ( \R^N)$ 
for some constant $ z_0 \in \C$, where $ 2^* = \frac{ 2N}{N-2}$
(see, e.g., Lemma 7 and Remark 4.2 pp. 774-775 in \cite{PG}).
Since $|\psi | \lra 1$ as $ |x| \lra \infty$,  necessarily  $|z_0| = 1$. 
If $\Phi $ is a solution of (\ref{1.1}) and
$\al \in \R$, then $e^{i \al } \Phi$ is also a solution; hence 
 we may assume that $ z_0 = 1$, thus $ \psi - 1 \in L^{2^*} ( \R^N)$.
Let $ u =  \psi - 1$.
We say that $u$ has finite energy if $\psi $ does so. 
Then $u$ satisfies the equation 
\beq
\label{1.3}
- i c \frac{\p u}{\p x_1} + \Delta u + F(| 1+u|^2) (1+u) = 0 
\qquad \mbox{ in  } \R^N.
\eeq
It is obvious that a function $ u $ satisfies (\ref{1.3}) for some velocity $c$ if and only if 
$u( - x_1, x')$ satisfies (\ref{1.3}) with $c$ replaced by $-c$.
Hence it suffices to consider the case $ c > 0$. 
This assumption will be made throughout the paper.

\begin{figure}[H]
\label{curvesJPR}
\includegraphics[width=7cm]{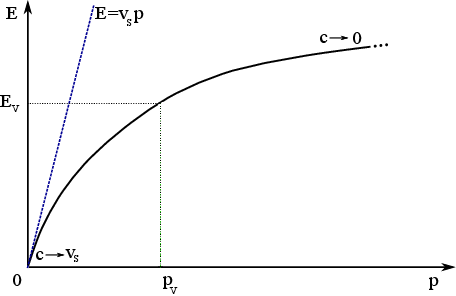}
\qquad \quad
\includegraphics[width=7cm]{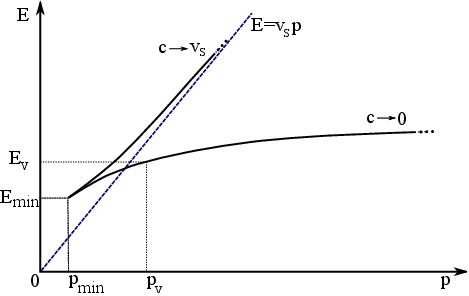}

\caption{{\small Energy ($E$) - momentum ($p$)  diagrams for (GP): (a) in dimension $2$; (b) in dimension $3$.}}
\end{figure}

For the Gross-Pitaevskii equation, C. A. Jones, C. J. Putterman and P. H. Roberts 
computed the energy and the momentum of the traveling waves they had found numerically. 
In space dimension two and three, they obtained the  curves given in figure \ref{curvesJPR}.


Formally, traveling waves are critical points of the energy $E$ 
when the momentum (with respect to the direction  of propagation  $Ox_1$) is fixed, say $ P_1 = p$. 
Equation (\ref{1.3}) is precisely the Euler-Lagrange equation associated to this variational problem, 
and the speed $c$ is the Lagrange multiplier. Note also that, formally, $ c = \frac{ \p E}{\p p }$. 

The first  conjecture in the Roberts programme  asserts that finite energy traveling waves 
of speed $c$ exist if and only~if~$|c| < v_s$.

In space dimension $N=1$, in many interesting applications 
equation  (\ref{1.3}) can be integrated explicitly and one obtains traveling  waves for 
all subsonic speeds. 
The nonexistence of such solutions for supersonic speeds has also been 
proven under general conditions (cf.  Theorem 5.1,  p. 1099 in \cite{M8}). 

Despite of many attempts, 
a rigorous proof of the existence of traveling waves in higher dimensions
has been a long lasting  problem. 
In the particular case of the Gross-Pitaevskii (GP) equation, 
this problem was considered in a series of papers. 
In space dimension $ N =2$, the existence of traveling waves was proven 
in \cite{BS} for all speeds in some interval $(0, \e)$, where $ \e$ is small. 
In space dimension $ N \geq 3$, the existence was proven in \cite{BOS} 
for a sequence of speeds $ c_n \lra 0$ by using constrained minimization; 
a similar result was established in \cite{chiron} for all sufficiently 
small speeds by using a mountain-pass argument. 
In  \cite{BGS}, the existence of traveling waves for (GP) has been proven
in space dimension $ N =2$ and $ N =3 $ by minimizing the energy at fixed momentum.   
The  propagation speed is then  the Lagrange multiplier  associated to minimizers. 
If $ N =2$, this gives solutions for any speed in a set $ A \subset (0, v_s)$, where 
$A$ contains points arbitrarily close to $0$ and to $v_s$ (although it is not clear that 
$A = (0, v_s)$). 
It was shown later in \cite{CM1} that the minimization of the energy at fixed momentum can be used 
in any dimension $N \geq 2$ 
for general nonlinearities such that the nonlinear potential $V$ appearing in the energy is nonnegative. 
Moreover, the set of solutions that it gives   is orbitally stable. 
However, this method has two disadvantages. 
Firstly, it is not clear that the set of speeds of traveling waves constructed in this way form an interval. 
Secondly, it has been proved  in \cite{BGS} and \cite{deLaire} that in space dimension $N \geq 3$ there exists $ v_0 \in (0, v_s)$  
such that minimizing the energy at fixed momentum cannot give traveling waves of speed $ c \in (v_0, v_s)$. 
In particular, the "upper branch" in figure \ref{curvesJPR} (b) cannot be obtained in this way.

In the case of  cubic-quintic type nonlinearities, it was proved in \cite{M3}
that traveling  waves exist for any sufficiently small speed if $ N \geq 4$. 

To our knowledge, even for  specific nonlinearities there are no existence results 
in the literature that cover the whole range $(0, v_s)$ of possible speeds. 

The nonexistence of traveling waves for supersonic speeds ($c > v_s$)
was proven in \cite{gravejat1}
in the case of  the Gross-Pitaevskii equation, respectively in \cite{M8} for a large class of nonlinearities.

\medskip

It is the aim of this paper  to prove the existence of nontrivial finite energy 
traveling waves of (\ref{1.1}) in space dimension $ N \geq 3$, under general 
conditions on the nonlinearity $F$ and for any speed $ c \in (0, v_s)$.  

\medskip

The qualitative properties of traveling waves  have been extensively investigated. 
It turns out that these solutions have the best regularity allowed by the nonlinearity $F$
(see, for instance, \cite{farina1}, \cite{farina2}, \cite{M8}). 
It was proved in \cite{BGS} that the traveling waves to the (GP) equation are analytic functions. 
In view of formal computations, Jones, Putterman and Roberts (\cite{JPR}) predicted the asymptotic behavior 
of traveling waves as $|x| \lra \infty$. 
For the (GP) equation the asymptotics have been computed by P. Gravejat (see \cite{gravejat2} and references therein). 
It is likely that the proofs of Gravejat can be adapted to general nonlinearities. 

Even for specific nonlinearities (such as (GP)), the vortex  structure of traveling waves 
is not yet completely understood.
It was conjectured in \cite{JR, JPR} that there is a critical speed $c_v$
(correesponding to the energy $E_v$ and momentum $p_v$) 
such that traveling waves of speed less than $ c_v$ present vortices,
while those of speed greater than $c_v$ do not. 
The small velocity solutions to (GP) constructed in \cite{BS, BOS, chiron} have vortices. 
For other nonlinearities the behavior may be different. 
For instance,  small speed traveling waves constructed in \cite{M3} in the case of nonlinearities of cubic-quintic  type
do not have vortices. We suspect that traveling waves constructed in the present paper 
develop vortices in the limit $ c \lra 0$ if and only if (\ref{1.1}) does not admit finite energy stationary solutions. 
For general nonlinearities it was recently proved in \cite{CM2}  that traveling waves do not have vortices  
if $c$ is close to $ v_s $ and  $ N \in \{ 2, 3 \}$.

\medskip

The energy-momentum diagrams for (GP) suggest that there are traveling waves of arbitrarily small energy and momentum 
in dimension two. 
Such solutions were obtained in \cite{BGS} 
by minimizing the energy at fixed (and small) momentum; their velocities are close to $ v_s$.
A similar result holds for general  nonlinearities, see \cite{CM1}.
A scattering theory for small energy solutions to (\ref{1.1}) in dimension two is therefore excluded. 

The situation is completely different in higher dimensions. 
It was noticed in \cite{JPR} that the energy and the momentum of the three dimensional traveling waves for (GP)
are bounded from below by positive constants $E_{min}$ and $ p_{min}$, respectively. 
It was proved in \cite{BGS} that the three dimensional (GP) equation does not admit small energy traveling waves, 
and the proof was later extended to higher dimensions in \cite{deLaire}. 
It turns out that this result is true for  general nonlinearities: for any $ N \geq 3$ there is $ k_N > 0 $ such that 
any traveling wave $U$ to (\ref{1.1}) satisfying $ \| \nabla U \|  _{L^2(\R^N)}  < k_N$ must be constant
(see \cite{CM1}, Proposition 1.4). 
This result can be  further improved in dimension $ N \geq 6 $ (see \cite{CM2}, Proposition 18). 
Moreover, S. Gustafson, K. Nakanishi and T.-P. Tsai \cite{GNT1, GNT2} established a scattering theory of small solutions to (GP)
in dimension $N\geq 4$ (in the energy space) and $N=3$ (in some weighted space). 

\medskip

In view of formal computations, Jones, Putterman and Roberts conjectured that after a suitable rescaling, 
the modulus and the phase of 
traveling waves converge in the transonic limit $ c \lra v_s$ 
to the solitary waves of the Kadomtsev-Petviashvili I (KP-I) equation.
The present paper is the first to provide finite energy traveling waves to (\ref{1.1}) 
of speed close to $ v_s$ in dimension $ N \geq 3$. 
Very recently it was  proved  that for general nonlinearities, 
the three-dimensional traveling waves found here have modulus close to 1 (thus can be lifted) and their 
phase and modulus tend, after rescaling, to ground states of the KP-I equation (see  \cite{CM2}, Theorem 6). 
Precise estimates on their energy and momentum have also been established in \cite{CM2} and are in full agreement with those in \cite{JR} and \cite{JPR}.
Hence the conjecture concerning the existence of the "upper branch"  of travelling waves has been proven in dimension three  
(with one exception: it is not clear  that we have a continuum of solutions). 
Quite unexpectedly, a similar asymptotic  behavior of traveling waves in the transonic limit 
cannot be true in dimension $ N \geq 4$  (cf. \cite{CM2}, Proposition 19).

\medskip

A much more difficult problem is to understand the stability of traveling waves and, more generally, their role in the 
dynamics of (\ref{1.1}). 
The guess formulated in \cite{JPR} is that the two-dimensional traveling waves represented in the momentum-energy diagram 
in figure 1  should be stable.
In the three dimensional case solutions on the "lower branch" should be stable, while those on the "upper branch" should be unstable. 
It is also suggested in \cite{JPR} p. 3000 that a solution of (GP) starting in a neighborhood of the upper branch could eventually 
"collapse" onto the lower branch, generating "sound waves that radiate the excess energy... to infinity."

Before even speaking of stability, one has to understand the well-posedness of the Cauchy problem associated to (\ref{1.1}).
Important progress has been achieved  in this direction  during the last years; 
we refer to the survey paper \cite{PG2} (see also \cite{gallo}).
It was proved in \cite{PG}, \cite{PG2} that in the subcritical case $N \in \{ 1, 2, 3 \}$, the Cauchy problem for (GP) 
is globally well-posed for all initial data in the energy space, and that in the critical case $N=4$ 
it is globally well-posed for initial data with small energy. 
The method in \cite{PG, PG2}  adapts to other nonlinearities, including the cubic-quintic case
(note that  the cubic-quintic NLS becomes critical in dimension three).
Global well-posedness of the four-dimensional (GP) and of the three-dimensional cubic-quintic NLS
for all initial data in the energy space has been recently proven in \cite{KOPV}.

In dimension one, traveling waves to (GP) are known explicitly. 
Their orbital stability has been studied and proven in \cite{lin}, 
\cite{BGSS},  \cite{GZ}. Other nonlinearities are also considered in \cite{lin}.
The asymptotic stability of these solutions is not known. 

If the nonlinear potential $V$ is nonnegative, traveling waves can be obtained by minimizing the energy at fixed momentum. 
Moreover, all minimizing sequences are precompact, 
and this gives the orbital stability of the set of solutions constructed in this way  (cf. \cite{CM1}, Theorem 6.2). 
For the (GP) equation, the results in \cite{CM1} imply the orbital stability of the full branch of traveling waves 
in dimension 2, and of traveling waves situated on the "lower branch" below the line $E = v_s p$ in dimension 3. 
If $V$ changes sign, a local minimization of the energy at fixed momentum is still possible in dimension 2 and gives
a branch of orbitally stable traveling waves. 
To our knowledge, the orbital stability/instability  of solutions  corresponding to the "upper branch" in dimension 3 
as well as the asymptotic stability of traveling waves in any dimension $ N \geq 2$ are still open problems.

\medskip

{\bf Main results. }  We will consider the following set of assumptions: 

\medskip

{\bf A1. } The function $F$ is continuous on $[0, \infty)$, 
$C^1$ in a neighborhood of $ 1$, $F(1) = 0$ and $F'(1)=-1$.

\smallskip

{\bf A2. } 
There exist $C > 0$ and $ p _0 < \frac{2}{N-2} $ such that $|F(s) | \leq C(1 + s^{p_0}) $
for any $ s \geq 0$. 

\smallskip

{\bf A3. } There exist $C, \,  \al _0> 0$ and $r_* > 1 $  such that $F(s) \leq - C s^{\al _0} $ 
for any $ s \geq r_*$. 

\medskip 

\medskip

Our main  results can be summarized as follows.

\begin{Theorem}
\label{T1.1}
Assume that $ N \geq 3$, $0< c< v_s$,    and  the conditions (A1) and (A2)   are satisfied. 
Then equation (\ref{1.3}) admits a nontrivial finite energy solution $ u  $.
Moreover, $ u \in W_{loc}^{2,p}  ( \R^N)$ for any $ p \in [1, \infty ) $
and, after a translation, $u$ is axially symmetric with respect to $Ox_1$. 
\end{Theorem}

\begin{Corollary}
\label{C1.2}
Suppose that $ N \geq 3$, $0< c< v_s$,    and  the conditions (A1) and (A3)   are verified. 
Then equation (\ref{1.3}) admits a nontrivial finite energy solution $ u$   such that  
 $ u \in W_{loc}^{2,p}  ( \R^N)$ for any $ p \in [1, \infty ) $
and, after a translation, $u$ is axially symmetric with respect to $Ox_1$. 
\end{Corollary}

It is easy to see how Corollary \ref{C1.2} follows from  Theorem \ref{T1.1}.
Indeed, suppose that Theorem \ref{T1.1} holds.
Assume that (A1) and (A3) are satisfied. 
Let  $C,\;  r_*, \;  \al _0$ be as in (A3).
There exist $ \beta \in (0, \frac{2}{N-1}) $,  $ \tilde{r} > r_* $, and $C_1 > 0$ such that 
$$
Cs^{ 2 \al _0} - \frac{1}{2} \geq C_1 ( s - \tilde{r}) ^{2 \beta } 
\qquad \mbox{ for any } s \geq \tilde{r}. 
$$
Let $\tilde{F}$ be a function with the following properties: 
$ F = \tilde{F}$ on $[0, 4 \tilde{r}^2]$, 
$\tilde{F} ( s) = - C_2 s^{ \beta }$ for $s$ sufficiently large, 
and $ \tilde{F} ( s ^2 ) + \frac{1}{2} \leq - C_3 ( s - \tilde{r}) ^{ 2 \beta} $ 
for any $ s \geq \tilde{r}$,
where $ C_2$, $C_3$ are some positive constants. 
Then $ \tilde{F}  $ satisfies (A1), (A2), (A3) and  Theorem \ref{T1.1} implies 
that equation  (\ref{1.3}) with $ \tilde{F} $ instead of $F$ 
has a nontrivial finite energy solution $ u $. 
From  the proof of Proposition 2.2 (i) p. 1079-1080  in \cite{M8} it follows  that any such solution satisfies 
$|1+u |^2 \leq 2 \tilde{r} ^2$, and consequently 
$F(|1+u|^2) =  \tilde{F} (|1+ u|^2)$. 
Thus $u$ satisfies (\ref{1.3}). 

We have to mention that the traveling waves in Theorem \ref{T1.1} are obtained 
as minimizers of the functional $E + c P_1$ under a Pohozaev constraint (see below), 
where $E$ is the energy and $P_1$ is the momentum with respect to $x_1$. 
Of course, if (A1) and (A3) are satisfied but (A2) does not hold,
 we do not claim that the solutions given by Corollary \ref{C1.2} 
 still solve  the same minimization problem. 
In fact,  assumptions (A1) and (A3) alone do not imply that $E$ is 
 well-defined on a convenient  function space.

In particular, for $F(s) = 1-s $ the conditions (A1) and (A3) are satisfied and
it follows that the Gross-Pitaevskii equation admits   nontrivial 
traveling waves of finite energy in any space dimension $ N \geq 3$ and 
for any speed $ c \in (0,  v_s )$
(although (A2) is not true  for $N > 3$: the (GP) equation is critical if $ N =4$,  
and supercritical if $ N \geq 5$). 
A similar result holds for the cubic-quintic  NLS.

\medskip

{\bf Notation and function spaces. } 
Throughout the paper,  $ \Lo ^N$ is the Lebesgue measure on $ \R^N$ and  
 $ \omega_N = \Lo ^N( B(0,1)) $ is the Lebesgue measure of the unit ball.
For $ x = ( x_1, \dots, x_N) \in \R^N$, 
we denote  $ x' = ( x_2, \dots , x_N ) \in \R^{N-1}$. 
We write $\langle z_1, z_2 \rangle $ for the scalar product 
of two complex numbers $z_1, z_2$. 
Given a  function $f$ defined on $ \R^N$ and $ \la , \, \si  > 0$, 
we denote~by 
\beq
\label{1.7}
f_{\la, \si }= f\left( \frac{x_1}{\la }, \frac{x'}{\si} \right)
\eeq
the dilations of $f$.
The behavior of functions and of functionals with respect to dilations in 
$ \R^N$ will be very important.
For $ 1 \leq p <N$, we denote by $ p^*$ the Sobolev exponent 
associated to $p$, that is $ \frac{1}{p^*} = \frac 1p - \frac 1N$.

If  (A1) is satisfied,   let  $ V(s) = \ds \int_s ^{1} F(\tau ) \, d \tau$.
Then the sound velocity at infinity associated to (\ref{1.1}) 
is $ v_s = \sqrt{2} $ and using Taylor's formula for $s$ in a neighborhood 
of $1 $ we have
\beq
\label{1.4}
V(s) = \frac 12 V''(1) ( s - 1)^2 + ( s - 1)^2 \e (s - 1)
= \frac 12  ( s - 1)^2 + ( s - 1)^2 \e (s - 1) , 
\eeq
 where $ \e (t) \lra 0  \mbox{ as } t \lra 0.$ 
Hence  $V(|\psi |^2) $ can be approximated by 
$\frac 12  ( |\psi |^2 - 1)^2$  for $|\psi |$ close to $ 1$. 

We fix an odd function $ \ph \in C^{\infty } (\R)$ such that $ \ph (s) = s $ 
for $ s \in [0, 2  ]$, $ 0 \leq \ph ' \leq 1 $ on $ \R$ and 
$ \ph (s) = 3  $ for $ s \geq 4  $. 
If assumptions (A1) and (A2) are satisfied, 
it is not hard to see that 
there exists $C_1 > 0$ 
such that
\beq
\label{i1}
\begin{array}{l}
|V(s)| \leq C_1 ( s - 1)^2 \quad \mbox{ for any } s \leq 9 ; 
\\
\mbox{in particular, } 
| V(\ph ^2(\tau )) | \leq C_1 (\ph ^2(\tau ) - 1) ^2 \mbox{ for any  } \tau. 
\end{array}
\eeq
Given  $ u \in H_{loc}^1 ( \R^N) $ and  an open set  $\Om \subset  \R^N$, 
the modified  Ginzburg-Landau  energy of $u$ in $ \Om $ is defined by 
\beq
\label{1.5}
E_{GL }^{\Om } (u) = \ds \int_{\Om } |\nabla u |^2 \, dx 
+ \frac 12 \int_{\Om } \left( \ph ^2(|1+ u|) - 1 \right)^2 \, dx .
\eeq
We  simply write $ E_{GL}(u) $ instead of $E_{GL}^{\R^N} (u)$.
The modified Ginzburg-Landau energy will  play a central role in our analysis. 
We consider the function space
\beq
\label{X} 
\begin{array}{rcl}
\Xo &  = & \{ u \in \DR \; | \; \ph ^2(|1+ u|) - 1 \in L^2( \R^N) \} 
\\
 & = & \{ u \in \dot{H}^1(\R^N) \; | \; 
u \in L^{2^*} (\R^N), \; E_{GL}(u) < \infty \},  
\end{array}
\eeq
where $\DR$ is the completion of $ C_c^{\infty } (\R^N)  $ 
for the norm $\|v \| = \|\nabla v \| _{L^2}$. 

Since $  \ph ^2(|1+ u|) - 1 = (  \ph (|1+ u|) + 1)( \ph (|1+ u|) - 1) $ and 
$ 1 \leq  \ph (|1+ u|) + 1 \leq 4$, 
it is obvious that $ \ph ^2(|1+ u|) - 1  \in L^2( \R^N)$ if and only if $\ph (|1+ u|) - 1  \in L^2( \R^N)$. Let $ N \geq 3$. 
We claim that for $ u \in \DR$ there holds $ \ph (|1+ u|) - 1  \in L^2( \R^N)$ if and only if $|1+ u| - 1  \in L^2( \R^N)$.
Indeed, if $ |1 + u | \leq 2 $ then $ \ph(|1 + u |) = |1 + u|$. 
If $ |1 + u | > 2 $ then necessarily  $|u| > 1 $ and 
$$
0 \leq |1 + u | - \ph (|1 + u |) <  |1 + u | <  2 |u| < 2 |u|^{\frac{2^*}{2}}.
$$
For $ N \geq 3 $ we have $ \DR \subset L^{2^*}(\R^N)$ by the Sobolev embedding, hence $|u|^{\frac{2^*}{2}} \in L^2( \R^N)$ 
and the claim follows. We have thus proved that
$$
\Xo   =  \{ u \in \DR \; \big| \; |1+ u| - 1 \in L^2( \R^N) \} .
$$

If $ N \geq 3$ and (A1), (A2) are satisfied, it is not hard to see that the function $\psi = 1 + u$ 
satisfies $ \nabla \psi \in L^2( \R^N)$ and $V(|\psi|^2) \in L^1( \R^N)$
if and only if $ u \in \Xo$ (see Lemma \ref{L4.1} below). 
Note that for $ N =3$, $ \Xo $ is {\it not} a vector space. 
However, in any space dimension we have $ H^1(\R^N) \subset \Xo $. 
If $ u \in \Xo$, it is easy to see that for any $ w \in H^1(\R^N) $ with 
compact support we have $ u + w \in \Xo $. 
For $ N = 3, 4$ it can be proved that $ u \in \DR $ 
belongs to $ \Xo $ if and only if $ |1 + u|^2 - 1 \in L^2( \R^N)$, 
and  consequently $\Xo $ coincides with the space $F_{1}$ introduced by P. G\'erard 
in \cite{PG}, section 4.

\medskip

{\bf Some ideas in the proofs and outline of the paper.}

At least formally, the solutions of (\ref{1.3}) are critical points of the functional 
\beq
\label{Ec}
E_c(u) = E(u) + c Q(u) = \int_{\R^N} |\nabla u |^2 \, dx + c Q(u) + \int_{\R^N} V(|1+u |^2) \, dx,
\eeq
where $Q=P_1$ is the momentum with respect to the $x_1-$direction.

It is the aim of section 2 to give a convenient definition of the momentum on the whole space $ \Xo$ 
and to study its basic properties. For now,    the formal definition (\ref{mom}) is sufficient.

The existence of finite energy traveling waves has been conjectured for all subsonic speeds and the
nonexistence of such solutions is known for all spersonic speeds 
(at least under some additional technical assumptions, see \cite{M8}), 
thus it is important to understand what changes in the structure of $E_c$ as $c $ crosses the sound velocity. 

If $ 0 < c < \sqrt{2}$, we may choose $ \e, \, \de > 0$ such that $ c < \sqrt{2}(1 - 2 \e) (1 - \de)$. 
Assume that $ u \in \Xo $ satisfies $ 1 - \de \leq |1 + u| \leq 1 + \de$. 
Then there is a lifting $ 1 + u = \rho e^{i \theta}$  and a simple computation shows that 
$$
\begin{array}{l}
|\nabla u |^2 = |\nabla \rho |^2 + \rho ^2 |\nabla \theta |^2 , \qquad
Q(u) = - \ds \int_{\R^N} ( \rho ^2 - 1) \frac{ \p \theta}{\p x_1} \, dx \qquad \mbox{ and } 
\\
\\
V(|1 + u |^2) = V(\rho^2) = \frac 12 (\rho^2 -1) ^2 + o ((\rho^2 -1) ^2 ) \geq \frac{1-\e}{2} (\rho^2 -1) ^2
\end{array}
$$
provided that $ \de $ is sufficiently small.
Then we have 
\beq
\label{vanish0}
\begin{array}{l}
\ds |c Q(u)| \leq \sqrt{2} (1 - 2 \e) (1 - \de) \int_{\R^N} | \rho ^2 -1| \cdot \Big| \frac{ \p \theta}{\p x_1} \Big| \, dx 
\\
\\
\ds \leq (1 - 2 \e) \int_{\R^N} ( 1 - \de )^2 \Big| \frac{ \p \theta}{\p x_1} \Big| ^2 + \frac 12 ( \rho ^2 -1)^2 \, dx 
\\
\\
\ds \leq \int_{\R^N} ( 1 - 2 \e) \rho ^2 | \nabla \theta |^2 + V( \rho ^2) - \frac{\e}{2} ( \rho ^2 -1)^2  \, dx 
\leq E(u) - \e E_{GL}(u).
\end{array}
\eeq
Thus $E_c (u) \geq \e E_{GL}(u) $ if $ |1 + u|$ is sufficiently close to $1$ in the $L^{\infty}$ norm. 
Since $E_{GL}(u)$ measures, in some sense, the closeness of $ 1 + u $ to $1$, 
we would like to establish a similar estimate for all functions with  small Ginzburg-Landau energy.  
However, $E_{GL}(u)$ does not control $\| \, | 1 + u | - 1 \|_{L^{\infty}}$. Moreover, there are functions with 
arbitrarily small Ginzburg-Landau energy which present small-scale topological "defects" (e.g., dipoles). 
To get rid of these difficulties  we use a procedure of regularization by minimization, 
which is introduced and studied  in section 3. 
Given $ u \in \Xo$, we minimize the functional 
$v \longmapsto E_{GL}(v) + \ds \frac{1}{h^2} \int_{\R^N} \ph(|v-u| ^2) \, dx $ in the set 
$\{ v \in \Xo \; | \;  v - u \in H^1( \R^N) \}.$
It is shown that minimizers exist (but are perhaps not unique) and any minimizer $ v_h $ has remarkable properties, for instance:

$\bullet $ $E_{GL}(v_h) \leq E_{GL}(u)$, 

$\bullet $ $\| v_h - u \|_{L^2} \lra 0 $ as $ h \lra 0$, and 

$\bullet $  $ \| \, |1 + v_h | - 1 \|_{L^{\infty}} $ can be estimated in terms of $h$ and $E_{GL}(u)$ 
and  is arbitrarily small if $E_{GL}(u)$  is sufficiently small. 


In section 4 we describe the variational framework. 
Using the above regularization procedure we prove that for any $ \e \in (0, 1 - \frac{c}{v _s} )$ 
and for all $ u \in \Xo $ with $E_{GL}(u)$ sufficiently small there holds $E_c(u) \geq \e E_{GL}(u)$. 
Then we show that for all $ k > 0$, the functional $E_c$ is bounded on the set 
$\{ u \in \Xo \; | \; E_{GL}(u) \leq k \}$. 
Let
$$
E_{c, min}(k) =  \inf \{ E_c (u) \; | \; u \in \Xo, \; E_{GL}(u) = k \}.
$$
We prove that for $0 < c < v_s$ the function $E_{c, min} $ has the following properties:

i) For any $ \e \in (0,1 - \frac{c}{v _s} )$  
there is  $ k_{\e} > 0$ such that $ E_{c, min} (k) > \e k $ for $ k \in (0, k_{\e})$.

ii)  $ \ds \lim_{ k \ra \infty } E_{c, min} (k ) = - \infty$. 

iii) For any  $ k > 0$ we have $ E_{c, min} (k) < k$. 
\\
The situation is very different if $ c > v_s$: in that case it can be proved that $E_{c, min}$ is negative and decreasing on $(0, \infty)$.



In order to get critical points of  $E_c$, it is tempting to  minimize  $E_c(u)$ under the constraint 
$E_{GL}(u) = k$ or $Q(u) = k$, where $k$ is a constant, and then to search for those  $k$ which give minimizers 
with the associated Lagrange multiplier equal to zero. 
However, it is well-known that it is hard to control the Lagrange multipliers in minimization problems which
do not have appropriate scaling or  homogeneity properties.
In order to avoid that difficulty we adopt the following strategy. 
We introduce the  functionals:
\beq
\label{A}
A(u) = \ds \int_{\R^N} \sum_{ j =2}^N 
\Big\vert \frac{ \p u }{\p x _j } \Big\vert ^2 \, dx, 
\eeq
\beq
\label{Bc}
B_c (u) = \ds \int_{\R^N} 
\Big\vert \frac{ \p u }{\p x _1 } \Big\vert ^2 \, dx
+ c Q(u) 
+ \ds \int_{\R^N}  V(|1+u |^2 ) \, dx  \qquad  \mbox{ and } 
\eeq
\beq
\label{Pc}
P_c(u) = \frac{N-3}{N-1} A(u) + B_c (u) .
\eeq
It is obvious that $E_c(u) = A(u) + B_c (u) = \frac{2}{N-1} A(u) + P_c (u)$. 
If the assumptions (A1) and  (A2)  above are satisfied, 
it can  be proved (see Proposition 4.1 p. 1091-1092 in \cite{M8})
 that any  traveling wave $u \in \Xo $  of (\ref{1.1}) 
must satisfy the Pohozaev-type identity $P_c(u) = 0$.  Indeed, 
it is easy to see that for any $ u \in \Xo $ we have $ E_c( u_{1, \si}) = \si^{N-3} A(u) + \si ^{N-1} B_c (u)$. 
Formally,   a critical point $u$ of $E_c$ should satisfy
$\frac{d}{d \si}_{| \si = 1} (E_c(u_{1, \si })) =0$, which gives  precisely $P_c(u) = 0$. 
We will prove the existence of traveling waves by showing that the problem of minimizing 
$E_c$ in the set 
$$
 \Co = \{ u \in \Xo \; | \; u \neq 0, P_c (u ) = 0 \} 
$$
admits solutions. 
It turns out that 
 minimizing a functional under a Pohozaev constraint (almost) automatically generates critical points of that functional; 
 that is, the Lagrange multiplier is fixed. 
This is a very general observation which seems to work in many problems in Calculus of Variations. 
To  our knowledge, it is used here for the first time. 
Let us explain how it works for $E_c$ in dimension $ N \geq 4$.  
Assume that $ u \in \Co $  satisfies the Euler-Lagrange equation 
$ E_c '(u) = \al P_c'(u)$. 
Then $u$ is a critical point of the functional $ E_c - \al P_c$. 
Formally we have 
$\frac{d}{d \si}_{| \si = 1} (E_c(u_{1, \si })- \al P_c ( u_{1, \si })) =0$, which is equivalent to 
$$
P_c (u) - \al \left[ \left(\frac{N-3}{N-1} \right) ^2 A(u) + B_c (u) \right] = 0. 
$$
Since $P_c(u) = 0$, the above identity implies $ \al \frac{N-3}{N-1} \cdot \left( \frac{N-3}{N-1} -1 \right) A(u) =0$, 
thus either $A(u) =0 $ (and $u$ is constant), or $ \al = 0$. 

The next step is to prove that $ \Co $ is not empty and $\ds \inf \{ E_c (u) \; | \; u \in \Co \} > 0$. 
Let us  present here the arguments in dimension $N \geq 4$. If $ u \in \Co $ we have $B_c(u) = - \frac{N-3}{N-1} A(u) < 0$. 
Then it is easy to see that the function 
$ \si \longmapsto E_c( u_{1, \si}) = \si^{N-3} A(u) + \si ^{N-1} B_c (u) $ 
is increasing on $(0,1)$ and decreasing on $(1, \infty)$, thus achieves a maximum at $ \si = 1$. Hence
$$
E_c(u) = E_c(u_{1,1}) \geq E_c(u_{1, \si }) \geq E_{c, min} (E_{GL}(u_{1, \si })) \qquad \mbox{ for all } \si > 0.
$$
Since $ \si \longmapsto E_{GL}(u_{1, \si }) $ takes all values in $(0, \infty)$ we infer that 
$$
T_c :=\ds \inf \{ E_c (u) \; | \; u \in \Co \} \geq \sup \{ E_{c, min}(k) \; | \; k > 0 \} > 0. 
$$

In section 5 we consider the case $ N \geq 4 $  and we prove that the functional 
$E_c$ has minimizers 
in $ \Co $ and these minimizers are solutions of (\ref{1.3}). 
To show the existence of minimizers we use the concentration-compactness principle 
and the regularization procedure developed in section 3.
The most difficult part is to show that minimizing sequences do not "vanish,"
that is their Ginzburg-Landau energy does not spread over $ \R^N$. 
Assume that $ N \geq 4 $ and $(u_n)_{n \geq 1} $ is a minimizing sequence for $E_c $ on $ \Co $ that "vanishes."
Letting $ \si _0 = \sqrt{\frac{2(N-1)}{N-3}} $ and $ \tilde{u}_n = (u_n)_{1, \si _0}$, 
we see that $(\tilde{u}_n)_{n \geq 1}$  also vanishes and 
$A(\tilde{u}_n ) + E_c (\tilde{u}_n) = \si _0 ^{N-1} P_c (u_n) = 0.$
Since $ A(\tilde{u}_n) = \si _0 ^{N-3}  A(u_n)$ and
$A(u_n) = \frac{N-1}{2} (E_c(u_n) - P_c(u_n)) \geq  \frac{N-1}{2}  T_c >0$, we get 
\beq
\label{vanish1}
\limsup_{n \ra \infty} E_c (\tilde{u}_n ) < 0.
\eeq

On the other hand, the vanishing of $ (\tilde{u}_n)_{n \geq 1}$ implies that 
$$
\int_{\R^N} V(|1 + \tilde{u}_n|^2) \, dx = \frac 12 \int_{\R^N} \left( \ph ^2(|1 + \tilde{u}_n |) -1 \right)^2 \, dx + o (1). 
$$
Using the regularization procedure in section 3 (see Lemma \ref{vanishing} there) we construct a sequence $ h_n \lra 0 $
 and for each $n$ we find a minimizer $v_n$ of the functional 
 $  E_{GL}(v) + \ds \frac{1}{h_n ^2} \int_{\R^N} \ph(|v - \tilde{u}_n |^2) \, dx $ 
 such that $\| \, | 1 + v_n | - 1 \|_{L^{\infty}} \lra 0 $ as $ n \lra \infty$. 
Then we have $ Q(  \tilde{u}_n) = Q( v_n) + o(1) $ and 
\beq
\label{vanish2}
E_c ( \tilde{u}_n) = E_{GL} ( \tilde{u}_n ) + cQ( \tilde{u}_n) + o(1) 
\geq E_{GL}(v_n) + c Q( v_n) + o(1) \geq 0 
\eeq
for all  $n$ sufficiently large by (\ref{vanish0}).
It is clear that (\ref{vanish1}) and (\ref{vanish2}) give a contradiction and this rules out vanishing.

If "dichotomy" occurs, the Ginzburg-Landau energy of $u_n$ is located in two regions which are far away from each other as $ n \lra \infty$. 
Using again the regularization procedure we show that there are  functions 
 $ u_{n,1}$,  $u_{n,2}$ such that $(E_{GL}(u_{n,i}))_{n \geq 1}$ is bounded and stays away from zero for $i=1,2$,  and
\beq
\label{dichotomy} 
| A(u_n) - A(u_{n,1}) -  A(u_{n,2}) | \lra 0 \qquad \mbox{ and } \qquad 
| P_c (u_n) - P_c(u_{n,1}) -  P_c(u_{n,2}) | \lra 0 
\eeq
as $ n \lra \infty$.
It is easy to see that  $(P_c (u_{n, i}))_{n \geq 1}$ is bounded for $ i=1,2$. 
Passing again to a subsequence, we may assume that
$P_c (u_{n, 1}) \lra p_1 $ and $ P_c (u_{n, 2}) \lra p_2$, 
 where  $ p_1 + p_2 = 0 . $

If $ p_1 = p_2 = 0$ we show that ${\ds \liminf_{n \ra \infty}} E_c (u_{n,i}) \geq T_c$ for  $ i =1,2$, and then 
$$
\liminf_{n \ra \infty}E_c(u_n) =  \liminf_{n \ra \infty} (E_c(u_{n, 1}) + E_c(u_{n, 2}) ) \geq 2T_c , 
\qquad \mbox{ a contradiction.}
$$

If $ p_1 < 0$,   we use Lemma \ref{L4.8} (ii) which asserts that for any bounded sequence 
$ (v_n) _{n \geq 1} \subset \Xo$ satisfying ${\ds \lim_{n \ra \infty}} P_c(v_n) <0$, there holds 
$ {\ds \liminf_{n \ra \infty}} A( v_n) > \frac{N-1}{2} T_c$. Hence
$$
\liminf_{n \ra \infty} E_c( u_n) 
= \frac{2}{N-1} \liminf_{n \ra \infty} A( u_n) 
\geq \frac{2}{N-1} \liminf_{n \ra \infty} A( u_{n, 1} ) > T_c, 
$$
 again a contradiction. A similar argument is valid  if $ p_2 < 0$. 

Since "vanishing" and "dichotomy" are excluded, necessarily "concentration" occurs 
and then 
we show that $(u_n)_{n \geq 1}$ has a subsequence which converges to a minimizer of $E_c$ in $\Co $. 

There are some important differences in the case $ N =3$ with respect to the case $ N \geq 4$, 
most of them due to different scaling properties. 
For instance, for any $ v \in \Xo $ we have $ A( v_{1, \si }) = A(v)$ and $B_c ( v_{1, \si}) = \si ^2 B_c (v)$. 
If $ v \in \Co$, for all $ \si > 0 $ we have 
$$
P_c ( v_{1, \si }) = B_c ( v_{1, \si }) = \si ^2 B_c (v) = 0 \qquad \mbox{ and } \qquad 
E_c ( v_{1, \si }) = A ( v_{1, \si }) = A(v) = E_c (v). 
$$
It is then clear that one may expect convergence of minimizing sequences for $E_c$ in $ \Co $ only after scaling. 
The proofs that vanishing and dichotomy do not occur are also a  bit more involved. 
We treat separately the case $ N=3$ in section 6. 

Next we have to prove that any minimizer  $u$  of $ E_c$ in $ \Co $ is nondegenerate and satisfies an Euler-Lagrange equation 
$ E_c ' (u) = \al P_c' (u)$ (then necessarily $ \al = 0$, as explained above). 
This is done in Proposition \ref{P5.6} in the case $ N \geq 4$, 
respectively in Lemma \ref{L6.4} and Proposition \ref{P6.4} in the case $ N=3$. 

Finally, we prove that traveling  waves found by minimization in sections 5 and 6 are axially symmetric
(as one would expect from physical considerations, see \cite{JR}).

\medskip

In space dimension two the situation is different,  mainly because of different scaling properties.  
Indeed, if $N=2$ it is easy to see that for any nonconstant function $u$ satisfying $P_c(u) = 0$, 
the mapping $ \si \longmapsto E_c ( u_{1, \si}) $ is decreasing on $(0, 1]$ and increasing on $[1, \infty)$, 
hence achieves its minimum at $ \si = 1$. 
This is exactly the opposite of what happens in the case $ N \geq 4$, when $E_c ( u_{1, \si } )$
reaches its maximum at $ \si = 1$. 
 It can be proved that for $ N =2 $ we have  $\ds \inf \{ E_c (u) \; | \; u \in \Xo, u \neq 0, P_c (u) = 0 \} = 0$ 
 and there are no 
minimizers of $E_c$ subject to the constraint $P_c = 0$. 
By using   different approaches,   the existence of two-dimensional 
traveling waves   has recently been proven in \cite{CM1} 
 for a set of speeds that contains elements arbitrarily close to zero and to $ v_s$. 
The existence for all speeds $ c \in (0, v_s)$ is still an open problem.
Although some of the results in sections 2$-$4 are also valid in space dimension $N=2$
(with straightforward modifications in the proofs), for simplicity 
 we assume throughout that $ N \geq 3$. 

\medskip

If $c =0$ and assumptions (A1) and (A2) are satisfied, 
equation (\ref{1.3}) has finite energy solutions if and only if the nonlinear potential $V$ achieves negative values. 
The existence follows, for instance, from Theorem 2.1 p. 100 and Theorem 2.2 p. 103 in \cite{brezis-lieb}
(see also \cite{AdB}). 
On the other hand, 
any  finite energy solution $ \psi $ of the equation $ \Delta \psi + F(|\psi |^2 )\psi = 0 $ in $ \R^N$  satisfies the Pohozaev identity
$$
(N-2) \int_{\R^N} |\nabla \psi |^2 \, dx + N \int_{\R^N} V(|\psi |^2) \, dx = 0 
$$
(see, e.g., Lemma 2.4 p. 104 in \cite{brezis-lieb}), and then it is clear that $ \psi $ must be constant if $V$ is nonnegative. 
In the case $ c =0$, our proofs imply that $E_0$ has a minimizer in the set 
$ \{u \in \Xo \; | \; u \neq 0, P_0 (u) = 0 \}$ 
whenever this set  is not empty. 
Then it is not hard to prove that minimizers satisfy (\ref{1.3}) for $c =0$ 
(after a scale change  if $N=3$). 
However, for simplicity we assume throughout (unless the contrary is explicitly mentioned)
that $0< c< v_s$.

\section{The  momentum} 

A good definition of the momentum 
is essential in any attempt 
to find solutions of (\ref{1.3}) by using a variational approach. 
Roughly speaking, the momentum (with respect to the $x_1-$direction) 
should be a functional with derivative $ 2 i u_{x_1}$. 
Various definitions have been given in the literature 
(see \cite{BS}, \cite{BGS}, \cite{BOS}, \cite{M3}),  any 
of them having its advantages and its inconvenients.  
Unfortunately, none of them is valid for all functions in $ \Xo$. 
We propose a new and more general definition in this section. 

\medskip

It is clear that for functions $ u \in H^1(\R^N)$, the momentum should be given by 
\beq
\label{2.1}
Q_1(u) = \ds \int_{\R^N} \langle iu_{x_1}, u \rangle \, dx, 
\eeq
and this is indeed a nice functional on $ H^1(\R^N)$.
The problem is that there are functions $ u \in \Xo \setminus  H^1(\R^N)$ such that 
 $ \langle iu_{x_1}, u \rangle \not\in L^1( \R^N)$. 

If $ u \in \Xo $ is such that $ 1+ u$ admits a lifting  $ 1+ u = \rho e^{ i \theta}$, 
a formal computation gives
\beq
\label{2.2}
\ds \int_{\R^N} \langle iu_{x_1}, u \rangle \, dx = 
- \int_{\R^N} \rho ^2 \theta _{x _1}\, dx  = - \int_{\R^N} (\rho ^2 - 1)
 \theta _{x _1}\, dx.
 \eeq
 It is not hard to see that if $ u \in \Xo $ is as above, 
 then $(\rho ^2 - 1 ) \theta _{x_1}  \in L^1( \R^N)$. 
However, there are many "interesting" functions  $u \in \Xo $ such that 
$ 1+ u $ does not admit a lifting. 

Our aim is to define the momentum on $ \Xo $ in such a way that it agrees with 
(\ref{2.1}) for functions in $ H^1(\R^N) $ and with (\ref{2.2}) when a lifting 
as above exists. 

\begin{Lemma}\label{lifting} 
Let $ u \in \Xo $ be such that $ m \leq |1+ u (x)| \leq 2  $ 
a.e. on $ \R^N$, where $ m > 0$. 
There exist two real-valued functions $ \rho , \theta$ such that 
$ \rho - 1 \in H^1(\R^N)$, $ \theta \in \DR $, $ 1+ u = \rho e^{i \theta } $ a.e. 
on $ \R^N$ and 
\beq 
\label{2.3} 
\langle iu_{x_1}, u \rangle  =  \frac{ \p }{\p x_1 } ( \I (u) -  \theta) 
- ( \rho ^2 - 1 ) \frac{ \p \theta}{\p x _1} 
\qquad \mbox{ a.e. on } \R^N .
\eeq
Moreover, we have 
$ \ds \int _{\R^N} \Big\vert ( \rho ^2 - 1) \theta _{x_1} \Big\vert \, dx 
\leq \frac{1}{\sqrt{2} m} E_{GL} (u)$. 
\end{Lemma}

{\it Proof. } 
Since $ 1+u \in H_{loc}^1 (\R^N) $, the fact that there exist 
$ \rho, \theta \in H_{loc}^1(\R^N)$ such that $ 1+u = \rho e^{i \theta } $ a.e. 
is standard and follows from Theorem 3 p. 38 in \cite{BBM}. 
We have
\beq
\label{2.4}
\bigg\vert \frac{\p u}{\p x _j } \bigg\vert ^2 = 
\bigg\vert \frac{\p \rho}{\p x _j } \bigg\vert ^2 + 
\rho ^2 \bigg\vert \frac{\p \theta}{\p x _j } \bigg\vert ^2
\qquad \mbox{ a.e. on } \R^N \mbox{ for } j = 1, \dots, N.
\eeq
Since $ \rho = |1+ u| \geq m $ a.e., it follows that 
$ \nabla \rho, \nabla \theta \in L^2( \R^N)$. 
If $ N \geq 3$, we infer that there exist $ \rho _0, \theta _0 \in \R$ such that 
$ \rho - \rho _0 $ and $ \theta - \theta _0$ belong to $L^{2^*} (\R^N)$.
Then it is not hard to see that $ \rho _0 = 1 $ and $ \theta _0 = 2k _0\pi $, 
where $ k _0\in \Z$.  
Replacing $ \theta $ by $ \theta - 2 k _0\pi$, we have $ \rho - 1, \theta \in \DR$. 
Since $ \rho  \leq 2  $ a.e., we have 
$ \rho ^2 - 1 = \ph ( | 1+ u | ) ^2 - 1 \in L^2 ( \R^N) $ because
$ u \in \Xo $. Clearly $ | \rho - 1 | = \frac{ | \rho ^2 - 1 | }{\rho + 1 } 
\leq  | \rho ^2 - 1 | $, hence $ \rho - 1 \in L^2( \R^N)$. 

A straightforward computation gives
$$
\langle iu_{x_1}, u \rangle  = \langle iu_{x_1}, -1 \rangle   - \rho ^2 \theta_{x_1}
=  \frac{ \p }{\p x_1 } ( \I (u) - \theta) 
- ( \rho ^2 - 1 ) \frac{ \p \theta}{\p x _1}.  
$$
By (\ref{2.4}) we have $ \big\vert \frac{\p \theta}{\p x _j } \big\vert
\leq \frac{1}{\rho} \big\vert \frac{\p u}{\p x _j } \big\vert
\leq \frac 1m \big\vert \frac{\p u}{\p x _j } \big\vert$ 
and  the Cauchy-Schwarz inequality gives
$$
\ds \int _{\R^N} \Big\vert ( \rho ^2 - 1) \theta _{x_1} \Big\vert \, dx 
\leq \| \rho ^2 - 1 \|_{L^2} \|\theta _{x_1} \|_{L^2} 
\leq \frac 1m \|\rho ^2 - 1 \|_{L^2} \| u _{x_1}\| _{L^2} 
\leq \frac{1}{m\sqrt{2}} E_{GL}(u).
$$
\hfill $\Box $

\begin{Lemma} 
\label{L2.2}
Let $ \chi \in C_c^{\infty } (\C, \R)$ be a function such that 
$\chi = 1 $ on $ B(0, \frac{1}{4})$, $0 \leq \chi \leq 1 $ and
$\mbox{\rm supp}(\chi) \subset B(0, \frac{1}{2})$.
For an arbitrary $ u \in \Xo$, denote $ u_1 = \chi (u) u$ and $ u_2 = ( 1 - \chi (u)) u$. 
Then $ u_1 \in \Xo$, $ u_2 \in H^1 ( \R^N)$ and the following estimates hold: 
\beq
\label{2.5}
| \nabla u_i | \leq C |\nabla u | \quad \mbox{ a.e. on } \R^N \mbox{ for }  i=1, 2, \mbox{ where } 
C \mbox{ depends only on } \chi, 
\eeq
\vspace*{-10pt}
\beq
\label{2.6}
\| u_2 \|_{L^2(\R^N)} \leq C _1 \| \nabla u \|_{L^2(\R^N)} ^{\frac{2^*}{2}} 
\; \mbox{ and } \;
\|(1 - \chi ^2 (u)) u \|_{L^2(\R^N)} \leq C _1 \| \nabla u \|_{L^2(\R^N)} ^{\frac{2^*}{2}} , 
\eeq
\vspace*{-10pt}
\beq
\label{2.7}
\ds \int_{\R^N} \left( \ph ^2( | 1+ u_1|) - 1 \right)^2 \, dx 
\leq 
\int_{\R^N} \left( \ph ^2( | 1+ u|) - 1 \right)^2 \, dx 
+ C_2  \| \nabla u \| _{L^2(\R^N)} ^{2^*}, 
\eeq
\vspace{-10pt}
\beq
\label{2.8}
\int_{\R^N} \left( \ph ^2( | 1+ u_2|) - 1 \right)^2 \, dx 
\leq 
 C_2  \| \nabla u \|_{L^2(\R^N)} ^{2^*}.
\eeq
Let $ 1+ u_1 = \rho e^{i \theta} $ be the lifting of $ 1+ u_1$, as given by 
Lemma  \ref{lifting}. Then we have
\beq
\label{2.9}
\langle iu_{x_1}, u \rangle  = 
( 1 - \chi ^2 (u)) \langle iu_{x_1}, u \rangle 
- ( \rho ^2 - 1 ) \frac{ \p \theta}{\p x _1} 
+  \frac{ \p }{\p x_1 } \left( \I (u ) \right) - \frac{\p \theta}{\p x_1}. 
\eeq
a.e. on $ \R^N$.
\end{Lemma}

{\it Proof. } 
Since $ |u_i | \leq |u|$, we have $ u_i \in L^{2^*}( \R^N)$ for  $ i=1, 2$. 
It is standard to prove that $ u_i \in H_{loc}^1 ( \R^N) $ (see, e.g., Lemma C1 p. 66 
in \cite{BBM}) and we have
\beq
\label{2.10}
\frac{\p u _1}{\p x_j } = \left( \p_1 \chi (u) \frac{\p ( \Real (u))}{\p x_j } + 
\p_2 \chi (u) \frac{\p ( \I (u))}{\p x_j } \right) u + \chi (u) 
\frac{\p u}{\p x_j }.
\eeq
A similar formula holds for $u_2$. 
Since the functions $ z \longmapsto \p_i \chi (z) z $, $ i=1, 2$, are bounded on $ \C$, 
(\ref{2.5}) follows immediately from (\ref{2.10}). 

Using the Sobolev embedding we have
$$
\| u_2 \|_{L^2}^2 \leq {\ds \int_{\R^N} } |u|^2 \1_{\{ |u| > \frac{1}{4}\} } (x) \, dx 
\leq  4 ^{2^* -2} 
{\ds \int_{\R^N} }
 |u|^{2^*} \1_{\{|u| > \frac{1}{4}\} } (x) \, dx 
\leq C_1 \|\nabla u \|_{L^2} ^{2^*}. 
$$
This gives the first estimate in (\ref{2.6}); the second one is similar. 

For $|u| \leq \frac{1}{4}$ we have $ u_1 (x) = u(x)$, hence
$$
{\ds \int_{ \{|u| \leq \frac{1}{4} \} } } \left( \ph ^2( |1+ u_1|) - 1 \right)^2 \, dx 
= { \ds \int_{\{|u| \leq \frac{1}{4} \} } } 
\left( \ph ^2( |1+ u|) - 1 \right)^2 \, dx.
$$
There exists $ C' > 0 $ such that 
$\left( \ph ^2( |1+z|) - 1 \right)^2 \leq C' |z|^2$  if $ |z| \geq \frac{1}{4}$.
Proceeding as in the proof of (\ref{2.6}) we have for $ i =1, 2$
$$
\ds \int_{ \{ |u| > \frac{1}{4} \} } \left( \ph ^2( |1+ u_i|) - 1 \right)^2 \, dx 
\leq C' 
\ds \int_{\{ |u| > \frac{1}{4} \} } |u_i|^2  \, dx  
\leq 
C_2 \| \nabla u \| _{L^2} ^{2^*}. 
$$
This clearly implies (\ref{2.7}) and (\ref{2.8}).

Since $ \p_1 \chi (u) \frac{\p ( \Real (u))}{\p x_j } + 
\p_2 \chi (u) \frac{\p ( \I (u))}{\p x_j } \in \R$, 
using (\ref{2.10}) we get
$ \langle i\frac{\p u_1}{\p x_1}, u _1 \rangle = \chi ^ 2 (u) \langle iu_{x_1}, u \rangle $
a.e. on $ \R$. 
Then (\ref{2.9}) follows from Lemma \ref{lifting}. 
\hfill $\Box$

\medskip

We consider the space $ \Yo = \{ \p _{x_1} \phi \; | \; \phi \in \dot{H}^1 (\R^N) \}$.
It is clear that $ \phi _1, \phi _2 \in \dot{H}^1 (\R^N)  $ and 
$ \p _{x_1} \phi _1 = \p _{x_1} \phi _ 2$ imply that $ \phi _1 - \phi _2$ is constant, hence $ \nabla \phi _1 = \nabla \phi _2$.  
Defining 
$$
\| \p _{x_1} \phi \|_{\Yo } = \|\phi \|_{\dot{H}^1 (\R^N)  } = \|\nabla \phi \|_{L^2(\R^N)}, 
$$
it is easy to see that $\|\cdot \|_{\Yo } $ is a norm on $ \Yo $ and 
$(\Yo, \|\cdot \|_{\Yo }  )$ is a Banach space. 
The~following~holds. 

\begin{Lemma} 
\label{L2.3}
Let $ N \geq 2$. For any $ v \in L^1( \R^N) \cap \Yo $ we have $ \ds \int_{\R^N} v(x) \, dx = 0$. 
\end{Lemma} 

{\it Proof. } 
Take $ \phi \in \dot{H}^1 (\R^N) $  such that $ v = \p _{x_1} \phi $. 
Then $ \phi \in \So ' ( \R^N) $ and $ |\xi | \wh{ \phi } \in L^2 ( \R^N)$. 
Hence $ \wh{\phi} \in L_{loc}^1( \R^N \setminus \{ 0 \} ) $. 
On the other hand we have 
$ v =  \p _{x_1} \phi \in L^1 \cap L^2 (\R^N) $ by hypothesis, hence 
$ \wh{v} = i \xi _1 \wh{ \phi } \in L^2 \cap C_b ^0 ( \R^N)$. 

We prove that $ \wh{v}(0) = 0$. 
We argue by contradiction and assume that $ \wh{v}(0) \neq 0$.
By continuity,  there exists $ m > 0$ and $ \e > 0$ such that 
$| \wh{v}(\xi )| \geq m $ for $ |\xi | \leq \e$. 
For  $j =2, \dots N$ we get 
$$
| i \xi _j \wh{\phi }(\xi )| = \frac{|\xi _j| }{|\xi _1|} | \wh{v}(\xi )| 
\geq m \frac{|\xi _j| }{|\xi _1|}
\qquad \mbox{ for a.e. } \xi \in B(0, \e).
$$
But this contradicts the fact that
 $ i \xi _j \wh{\phi } 
 \in L^2 ( \R^N).$ 
Thus necessarily $ \wh{v}(0) = 0$ and this is exactly the conclusion of Lemma \ref{L2.3}.
\hfill 
$\Box $

\medskip

It is obvious that $L_1(v) = \ds \int_{\R^N} v(x) \, dx $ and $ L_2(w) = 0$ 
are continuous linear functionals  on $L^1( \R^N) $ and on $ \Yo$, respectively. 
Moreover, by Lemma \ref{L2.3} we have $ L_1 = L_2 $ on $ L^1( \R^N ) \cap \Yo$. 
Putting
\beq 
\label{2.11}
L(v + w) = L_1 (v) + L_2 (w) = \ds \int_{\R^N} v(x) \, dx
\qquad \mbox{ for } v \in L^1( \R^N) \mbox{ and } w \in \Yo, 
\eeq
we see that $ L$ is well-defined and is a continuous linear functional  on 
$L^1( \R^N) + \Yo$. 

It follows from (\ref{2.9}) and Lemmas \ref{lifting} and \ref{L2.2} that 
for any $ u \in \Xo $ we have $ \langle iu_{x_1}, u \rangle \in L^1( \R^N) + \Yo$.
This enables us to give the following

\begin{Definition}
\label{D2.4}
Given $ u \in \Xo$, the momentum of $u$ (with respect to the
$x_1-$direction) is 
$$
Q(u) = L( \langle iu_{x_1}, u \rangle  ).
$$
\end{Definition}

If $ u \in \Xo $ and $ \chi, u_1, u_2, \rho, \theta $ are as in Lemma \ref{L2.2}, 
from (\ref{2.9}) we get 
\beq
\label{2.12}
Q(u) = \ds \int_{\R^N} ( 1 - \chi ^2(u))  \langle iu_{x_1}, u \rangle  - 
(\rho ^2 - 1 ) \theta _{x_1} \, dx . 
\eeq
It is easy to check that the right-hand side of (\ref{2.12}) does not depend on the  
choice of the cut-off function $ \chi$, provided that $ \chi $ is as in Lemma \ref{L2.2}.

It follows directly from (\ref{2.12}) that the functional $Q$ 
has a nice behavior with respect to 
dilations in $ \R^N$: for any $ u \in \Xo $ and $ \la, \, \si  > 0 $ we have 
\beq
\label{2.13}
Q(u_{\la,  \si}) = \si ^{N-1} Q(u).
\eeq

The next lemma will enable us to perform "integrations by parts". 

\begin{Lemma}
\label{L2.5}
For any $ u \in \Xo$ and $ w \in H^1 ( \R^N) $ we have 
$\langle iu_{x_1}, w \rangle \in L^1( \R^N)$, 
$\langle iu, w_{x_1}  \rangle  \in L^1( \R^N) + \Yo $ and 
\beq
\label{2.13ipp}
L(\langle iu_{x_1}, w \rangle +  \langle iu, w_{x_1}  \rangle) = 0.
\eeq
\end{Lemma}

{\it Proof.}
 
Since $ w, u_{x_1}  \in L^2( \R^N)$, the Cauchy-Schwarz inequality implies 
$\langle iu_{x_1}, w \rangle \in L^1( \R^N)$.

Let $ \chi, \; u_1, \; u_2 $ be as in Lemma \ref{L2.2} and denote 
$ w_1 = \chi(w) w$, $w_2 = ( 1 - \chi(w)) w$. 
Then $ u = u_1 + u_2$, $ w = w_1 + w_2 $ and it follows from Lemma 
 \ref{L2.2} that $ u_1 \in \Xo \cap L^{\infty }( \R^N) $ and 
$ u_2,\;  w_1, \; w_2 \in H^1 ( \R^N)$.

As above we have $\langle i \frac{ \p u_2}{\p x_1}, w \rangle, \;
\langle i u_2 , \frac{ \p w }{\p x_1} \rangle \in L^1 ( \R^N)$ 
by the Cauchy-Schwarz inequality. 
The standard integration by parts formula for  functions in $H^1(\R^N)$ 
(see, e.g., \cite{brezis}, p. 197) gives
\beq
\label{2.14}
\ds \int_{\R^N} \langle i \frac{ \p u_2}{\p x_1}, w \rangle
+ \langle i u_2 , \frac{ \p w }{\p x_1} \rangle \, dx = 0.
\eeq
Since $ u_1 \in \D \cap L^{\infty}( \R^N) $ and $ w_1 \in H^1 \cap L^{\infty}( \R^N) $, 
it is standard to prove that 
$ \langle iu_1, w_1 \rangle \in \D \cap L^{\infty}( \R^N) $ and
\beq
\label{2.15}
\langle i\frac{\p u_1}{\p x_1} , w_1 \rangle
+ \langle iu_1, \frac{\p w_1}{\p x_1}  \rangle
= \frac{\p}{\p x_1} \langle iu_1, w_1 \rangle
\qquad \mbox{ a.e. on } \R^N.
\eeq

Let $ A_w = \{ x \in \R^N \; | \; |w(x)| \geq \frac{1}{4} \}$. 
We have 
$ 
\frac{1}{16} \Lo ^N ( A_w) 
\leq {\ds  \int_{A_w}} |w|^{2} \, dx  \leq \| w\|_{L^{2}}^2$, 
and consequently $ A_w$ has finite measure. 
It is clear that $ w_2 = 0 $ and $ \nabla w _2 = 0$ a.e. on $ \R^N \setminus A_w$. 
Since $ w_2 \in L^{2^*}( \R^N) $ and $ \nabla w_2 \in L^2( \R^N)$, 
we infer that
$ w_2 \in L^1 \cap L^{2^*} ( \R^N)$ and $ \nabla w _2 
\in  L^1 \cap L^{2} ( \R^N)$.
Together with the fact that $ u_1 \in L^{2^*} \cap L^{\infty }(\R^N)$
and $ \nabla u_1 \in L^2( \R^N)$, this gives
$\langle iu _1, w_2 \rangle \in L^1 \cap L^{2^*} ( \R^N)$ and
$$
\langle i\frac{\p u _1}{\p x_j}, w_2 \rangle \in L^1 \cap L^{\frac{N}{N-1}}(\R^N), 
\qquad
\langle iu _1, \frac{\p w_2}{\p x_j}  \rangle \in L^1 \cap L^2 ( \R^N)
\qquad
\mbox{ for } j=1, \dots, N.
$$
It is easy to see that 
$$
\frac{ \p }{\p x_j } \langle iu _1, w_2 \rangle = 
\langle i\frac{\p u _1}{\p x_j}, w_2 \rangle + 
\langle iu _1, \frac{\p w_2}{\p x_j}  \rangle \qquad \mbox{ in } \Do '(\R^N).
$$
From the above we infer that $  \langle iu _1, w_2 \rangle \in W^{1,1}( \R^N)$. 
It is obvious that $ \ds \int _{\R^N} \frac{\p \psi}{\p x_j} \, dx = 0$ for any 
$ \psi \in W^{1,1}( \R^N)$ 
(indeed, let $(\psi _n)_{n \geq 1} \subset C_c^{\infty }(\R^N)$ be a sequence such that 
$ \psi _n \lra \psi $ in $ W^{1,1}( \R^N)$ as $ n \lra \infty$; 
then  $ \ds \int _{\R^N} \frac{\p \psi_n }{\p x_j} \, dx = 0$ for
each $n$ and     
$ \ds \int _{\R^N} \frac{\p \psi_n }{\p x_j} \, dx \lra  
 \int _{\R^N} \frac{\p \psi }{\p x_j} \, dx $ as $ n \lra \infty$). 
Thus we have 
$\langle i\frac{\p u _1}{\p x_1}, w_2 \rangle, \; 
\langle iu _1, \frac{\p w_2}{\p x_1}  \rangle \in L^1( \R^N) $ and 
\beq
\label{2.16}
\ds \int_{\R^N} \langle i\frac{\p u _1}{\p x_1}, w_2 \rangle + 
\langle iu _1, \frac{\p w_2}{\p x_1}  \rangle\, dx 
= \int_{\R^N} \frac{ \p }{ \p x_1} \langle iu _1, w_2 \rangle \, dx = 0 .
\eeq
Now (\ref{2.13ipp}) follows from (\ref{2.14}), (\ref{2.15}), (\ref{2.16})
and Lemma \ref{L2.5} is proven. 
\hfill 
$ \Box $ 

\begin{Corollary} 
\label{C2.6}
Let $ u, \; v \in \Xo $ be such that $ u -v \in L^2( \R^N)$. 
Then 
\beq
\label{2.17}
|Q(u) - Q(v) | \leq \| u -v \|_{L^2(\R^N)}
\left( \Big\| \frac{\p u}{\p x_1} \Big\| _{L^2(\R^N)} 
+ \Big\| \frac{\p v}{\p x_1} \Big\| _{L^2(\R^N)} \right) .
\eeq
\end{Corollary}

{\it Proof. } 
It is clear that $ w = u -v \in H^1( \R^N)$ and using (\ref{2.13ipp}) we get 
\beq
\label{2.18}
\begin{array}{rcl}
Q(u) - Q(v) & = & 
L( \langle i (u-v)_{x_1} , u \rangle + \langle i v_{x_1}, u-v \rangle ) 
\\
& = & L( \langle i u_{x_1} , u -v \rangle + \langle i v_{x_1}, u-v \rangle ) 
\\
& = &  \ds \int_{\R^N}  \langle i u_{x_1} + i v_{x_1}, u -v \rangle \, dx.
\end{array}
\eeq
Then (\ref{2.17}) follows from (\ref{2.18}) and the Cauchy-Schwarz inequality. 
\hfill
$\Box $

\medskip

The next result will be useful to estimate the contribution to the momentum 
of a domain where the modified Ginzburg-Landau energy is small. 

\begin{Lemma}
\label{L2.7}
Let $ M > 0$ and let $ \Om $ be an open subset of $ \R^N$. 
Assume that  $ u \in \Xo $ satisfies $ E_{GL}(u) \leq M $ and let 
$ \chi, \; \rho, \;  \theta $ be 
as in Lemma \ref{L2.2}. Then we have
\beq
\label{2.19}
\ds \int_{\Om } \Big\vert ( 1 - \chi ^2 (u)) \langle iu_{x_1}, u \rangle 
- ( \rho ^2 - 1 )  \theta _{ x _1} \Big\vert \, dx 
\leq  C( M^{\frac 12} + M^{\frac{2^*}{4} } )  \left( E_{GL}^{\Om }(u) \right)^{\frac 12}.
\eeq
\end{Lemma}

{\it Proof. } Using (\ref{2.6}) and the Cauchy-Schwarz inequality we get 
\beq
\label{2.20}
\begin{array}{rcl}
\ds \int_{\Om } \Big\vert ( 1 - \chi ^2 (u)) \langle iu_{x_1}, u \rangle  \Big\vert \, dx 
& \leq & 
\|u_{x_1}\|_{L^2 ( \Om )}  \|(1 - \chi ^2(u))u  \|_{L^2 ( \Om ) }
\\
& \leq & 
C_1 \|u_{x_1}\|_{L^2( \Om )} \|\nabla u \|_{L^2( \R^N)} ^{\frac{2^*}{2}}.
\end{array}
\eeq
We have  $|u_1| \leq \frac{1}{2}$, hence $|1+ u_1  | \leq \frac{3 }{2} $ and 
$ \ph ( |1+ u_1|) = |1+ u_1| = \rho$. 
Then (\ref{2.7})  gives
\beq
\label{2.21}
\|\rho ^2 - 1 \| _{L^2( \R^N)} ^2 \leq  C' (E_{GL}(u)+   E_{GL}(u) ^{\frac{2^*}{2}} ) 
\leq 
C'( M + M^{\frac{2^*}{2}}).
\eeq
From (\ref{2.4}) and (\ref{2.5}) we have 
$\big\vert \frac{\p \theta}{\p x_j } \big\vert
\leq \frac{1}{\rho } \big\vert \frac{\p u_1}{\p x_j } \big\vert 
\leq C'' \big\vert \frac{\p u }{\p x_j } \big\vert $ a.e. on $\R^N$. Therefore
\beq
\label{2.22}
\begin{array}{l}
\ds \int_{\Om } \Big\vert  ( \rho ^2 - 1 )  \theta_{ x _1} \Big\vert \, dx 
  \leq   
\| \rho ^2 - 1 \|_{L^2(\Om )} \|\theta _{x_1} \|_{L^2(\Om )} 
\\
\leq C '' \| \rho ^2 - 1 \|_{L^2(\R^N)} \| u_{x_1} \|_{L^2(\Om )} 
 \leq   C '''  \left( M +  M^{\frac{2^*}{2}} \right)^{\frac 12}
\left( E_{GL}^{\Om }(u) \right)^{\frac 12}. 
\end{array}
\eeq
Then (\ref{2.19}) follows from (\ref{2.20}) and (\ref{2.22}). 
\hfill
$\Box$

\section{A regularization procedure} 

Given a function $ u \in \Xo $ and a set $ \Om \subset \R^N$ 
such that  $E_{GL}^{\Om }(u)$ is small, we would like to get a fine  estimate of 
the contribution of $ \Om $ to the momentum of $u$. 
To do this, we will use a kind of "regularization" procedure for arbitrary functions 
in $ \Xo$. 
A similar device has been introduced in \cite{AB} to get rid of 
small-scale topological defects of functions; variants of it have been used for 
various purposes in \cite{BS}, \cite{BOS}, \cite{BGS}.

Throughout this section, $\Om $ is an open set in $ \R^N$. 
We do not assume $ \Om $ bounded, nor connected. 
If $ \p \Om \neq \emptyset$, we assume that $ \p \Om $ is $ C^2$. 
Let $\ph $ be as in the introduction.
Fix $ u \in \Xo $ and  $  h > 0$. 
We consider the functional 
$$
G_{h, \Om }^u (v) = E_{GL}^{\Om } (v) + 
\frac{ 1}{h^2} \ds \int_{\Om } \ph \left( |v-u|^2  \right) \, dx. 
$$
Note that $G_{h, \Om }^u (v) $ may equal  $ \infty $ for some  $ v \in \Xo $;
however, $G_{h, \Om }^u (v) $ is finite whenever $ v \in \Xo $ and 
$ v - u \in L^2( \Om )$. 
If there is no risk of confusion, we will simply write $G(v)$ instead of $G_{h, \Om }^u (v) $.
We denote 
$ H_0^1(\Om ) = \{ u \in H^1( \R^N) \; | \; u = 0 \mbox{ on } \R^N \setminus \Om \}$ 
and 
$$
H_u ^1 ( \Om ) = \{ v \in \Xo \; | \; v -u \in H_0^1 ( \Om ) \}.
$$
The next lemma gives the properties of functions that minimize 
$G $ in the space $ H_u ^1 ( \Om ) $. 

\begin{Lemma}
\label{L3.1}
i) The functional $ G $ has a minimizer in $H_u ^1 ( \Om ) $. 

ii) Let $ v_h$ be a minimizer of $ G $ in $H_u ^1 ( \Om ) $. 
There exist constants $ C_1, \, C_2 > 0$, depending only on 
$N,$  such that 
 $v_h$ satisfies: 
\beq
\label{3.1}
E_{GL}^{\Om } ( v_h )\leq E_{GL}^{\Om } ( u ); 
\eeq
\vspace{-10pt}
\beq
\label{3.2}
\| v_h - u \|_{L^2( \Om ) } ^2  \leq  h^2 E_{GL}^{\Om } ( u ) 
 + C _1 \left( E_{GL}^{\Om } ( u ) \right)^{1 + \frac 2N} h^{\frac 4N}; 
\eeq
\vspace{-20pt}
\beq
\label{3.3}
\ds \int_{\Om } \Big\vert \left( \ph ^2( |1+ u|) - 1 \right)^2 - 
  \left( \ph ^2( |1+ v_h |) - 1 \right)^2 \Big\vert \, dx 
\leq 36 h E_{GL}^{\Om } ( u ) ; 
\eeq
\vspace{-20pt}
\beq
\label{3.4} 
| Q(u) - Q( v_h) | \leq C_2 \left(  h^2 +  \left( E_{GL}^{\Om } ( u ) \right)^{ \frac 2N} 
h^{\frac 4N} \right)^{\frac 12}
E_{GL}^{\Om } ( u ).
\eeq

iii) For $ z \in \C$, denote 
$H(z) = \left( \ph ^2( |z +1|) - 1 \right) \ph ( |z +1|) \ph '( |z +1|) 
\frac{z +1}{ |z +1|}$  if $ z \neq -1$ and $H(-1) = 0$. 
Then any minimizer $v_h $ of $G $ in $H_u^1( \Om )$ satisfies the equation
\beq
\label{3.5}
- \Delta v_h +  H( v_h ) 
+ \frac{ 1}{ h^2 } \ph '  \left(|v_h -u|^2 \right)
(v_h -u)= 0 
\qquad \mbox{ in } \Do ' ( \Om ) . 
\eeq
Moreover, for any $ \omega \subset \subset \Om $ we have $ v_h \in W^{2, p} ( \omega ) $ 
for   $ p \in [1, \infty )$; thus, in particular, $ v_h \in C^{1, \al } (\omega )$
for $ \al \in [0, 1)$. 

iv) For any $ h>0$, $ \de >0$ and $ R> 0$ there exists a constant 
$K = K (  N, h, \de, R) > 0 $ such that  for any $ u \in \Xo $ with 
$E_{GL}^{\Om } ( u) \leq K$ and for any minimizer 
$v_h$ of $G $ in $H_u^1( \Om )$ there holds
\beq
\label{3.6}
1 - \de < |1+ v_h(x) | < 1 + \de 
\qquad \mbox{ whenever } x \in \Om \mbox{ and } dist(x, \p \Om ) \geq 4R.
\eeq
\end{Lemma}

{\it Proof. } 
i) It is obvious that $ u \in H_u^1( \Om )$. 
Let $ (v_n)_{n \geq 1}$ be a minimizing sequence for $G $ in $H_u^1( \Om )$.
We may assume that $G (v_n) \leq G (u) = E_{GL}^{\Om} (u)$ 
and this implies 
$\ds \int_{\Om } |\nabla v_n |^2 \, dx \leq E_{GL}^{\Om} (u)$.
It is clear that 
\beq
\label{3.7}
\ds \int_{\Om \cap \{ |v_n - u| \leq \sqrt{2}  \} } |v_n - u |^2 \, dx 
\leq  \ds \int_{\Om } \ph \left( |v_n - u|^2 \right) \, dx 
\leq  h^2 E_{GL}^{\Om} (u).
\eeq
Since $ v_n - u \in H_0^1 ( \Om ) \subset H^1( \R^N)$, by the Sobolev embedding 
we have 
$\|v_n - u\|_{L^{2^*}(\R^N)} \leq  C_S \| \nabla v_n - \nabla u \|_{L^2( \R^N)}$, where 
$C_S$ depends only on $N$.  Therefore 
\beq
\label{3.8} 
\begin{array}{l}
\ds \int_{ \{ |v_n - u| \geq 1  \} } |v_n - u |^2 \, dx 
\leq 
\ds \int_{ \{ |v_n - u| \geq 1 \} } |v_n - u|^{2^*} \, dx 
\\
\leq  \|  v_n - u \|_{L^{2^*}(\R^N)} ^{ 2^*}
\leq C ' \| \nabla v_n - \nabla u \|_{L^2( \R^N)} ^{2^*} 
\leq C \left( E_{GL}^{\Om} (u) \right)^{\frac{2^*}{2}}.
\end{array}
\eeq
It follows from (\ref{3.7}) and (\ref{3.8}) that $\|v_n - u \|_{L^2( \Om )} $ 
is  bounded, hence $v_n - u $ is bounded in $H_0^1( \Om )$. 
We infer that there exists a sequence (still denoted $(v_n)_{n \geq 1}$) 
and there is $ w \in H_0^1( \Om ) $ such that 
$ v_n - u \rightharpoonup w $ weakly in $H_0^1( \Om ) $, 
$  v_n - u \lra w $ a.e. and $  v_n - u \lra w $ in $L_{loc}^p ( \Om ) $ 
for $1 \leq p < 2^*$. 
Let $ v = u + w$.  Then $\nabla v_n \rightharpoonup \nabla v$ weakly in $ L^2( \R^N)$ 
and this implies
$$
\ds \int_{\Om } |\nabla v |^2 \, dx 
\leq  \liminf _{ n \ra \infty } \ds \int_{\Om } |\nabla v_n |^2 \, dx .
$$
Using the a.e. convergence and Fatou's Lemma we infer that 
$$
\ds \int_{\Om } \left( \ph ^2( |1+ v |) - 1 \right)^2 \, dx 
\leq  \liminf _{ n \ra \infty } \ds \int_{\Om } 
\left( \ph ^2( | 1+ v _n|) - 1 \right)^2 \, dx 
\qquad \mbox{ and } 
$$
\vspace{-5pt}
$$
\ds \int_{\Om } 
\ph \left( |v- u|^2 \right) \, dx 
\leq  \liminf _{ n \ra \infty }
\ds \int_{\Om } 
\ph \left( |v_n - u|^2 \right) \, dx.
$$
Therefore $G (v)
\leq  \ds \liminf _{ n \ra \infty }  G (v_n) $ 
and consequently $v$ is a minimizer of $G $ in $H_u^1 (\Om )$. 

\medskip

ii) 
Since $ u \in H_u^1 (\Om )$, we have 
$ E_{GL}^{\Om } (v_h) \leq G  (v_h) \leq E_{GL}^{\Om } (u) $; hence   (\ref{3.1}) holds. 
It is clear that $\ph \left( |v_h -u|^2 \right)  \geq 1  $ if 
$|v_h -u| \geq 1 $, thus 
\vspace{-7pt}
$$
 \Lo ^N( \{ | v_h -u | \geq 1  \}) 
\leq \int_{\R^N} \ph \left(  |v_h -u|^2 \right) \, dx 
\leq h^2 G ( v_h) \leq h^2 E_{GL}^{\Om } (u) .
$$
Using  H\"older's inequality, the above estimate and the Sobolev inequality   we get 
\beq
\label{3.9} 
\begin{array}{l}
\ds \int_{ \{ |v_h - u| \geq 1  \} } |v_h - u |^2 \, dx 
\\
\leq \| v_h - u\|_{L^{2^*} ( \{ |v_h - u| \geq 1 \} )}^ 2 
\left(\Lo ^N ( \{ | v_h -u | \geq 1  \}) \right) ^{1 - \frac{2}{2^*}} 
\\
\leq \|v_h - u\|_{L^{2^*}(\R^N)} ^2 
\left( \Lo ^N ( \{ | v_h -u | \geq 1  \}) \right) ^{1 - \frac{2}{2^*}} 
\\
\leq C_S \|\nabla v_h - \nabla u \|_{L^2(\R^N)} ^2 
\left( h^2  E_{GL}^{\Om } (u) \right) ^{1 - \frac{2}{2^*}}  
\leq 4 C_S h^{\frac 4N} \left(  E_{GL}^{\Om } (u) \right)^{1 + \frac 2N}.
\end{array}
\eeq
It is clear that (\ref{3.7}) holds with $v_h$ instead of $ v_n$ 
and then (\ref{3.2}) follows from (\ref{3.7}) and (\ref{3.9}). 

We claim that 
\beq
\label{3.10}
\begin{array}{c}
\Big\vert 
\ph(| z|) - \ph(| \zeta|) \Big\vert
\leq \left[ \frac 92 \ph \left(  |z - \zeta|^2 \right) \right]^{\frac 12}
\qquad \mbox{ for any } z, \; \zeta \in \C.
\end{array}
\eeq
Indeed, let $0\leq a \leq b$. 
If $ b \in [a, a + \sqrt{2}]$ we have $\ph ((b-a )^2) = (b-a)^2 $, hence 
$$ 0 \leq \ph(b) - \ph (a) \leq b-a = \left[ \ph( (b-a)^2) \right]^{\frac 12}.$$
If $ b > a +\sqrt{2}$ we have $ 0 \leq \ph(b) - \ph (a) \leq 3 $ and 
$\left[ \ph( (b-a)^2) \right]^{\frac 12} \geq \left( \ph(2) \right) ^{\frac 12} = \sqrt{2}$, thus
$ 0 \leq \ph(b) - \ph (a) \leq \frac{3}{\sqrt{2}} \left[ \ph( (b-a)^2) \right]^{\frac 12} .$ 
Assuming that $ |z| \leq |\zeta |$ we get 
$$
\begin{array}{c}
\Big\vert  \ph(| z|) - \ph(| \zeta|) \Big\vert =  \ph(| \zeta|) -  \ph(| z|) 
\leq \left[ \frac 92 \ph \left( ( |\zeta| - |z|) ^2 \right) \right]^{\frac 12}
\leq \left[ \frac 92 \ph \left(  | \zeta -z|^2 \right) \right]^{\frac 12}.
\end{array}
$$

It is obvious that 
\beq
\label{3.11}
\begin{array}{l}
\Big\vert \left( \ph ^2( |1+ u |) - 1 \right)^2 
- \left(\ph ^2( | 1+ v_h |) - 1 \right)^2 \Big\vert
\\
\leq 6  \Big\vert  \ph ( | 1+ u |) - \ph ( | 1+ v_h  |) \Big\vert
\cdot \Big\vert \ph ^2( | 1+ u |) + \ph ^2( | 1+ v_h |) - 2  \Big\vert.
\end{array}
\eeq
Using (\ref{3.11}), the Cauchy-Schwarz inequality and  (\ref{3.10})  we get 
\vspace{-3pt}
$$
\begin{array}{l}
\ds \int_{\Om } \Big\vert \left( \ph ^2( |1+ u|) - 1 \right)^2 - 
  \left( \ph ^2( |1+ v_h |) - 1 \right)^2 \Big\vert \, dx 
\\
\\
\leq 6 
\left( \ds \int_{\Om } 
\Big\vert  \ph ( |1+ u|) - \ph ( |1+ v_h|) \Big\vert ^2\, dx \right)^{\frac 12}
\left( \ds \int_{\Om } 
\Big\vert  \ph ^2( |1+ u|) + \ph ^2( |1+ v_h|) -2  \Big\vert ^2\, dx \right)^{\frac 12}
\\
\\
\leq 6  \left( {\ds \int_{\Om }} \frac 92
\ph \left(|v_h - u|^2 \right)\, dx \right)^{\frac 12}
\left( 2{ \ds \int_{\Om } }
\left( \ph ^2( |1+ u|) - 1 \right) ^2 
+\left(  \ph ^2( |1+ v_h|) -1 \right)^2  \,  dx \right)^{\frac 12}
\\
\\
\leq 18 \left( h^2 G ( v_h ) \right)^{\frac 12} 
\left( 2 E_{GL}^{\Om } (u) +  2 E_{GL}^{\Om }(v_h) \right)^{\frac 12} 
\leq 36 h E_{GL}^{\Om } (u)
\end{array}
$$
and (\ref{3.3}) is proven. 
Finally, (\ref{3.4}) follows directly from (\ref{3.1}), (\ref{3.2}) 
and Corollary \ref{C2.6}. 

\medskip

iii) The proof of (\ref{3.5}) is standard. For any 
$ \psi \in C_c^{\infty }(\Om )$ we have $ v + \psi \in H_u^1( \Om )$ 
and the function $ t \longmapsto G  ( v + t \psi )$ achieves its minumum 
at $ t =0$. 
Hence $\frac{d}{dt }_{\big\vert _{t =0}} \left(G ( v + t \psi )\right) =0$
for any $ \psi \in  C_c^{\infty }(\Om ) $ and this is precisely (\ref{3.5}).

For any $ z \in \C$ we have 
\beq
\label{3.12}
|H(z) | \leq 3  | \ph ^2 ( |z +1|) - 1| \leq 24  . 
\eeq
Since $ v_h \in \Xo$, we have $ \ph ^2 (| 1+ v_ h |) - 1 \in L^2 ( \R^N) $  
and (\ref{3.12})  gives $ H( v_h ) \in L^{2}\cap L^{\infty } ( \R^N)$.
We also have 
$\big\vert  \ph ' \left( |v_h -u|^2 \right) (v_h - u ) \big\vert 
\leq |v_h - u| $ and
$$
\big\vert \ph '\left( |v_h -u|^2  \right)(v_h - u ) \big\vert  
 \leq  {\ds \sup_{s \geq 0 } }\,  \ph ' \left( s^2    \right) s < \infty.
$$ 
Since $ v_h - u \in L^2( \R^N)$, we get 
$  \ph '\left( |v_h -u|^2 \right)(v_h - u ) \! \in \!  L^{2}\cap L^{\infty } ( \R^N)$.
Using  (\ref{3.5}) we infer  that $ \Delta v_h \in L^2 \cap L^{\infty }( \Om )$. 
Then (iii) follows from standard elliptic estimates 
 (see, e.g., Theorem 9.11 p. 235 in \cite{GT})
 and a straightforward 
bootstrap argument.

\medskip

iv) We use (\ref{3.5}), Sobolev and Gagliardo-Nirenberg inequalities  and elliptic regularity theory to prove that there is 
$r \leq R$ such that for all $x$ satisfying $B(x, 4R) \subset \Om$, one may estimate $\| \nabla v_h \|_{L^p (B(x, r)) }$ 
 in terms of $E_{GL}^{\Om} (u)$ (see (\ref{3.26}) below).
This estimate with $ p > N$ and the Morrey inequality imply that $v_h$ is uniformly H\"older continuous on 
$\{ x \in \Om \; | \; dist(x, \p \Om ) \geq 4 R \}$. 
In particular, if $\big| \, | 1 + v_h (x_0)| -1  \big| > \de $ for some $ x_0 $ verifying $B(x_0, 4R) \subset \Om$, 
then necessarily $\big| \, | 1 + v_h | -1  \big| > \frac{\de}{2} $ on a ball of fixed radius 
centered at $ x_0$ and this implies that $E_{GL}^{\Om} (v_h)$ (and, consequently,  $ E_{GL}^{\Om} (u)$) 
is bounded from below by a positive constant.

We start by estimating the nonlinear terms in (\ref{3.5}).
Using (\ref{3.12}) we get 
$$
\ds \int_{\Om}  |H(v_h)|^2 \, dx \leq 9  \int_{\Om } 
\left(  \ph ^2 ( |1+ v_h |) - 1 \right)^2\, dx 
\leq 18 E_{GL}^{\Om } (v_h) 
\leq 18 E_{GL}^{\Om } (u),  
$$
hence $\| H(v_h) \|_{L^2(\Om )} \leq 3 \sqrt{2} \left( E_{GL}^{\Om } (u) \right)^{\frac 12}$. 
By interpolation we find for any $p \in [2, \infty]$, 
\beq 
\label{3.13}
\| H(v_h)\|_{L^p(\Om )} \leq \|H(v_h)\|_{L^{\infty} (\Om )} ^{\frac{p-2}{p}} 
\|H(v_h)\|_{L^2(\Om )} ^{\frac 2p} \leq C   \left( E_{GL}^{\Om } (u) \right)^{\frac 1p}.
\eeq

It is easy to see that $ \big\vert \ph'(s^2) s \big\vert ^2 \leq 2 \ph ( s^2) $ and $\big\vert \ph'(s^2) s \big\vert \leq 2$
for  any $ s \geq 0$. 
Then we have 
$$
{\ds \int_{ \Om }} 
\big\vert \ph ' \left( |v_h - u |^2 \right) ( v_h - u ) \big\vert ^2 \, dx 
\leq 
2 {\ds \int_{ \Om }}  
  \ph  \left(  |v_h - u|^2   \right) \, dx 
  \leq 2 h^2   E_{GL}^{\Om } (u), 
$$
thus $ \| \ph ' \left( |v_h - u |^2 \right) ( v_h - u ) \| _{L^2( \Om )} 
\leq h \left(2 E_{GL}^{\Om } (u) \right)^{\frac 12}$. 
By  interpolation we get  
\beq
\label{3.14}
\begin{array}{l}
\big\| \ph ' \left( |v_h - u |^2 \right) ( v_h - u ) \big \| _{L^p( \Om )} 
\\
\leq 
\big\| \ph ' \left(  |v_h - u |^2  \right) ( v_h - u ) \big\| _{L^{\infty} ( \Om )} ^{\frac{p-2}{p}}
\big\| \ph ' \left(|v_h - u |^2 \right) ( v_h - u )  \big\| _{L^2( \Om )} ^{\frac 2p}
\\
\leq 
C h^{\frac 2p}
 \left( E_{GL}^{\Om } (u) \right)^{\frac 1p}
\end{array}
\eeq
 for any $ p \in [2, \infty]$.
From (\ref{3.5}), (\ref{3.13}) and (\ref{3.14}) we obtain 
\beq
\label{3.15}
\|\Delta v_h \|_{L^p( \Om )}  
\leq C ( 1 + h^{ \frac 2p - 2} ) \left( E_{GL}^{\Om } (u) \right)^{\frac 1p}
\qquad
\mbox{ for any } p \geq 2. 
\eeq

For a measurable set  $ \omega \subset \R^N$ with $ \Lo ^N ( \omega ) < \infty$ and for 
any $ f \in L^1( \omega)$, we denote by 
$m( f, \omega) = \frac{1}{\Lo ^N ( \omega ) } \ds \int_{\omega} f (x) \, dx $ 
the mean value of $f$ on $ \omega$.

Let $ x_0 $ be such that $B( x_0 , 4R) \subset \Om $. 
Using the Poincar\'e inequality and (\ref{3.1}) we have 
\beq
\label{3.16}
\| v_h - m( v_h, B( x_0, 4R)) \|_{L^2( B(x_0, 4R))} 
\leq C_P R \|\nabla v_h \|_{L^2( B(x_0, 4R))}  
\leq C_P R \left( E_{GL}^{\Om } (u) \right)^{\frac 12}.
\eeq
We claim that there exist $ k \in \N$, depending only on $N$, 
and $C_* = C_*( N, h, R)$ such that 
\beq
\label{3.17}
\| v_h - m( v_h, B( x_0, 4R)) \|_{W^{2, N} ( B(x_0, \frac{R}{2^{k-2}}))} 
\leq C _* \left( \left( E_{GL}^{\Om } (u) \right)^{\frac 12}
+ \left( E_{GL}^{\Om } (u) \right)^{\frac 1N} \right).
\eeq

It is well-known (see  Theorem 9.11 p. 235 in \cite{GT})
that for $ p \in (1, \infty )$ there exists $ C = C(N, r, p ) > 0$ such that for any 
$ w \in W^{2, p} ( B(a, 2r)) $  there holds
\beq
\label{3.18}
\|w\| _{W^{2, p}(B(a, r))} \leq C \left( 
\|w\| _{L^ p(B(a, 2 r))} + \|\Delta w\| _{L^ p(B(a, 2 r))} \right).
\eeq

From (\ref{3.15}), (\ref{3.16}) and (\ref{3.18}) we infer that 
\beq
\label{3.19}
\| v_h - m( v_h, B( x_0, 4R)) \|_{W^{2, 2} ( B(x_0, 2R ))} 
\leq C ( N, h, R)  \left( E_{GL}^{\Om } (u) \right)^{\frac 12}.
\eeq
If $\frac 12 - \frac 2N \leq \frac 1N$, from (\ref{3.19}) and the Sobolev embedding we find 
\beq
\label{3.20}
\| v_h - m( v_h, B( x_0, 4R)) \|_{L^N  ( B(x_0, 2R ))} 
\leq C ( N, h, R)  \left( E_{GL}^{\Om } (u) \right)^{\frac 12}.
\eeq
Then using (\ref{3.15}) (for $p = N$), (\ref{3.20}) and (\ref{3.18}) we infer that 
(\ref{3.17}) holds for $ k =2$. 

If $\frac 12 - \frac 2N >  \frac 1N$,  (\label({3.19}) and the Sobolev embedding imply 
\beq
\label{3.21}
\| v_h - m( v_h, B( x_0, 4R)) \|_{L^{p_1} ( B(x_0, 2R ))} 
\leq C ( N, h, R)  \left( E_{GL}^{\Om } (u) \right)^{\frac 12}, 
\eeq
where $ \frac{1}{p_1} = \frac 12 - \frac 2N$. 
Then (\ref{3.21}), (\ref{3.15}) and (\ref{3.18}) give 
\beq
\label{3.22}
\| v_h - m( v_h, B( x_0, 4R)) \|_{W^{2, p_1} ( B(x_0, R ))} 
\leq C ( N, h, R)  \left( \left( E_{GL}^{\Om } (u) \right)^{\frac 12}
+ \left( E_{GL}^{\Om } (u) \right)^{\frac 1N} \right).
\eeq
If $\frac{1}{p_1}  - \frac 2N \leq \frac 1N$,  
using (\ref{3.22}), the Sobolev embedding, (\ref{3.15}) and (\ref{3.18}) 
we get 
$$
\| v_h - m( v_h, B( x_0, 4R)) \|_{W^{2, N} ( B(x_0, \frac R2 ))} 
\leq C ( N, h, R)  \left( \left( E_{GL}^{\Om } (u) \right)^{\frac 12}
+ \left( E_{GL}^{\Om } (u) \right)^{\frac 1N} \right);
$$
otherwise we repeat the  process. 
After a finite number of steps we find $ k \in \N$ such that (\ref{3.17}) holds.

We will use the following variant of the Gagliardo-Nirenberg inequality: 
\beq
\label{3.23}
\| w - m( w,B ( a, r) ) \|_ {L^p ( B(a, r))} 
\leq C(p, q, N, r) \|w\| _ {L^q ( B(a,2 r))} ^{\frac qp} 
\|\nabla w\| _ {L^N ( B(a,2 r))} ^{1- \frac qp} 
\eeq
for any $ w \in W^{1, N}(B(a, 2r))$, where $ 1 \leq q \leq p < \infty $ 
(see, e.g., \cite{kavian} p. 78). 

Using   (\ref{3.23}) with $ w = \nabla v_h$ and (\ref{3.17}) we find
\beq
\label{3.24}
\begin{array}{l}
\|\nabla v_h - m( \nabla v_h, B(x_0, \frac{R}{2^{k-1}}))\|
_{L^p(B(x_0, \frac{R}{2^{k-1}}))}
\\
\leq C 
\|\nabla v_h \| _{L^2(B(x_0, \frac{R}{2^{k-2}}))}^{\frac 2p}
\|\nabla ^2 v_h \| _{L^N(B(x_0, \frac{R}{2^{k-2}}))}^{1-\frac 2p}
\\
\leq C \left( E_{GL}^{\Om } (u) \right)^{\frac 1p}
\left( \left( E_{GL}^{\Om } (u) \right)^{\frac 12}
+ \left( E_{GL}^{\Om } (u) \right)^{\frac 1N} \right)^{1-\frac 2p}
\end{array}
\eeq
for any $ p \in [2, \infty)$, where the constants depend only on 
$ N, \, p, \, h, \, R$. 

Using the Cauchy-Schwarz inequality and (\ref{3.1}) we have 
$$
\Big\vert m( \nabla v_h, B(x_0, \frac{R}{2^{k-1}})) \Big\vert 
\leq
\Lo ^N ( B(x_0, \frac{R}{2^{k-1}}) )^{- \frac 12} \|\nabla v _h \|_{L^2(B(x_0, \frac{R}{2^{k-1}}))}
\leq C  \left( E_{GL}^{\Om } (u) \right)^{\frac 12}
$$
and we infer that for any $ p \in [1, \infty] $  the following estimate  holds:
\beq
\label{3.25}
\begin{array}{l}
\|  m( \nabla v_h, B(x_0, \frac{R}{2^{k-1}})) \|_{L^p(B(x_0, \frac{R}{2^{k-1}}))}
\\
\leq
\Big\vert m( \nabla v_h, B(x_0, \frac{R}{2^{k-1}})) \Big\vert 
\left(\Lo ^N (B(x_0, \frac{R}{2^{k-1}})) \right)^{\frac 1p}
\leq C(N, p, R) \left( E_{GL}^{\Om } (u) \right)^{\frac 12}.
\end{array}
\eeq
From (\ref{3.24}) and (\ref{3.25}) we obtain for any $ p \in [2, \infty)$, 
\beq
\label{3.26}
\|\nabla v_h \|_{L^p(B(x_0, \frac{R}{2^{k-1}}))}
\leq C( N, p, h, R) \left( 
\left( E_{GL}^{\Om } (u) \right)^{\frac 12} 
+ \left( E_{GL}^{\Om } (u) \right)^{\frac 1p + \frac 1N ( 1 - \frac 2p)} \right).
\eeq

We will use the Morrey inequality which asserts that for any 
$ w \in C^0 \cap W^{1, p}(B(x_0, r))$ with $ p >N$ there holds
\beq
\label{3.27}
| w(x) - w(y) | \leq C(p, N) | x-y|^{1 - \frac Np} \|\nabla w \|_{L^p(B(x_0, r))}
\qquad \mbox{ for all } x, y \in B(x_0, r)
\eeq
(see, e.g., the proof of Theorem IX.12 p. 166 in \cite{brezis}). 
Using (\ref{3.26}) and the Morrey's inequality (\ref{3.27}) for $ p = 2N$
we get 
\beq 
\label{3.28}
|v_h (x) - v_h(y)| \leq C( N, h, R) |x - y |^{\frac 12 } 
\left(
\left( E_{GL}^{\Om } (u) \right)^{\frac 12} 
+ \left( E_{GL}^{\Om } (u) \right)^{ \frac 1N ( 1 + \frac{1}{2^*})}
\right)
\eeq
 for any $ x, y \in B(x_0, \frac{R}{2^{k-1}})$.

Let $ \de > 0$. Assume that there is $ x_0 \in \Om $ such that $B(x_0, 4R) \subset \Om $ and 
$\big| \, | v_h (x_0) +1| - 1 \big| \geq \de $.
Since 
$\big\vert \,  | \, | v_h (x) +1| - 1 | - |\,  | v_h (y) +1| - 1 | \,  \big\vert 
\leq | v_h (x) - v_h (y)|$, 
from (\ref{3.28}) we infer that 
$$
 \big| \, | v_h (x) +1| - 1 \big| \geq \frac{\de}{2} 
 \qquad 
\mbox{ for any } x \in B( x_0,  r_{\de}) , 
$$
 where
\beq
\label{3.29}  
r_{\de } = \min \left( \frac{R}{2^{k-1}}, 
\left(\frac{\de}{2 C(  N, h, R)} \right)^2 
\left( \left( E_{GL}^{\Om } (u) \right)^{\frac 12} 
+ \left( E_{GL}^{\Om } (u) \right)^{ \frac 1N ( 1 + \frac{1}{2^*})} \right)^{-2}
\right).
\eeq

Let 
\beq
\label{3.30}
\eta ( s) = \inf \{ ( \ph ^2 ( \tau ) - 1 )^2 \; | \; 
\tau \in ( - \infty, 1 - s ] \cup [1 + s, \infty) \}.
\eeq
It is clear that $ \eta $ is nondecreasing and positive on $(0, \infty )$. 
We have:
\beq
\label{3.31}
\begin{array}{l}
E_{GL}^{\Om} (u) \geq E_{GL}^{\Om} (v_h ) 
\geq \ds \frac 12 { \int _{B(x_0, r_{\de})} }
\left( \ph ^2(|1 + v_h | ) - 1 \right)^2 \, dx
\\
\geq \ds \frac 12 { \int _{B(x_0, r_{\de})} } \eta \left( \frac{\de}{2} \right) \, dx 
=  \frac 12 \Lo ^N (B(0, 1))  \eta ( \frac{\de}{2}) r_{\de }^N, 
\end{array}
\eeq
where $ r_{\de}$ is given by (\ref{3.29}). It is obvious that there exists a constant 
$K >0$, depending only on $  N, \,  h, \, R, \, \de $ 
such that (\ref{3.31}) cannot hold for $E_{GL}^{\Om } (u) \leq K$. 
We infer that 
$ \big| \, | 1+ v_h (x_0) | - 1 \big| < \de $ 
if $B(x_0, 4R) \subset \Om $ and  $E_{GL}^{\Om } (u) \leq K$.
This completes the 
proof of  Lemma \ref{L3.1}.
\hfill
$\Box$

\begin{Lemma}
\label{vanishing}
Let $ (u_n)_{n \geq 1} \subset \Xo $ be a sequence of functions satisfying: 

\smallskip

a) $E_{GL}(u_n)$ is bounded and

\smallskip

b) $ \ds \lim_{n \ra \infty } \Big( \sup_{y \in \R^N} E_{GL}^{B(y, 1)} (u_n) \Big) =0.$

\smallskip

There exists a sequence $ h_n \lra 0 $ such that for any minimizer $ v_n $ of 
$G_{h_n, \R^N}^{u_n}$ in $H_{u_n} ^1 (\R^N)$ we have 
$ \| \, | v_n +1 | - 1 \|_{L^{\infty }(\R^N)} \lra 0 $ as $ n \lra \infty $.
\end{Lemma}

{\it Proof. } The proof of Lemma \ref{vanishing} is quite tricky 
and we split it into four steps. 
First we explain the choice of the sequence $(h_n)_{n \geq 1}$.
Then we prove that there 
is $ C > 0 $ such that 
for any minimizer $ v_n $ of $G_{h_n, \R^N} ^{u_n} $ and for all $ x \in \R^N$ there holds 
$\| \Delta v_n \|_{L^N (B(x,1)) } \leq C$. 
To get this estimate we write (\ref{3.5}) in a convenient form, multiply it by appropriate cut-off functions, 
then perform integrations by parts and use elliptic regularity and a finite induction 
to prove that $u_n$ and $ v_n $ are locally sufficiently close 
(for instance, it follows from (\ref{3.36}) and (\ref{3.51}) below that 
$\| u_n - v_n \| _{L^2(B(x, 1))} \leq C h_n ^N$ for all $ x \in \R^N$). 
 Then we use (\ref{3.5}) again to get the desired bound on $ \Delta v_n$. 
In the third step we use Sobolev and Morrey inequalities to prove that $ v_n $ is uniformly  H\"older continuous. 
Finally, if $ \de $ is fixed and $ \big| \, | 1 + v_n ( x_0) | - 1 \big| \geq \de $ for some $ x_0 \in \R^N$, we  have necessarily
 $ \big| \, | 1 + v_n  | - 1 \big| \geq \frac{\de }{2}$ on a ball $B(x_0 , r)$, where $ r$ does not depend on $ n$, 
 thus $\| \ph ^2 (|1 + v_n|) - 1 \| _{L^2(B(x_0, 1))} $ is bounded from below by a positive constant. 
 This is impossible for large $n$ because 
$\| \ph ^2 (|1 + v_n|) - 1 \| _{L^2(B(x_0, 1))} $ is close to $\| \ph ^2 (|1 + u_n|) - 1 \| _{L^2(B(x_0, 1))} $
and the last quantity tends to zero by assumption (b).

\medskip

{\it Step 1. Choice of  $h_n$. }
Let $M = \ds \sup_{n \geq 1} E_{GL}(u_n)$. 
For $ n \geq 1 $ and $ x \in \R^N$ we denote
$$
m_n (x) = m(u_n, B(x, 1)) = \frac{1}{\Lo ^N (B(0,1))} \int_{B(x, 1)} u_n (y) \, dy. 
$$
By the Poincar\'e inequality, there exists $ C_0 > 0$ such that 
$$
\int_{B(x, 1)} | u_n (y) - m_n (x)| ^2  \, dy \
\leq C_0 \int_{B(x, 1)} | \nabla u_n (y) | ^2 \, dy. 
$$
From (b) and the Poincar\'e inequality it follows that 
\beq
\label{3.32}
\sup_{x \in \R^N} \| u_n - m_n (x) \|_{L^2(B(x, 1))} \lra 0 \qquad \mbox{ as } n \lra \infty .
\eeq
Let $H$ be as in Lemma \ref{3.1} (iii).
From (\ref{3.12}) and (b) we get 
\beq
\label{3.33}
\sup_{x \in \R^N} \| H(u_n) \| _{L^2(B(x, 1))} ^2 
\leq \sup_{x \in \R^N} 9  \int_{B(x, 1)} 
\left( \ph ^2(| 1+ u_n(y)|) - 1 \right)^2 \, dy
\lra 0 
\eeq
 as $n \lra \infty$. 
It is obvious that $H$ is Lipschitz on $ \C$. 
Using (\ref{3.32}) we find
\beq
\label{3.34}
\sup_{x \in \R^N} \| H(u_n) - H(m_n(x))\| _{L^2(B(x, 1))} 
\leq C_1 \sup_{x \in \R^N} \| u_n - m_n(x) \| _{L^2(B(x, 1))} 
\lra 0 
\eeq
as  $ n \lra \infty$. 
From (\ref{3.33}) and (\ref{3.34})  we infer that 
$
\sup_{x \in \R^N} \| H(m_n(x))\| _{L^2(B(x, 1))}   \lra 0 
$ as $ n \lra \infty.$ 
Since $\|H(m_n(x))\| _{L^2(B(x, 1))} = \left( \Lo ^N( B(0,1)) \right)^{\frac 12}  | H( m_n(x))|$, we have proved that 
\beq
\label{3.35}
\lim_{ n \ra \infty} \sup_{x \in \R^N} |H(m_n(x))| = 0.
\eeq
Let 
\beq
\label{3.36}
h_n = \max \left( \left(  \sup_{x \in \R^N} \| u_n - m_n (x) \|_{L^2(B(x, 1))} \right)^{\frac{1}{N+2}}, 
\left( \sup_{x \in \R^N} |H(m_n(x))| \right)^{\frac 1N} \right). 
\eeq
From (\ref{3.32}) and (\ref{3.35}) it follows that 
$ h_n \lra 0 $ as $ n \lra \infty $. 
Thus we may assume that $ 0 < h_n < 1$ for any $n$ 
(if $ h_n = 0$, we see that $ u_n$ is constant  a.e. 
and there is nothing to prove). 

Let $ v_n $ be a minimizer of $G_{h_n, \R^N}^{u_n}$ 
(such  minimizers exist by Lemma \ref{L3.1} (i)).
It follows from Lemma \ref{L3.1} (iii) that $v_n $ satisfies (\ref{3.5}). 

\medskip

{\it Step 2. 
We  prove that there exist $R_N > 0$ and $C> 0$,  independent on $n$, such that }
\beq
\label{3.37}
\| \Delta v_n \|_{L^N(B(x, R_N)) }   \leq C \qquad \mbox{ for any } x \in \R^N 
\mbox{ and } n \in \N^*. 
\eeq 
Clearly, it suffices to prove (\ref{3.37}) for $ x = 0$. 
Let $ m_n = m_n (0) $.  Then (\ref{3.5}) can be written as 
\beq
\label{3.38}
- \Delta v _n + \frac{1}{h_n ^2} {\ph} ' ( | v_n - m_n | ^2) (v_n - m_n)
= f_n, 
\eeq
where 
\beq
\label{3.39} 
\begin{array}{rcl}
f_n & = & -  \left( H( v_n) - H( m_n) \right) -  H( m_n) 
\\
\\
& &  + \frac{1}{ h_n ^2} 
\left( {\ph} ' ( | v_n - m_n | ^2) (v_n - m_n) 
- {\ph} ' ( | v_n - u_n | ^2) (v_n - u_n)  \right).
\end{array}
\eeq
In view of Lemma \ref{L3.1} (iii),  equality  (\ref{3.38}) 
holds in $L_{loc}^p ( \R^N) $ (and not only in $ \Do ' ( \R^N)$). 

The function $ z \longmapsto {\ph } ' ( |z|^2) z $ belongs to
 $ C_c^{\infty }(\C)$ 
and consequently it is Lipschitz. 
Using (\ref{3.36}), we see that there exists $C_2 > 0 $ such that 
\beq
\label{3.40}
\begin{array}{l}
\|  {\ph} ' ( | v_n - m_n | ^2) (v_n - m_n) 
- {\ph} ' ( | v_n - u_n | ^2) (v_n - u_n) \|_{L^2(B(0, 1))} 
\\
\\
\leq C_2 \| u_n - m_n \| _{L^2(B(0, 1))}  
\leq C_2 h_n ^{ N+2}. 
\end{array}
\eeq
By (\ref{3.36}) we have also 
$ \|H(m_n)\| _{L^2(B(0, 1))}  = \left( \Lo ^N( B(0, 1) \right)^{\frac 12} |H(m_n)|
\leq \left( \Lo ^N( B(0, 1) ) \right)^{\frac 12} h_n ^N$. 
From this estimate, (\ref{3.39}), (\ref{3.40}) and the fact that $H$ is Lipschitz we get 
\beq
\label{3.41}
\| f_n \| _{L^2(B(0, R))}
\leq C_3 \|v_n - m_n \|_{L^2(B(0, R))} + C_4 h_n ^N 
\qquad \mbox{ for any } R \in (0, 1].
\eeq
Let $ \chi \in C_c^{\infty}( \R^N, \R)$. 
Taking the scalar product (in $ \C$) of (\ref{3.38}) by 
$ \chi (x) ( v_n (x) - m_n)$ and integrating by parts we find 
\beq
\label{3.42}
\begin{array}{l}
{ \ds \int_{\R^N} } \chi |\nabla v_n |^2 \, dx 
+ \ds \frac{1}{h_n ^2}  { \ds \int_{\R^N} } 
\chi {\ph } ' ( |v_n - m_n|^2) |v_n - m_n|^2 \, dx 
\\
\\ 
= \ds \frac 12 { \ds \int_{\R^N} } (\Delta \chi ) |v_n - m_n|^2 \, dx 
+ { \ds \int_{\R^N} } \langle f_n (x) , v_n (x) - m_n \rangle \chi(x) \, dx.
\end{array}
\eeq
From (\ref{3.2}) we have $\|v_n - u_n \|_{L^2( \R^N ) } \leq C _5 h_n ^{\frac 2N}$, thus
\beq
\label{3.43}
\| v_n - m_n \|_{L^2 (B(0, 1))} 
\leq \| v_n - u_n \|_{L^2 (B(0, 1))} 
+ \| u_n - m_n \|_{L^2 (B(0, 1)) } \leq K_0 h_n^{\frac 2N}.
\eeq
We prove that 
\beq
\label{3.44}
\| v_n - m_n \|_{L^2 (B(0, \frac{1}{2^{j-1}}))} \leq K_j h_n ^{\frac{2j}{N}}
\qquad \mbox{ for } 1\leq j \leq \left[ \frac{N^2}{2} \right] +1, 
\eeq
where $K_j$ does not depend on $n$. 
We proceed by induction. From (\ref{3.43}) it follows that (\ref{3.44})
is true for $ j =1$. 

Assume that (\ref{3.44}) holds for some 
$ j \in \N^*$, $ j \leq \left[ \frac{N^2}{2} \right]$.
Let $ \chi_j \in C_c^{\infty}(\R^N)$ be a real-valued function such that 
$ 0 \leq \chi_j \leq 1$, $\mbox{supp}(\chi_j) \subset B(0, \frac{1}{2^{j-1}})$
and $ \chi _j = 1 $ on $B(0, \frac{1}{2^{j}})$.
Replacing $ \chi $ by $ \chi _j $ in (\ref{3.42}), then using 
the Cauchy-Schwarz inequality and (\ref{3.41}) we find
\beq
\label{3.45}
\begin{array}{l}
\ds \int_ {B(0, \frac{1}{2^{j}})} |\nabla v_n |^2 \, dx 
+ \frac{1}{h_n ^2} \ds \int_ {B(0, \frac{1}{2^{j}})}  
{\ph }' ( |v_n - m_n |^2) |v_n - m_n |^2 \, dx 
\\
\\
\leq \frac 12 \|\Delta \chi _j \|_{L^{\infty }(\R^N)} 
\| v_n - m_n \|_{L^2 ( B(0, \frac{1}{2^{j-1}}))} ^2 
+ \| f_n \|_{L^2 (B(0, \frac{1}{2^{j-1}}))} 
\| v_n - m_n \|_{L^2 (B(0, \frac{1}{2^{j-1}}))} 
\\
\leq A_j \| v_n - m_n \|_{L^2 (B(0, \frac{1}{2^{j-1}}))}  ^2 
+ C_4 h_n^N \| v_n - m_n \|_{L^2 (B(0, \frac{1}{2^{j-1}}))} 
\leq A_j ' h_n ^{\frac{4j}{N}}. 
\end{array}
\eeq
From (\ref{3.44}) and (\ref{3.45}) we infer that 
$\| v_n - m_n \|_{H^1 ( B(0, \frac{1}{2^{j}}) )} \leq B_j h_n ^{\frac{2j}{N}}$.
Then the Sobolev embedding implies
\beq
\label{3.46}
\| v_n - m_n \|_{L^{2^*} ( B(0, \frac{1}{2^{j}}) ) }
\leq D_j h_n ^{\frac{2j}{N}}.
\eeq
The function $ z \longmapsto {\ph } (|z|^2) $ is clearly Lipschitz 
on $ \C$, thus we have 
$$
\begin{array}{l}
\ds \int_{B(0,1)} 
| {\ph } (|v_n - u_n |^2) - {\ph } (|v_n - m_n |^2) | \, dx 
\leq C_6' \ds \int_{B(0,1)}  |u_n - m_n | \, dx 
\\
\leq C_6 \| u_n - m_n \|_{L^2 (B(0, 1))}
\leq C_6 h_n ^{N+2}. 
\end{array}
$$
It is clear that 
$ \ds \int_{B(0,1)} 
 {\ph } (|v_n - u_n |^2)  \, dx \leq h_n ^2 G_{h_n , \R^N}^{u_n} ( v_n) 
 \leq h_n ^2 E_{GL}(u_n) \leq h_n ^2 M$ and we obtain 
\beq
\label{3.47}
\ds \int_{B(0,1)} 
 {\ph } (|v_n - m_n |^2) \, dx \leq C_7 h_n ^2 .
\eeq
If $ | v_n (x) - m_n| \geq \sqrt{2} $ we have 
$ {\ph } (|v_n (x) - m_n |^2)  \geq 2 $, 
hence
\beq
\label{3.48}
 \Lo ^N ( \{ x \in B(0, 1) \; | \; | v_n (x) - m_n| \geq \sqrt{2}  \} ) 
\leq
\ds \frac 12 { \int_{B(0,1)}  } {\ph } \left( |v_n - m_n |^2 \right) 
\, dx 
\leq \frac{C_7}{2} h_n ^2 .
\eeq
By H\"older's inequality, (\ref{3.46}) and (\ref{3.48}) we have
\beq
\label{3.49}
\begin{array}{l}
\ds \int _{\{ |v_n - m_n | \geq \sqrt{2}   \} \cap B(0, \frac{1}{2 ^j} )}
|v_n - m_n |^2 \, dx 
\\
\leq
\| v_n - m_n \|_{L^{2^*} ( B(0, \frac{1}{2^{j}}) ) } ^2 
\left( \Lo ^N ( \{ x \in B(0, 1) \; | \; | v_n (x) - m_n| \geq \sqrt{2} \} ) \right)^{1 - \frac{2}{2^*}}
\\
\leq \left( D_j h_n ^{\frac{2j}{N}}\right)^2 
\left( C_7 h_n ^2 \right)^{1 - \frac{2}{2^*}}
\leq E_j h_n ^{\frac{4j+4}{N}}.
\end{array}
\eeq
From (\ref{3.45}) it follows that 
\beq
\label{3.50}
\begin{array}{l}
\ds \int _{\{ |v_n - m_n | < \sqrt{2}  \} \cap B(0, \frac{1}{2 ^j} )}
|v_n - m_n |^2 \, dx 
\leq 
\ds \int _{ B(0, \frac{1}{2 ^j} )}
{ \ph }' (|v_n - m_n | ^2) |v_n - m_n | ^2 \, dx
\\
\\
\leq A_j 'h_n ^{ 2 + \frac{4j}{N}} \leq A_j 'h_n ^{ \frac{4j+4}{N}}.
\end{array}
\eeq
Then (\ref{3.49}) and (\ref{3.50}) imply that (\ref{3.44}) 
holds for $ j +1$ and the induction is complete. Thus (\ref{3.44})  is established. 
Denoting $ j_N = \left[ \frac {N^2}{2} \right] + 1$ and 
$R_N = \frac{1}{2^{ j_N -1}}$, we have proved that 
\beq
\label{3.51}
\| v_n - m_n \|_{L^2( B(0, R_N ) ) }
\leq K_{j_N} h_n ^{\frac{2 j_N}{N}} \leq K_{j_N} h_n ^N. 
\eeq
It follows that 
\beq
\label{3.52}
\begin{array}{l}
\ds \int_{B(0, R_N)} \Big\vert \frac{1}{h_n^2} 
{\ph} ' ( | v_n - m_n |^2) ( v_n - m_n ) \Big\vert ^N \, dx
\\
\leq \ds \frac{1}{h_n ^{2N}} \sup_{z \in \C } \big\vert {\ph}'\left(|z|^2 \right) z \big\vert^{N-2}
 \ds \int_{B(0, R_N)} | v_n - m_n |^2 \, dx 
 \leq C_8. 
\end{array}
\eeq
Arguing as in (\ref{3.40}) and using (\ref{3.36}) we get 
\beq
\label{3.53}
\begin{array}{l}
\|  {\ph} ' ( | v_n - m_n | ^2) (v_n - m_n) 
- {\ph} ' ( | v_n - u_n | ^2) (v_n - u_n) \|_{L^N(B(0, 1))} ^N
\\
\\
\leq C_9  
\ds \sup_{z \in \C } \big\vert {\ph}'\left(|z|^2 \right) z \big\vert^{N-2}
\| u_n - m_n \| _{L^2(B(0, 1) )}  ^2
\leq C_{10} h_n ^{2 N+4}. 
\end{array}
\eeq
From (\ref{3.39}), (\ref{3.53}) and the fact that $H$ 
is  bounded  on $ \C$ it follows that $\|f_n \|_{L^N ( B(0, R_N ))} \leq C_{11}$, 
where $C_{11}$ is  independent of $n$. 
Using this estimate, (\ref{3.52}) and (\ref{3.38}), 
we infer that (\ref{3.37}) holds.

Since any ball of radius $1$ can be covered by a finite number of balls of radius
$ R_N$, it follows that there exists $ C>0$ such that 
\beq
\label{3.54}
\|\Delta v_n \|_{L^N(B(x, 1)) }   \leq C \qquad \mbox{ for any } x \in \R^N 
\mbox{ and } n \in \N^*. 
\eeq 

\medskip
{\it Step 3. The functions $v_n$ are uniformly H\"older continuous. } 
We will use (\ref{3.18}) and (\ref{3.54}) to prove that there exist
$ \tilde{R}_N \in ( 0, 1] $ and $C>0$ such that 
\beq
\label{3.55}
\| v_n - m_n ( x ) \|_{W^{2,N} (B(x, \tilde{R}_N) ) }  
 \leq C \qquad \mbox{ for any } x \in \R^N 
\mbox{ and } n \in \N^*. 
\eeq  
As previously, it suffices to prove (\ref{3.55}) for $ x_0 = 0$. 
From (\ref{3.54}) and H\"older's inequality it follows that for 
$1 \leq p \leq N$ we have
\beq
\label{3.56}
\|\Delta v_n \|_{L^p(B(x, 1)) }   
\leq \left( \Lo ^N (B(0,1)) \right)^{\frac 1p  - \frac 1N} 
\|\Delta v_n \|_{L^N(B(x, 1)) } 
\leq C(p). 
\eeq
Using (\ref{3.43}), (\ref{3.56}) with $ p =2 $ and (\ref{3.18}) we obtain 
\beq
\label{3.57}
\| v_n - m_n (0) \|_{W^{2,2} (B(x, \frac 12) )}  
 \leq C  .
\eeq
If $ \frac 12 - \frac 2N \leq \frac 1N$, (\ref{3.57}) and the Sobolev 
embedding give 
$$
\| v_n - m_n (0) \|_{L^N (B(x, \frac 12) )}   \leq C, 
$$
and this estimate together with (\ref{3.54}) and (\ref{3.18}) imply that 
(\ref{3.55}) holds for $ \tilde{R}_N = \frac 14$. 

If $ \frac 12 - \frac 2N > \frac 1N$, from (\ref{3.57}) and the Sobolev 
embedding we find 
$ 
\| v_n - m_n (0) \|_{L^{p_1} (B(x, \frac 12)) }   \leq C, 
$
where $ \frac{1}{p_1} = \frac 12 - \frac 2N $. 
This estimate, (\ref{3.56}) and (\ref{3.18}) imply 
$ \| v_n - m_n (0) \|_{W^{2, p_1} (B(x, \frac 14) )} \leq C$. 
If $\frac{1}{p_1} - \frac 2N \leq \frac 1N$, from the Sobolev embedding we obtain 
$
\| v_n - m_n (0) \|_{L^N (B(x, \frac 14) )}   \leq C, 
$ 
and then using (\ref{3.54}) and (\ref{3.18}) we infer that 
(\ref{3.55}) holds for $ \tilde{R}_N = \frac 18$. 
Otherwise we 
repeat the above argument.  After a finite number of steps we see that 
(\ref{3.55}) holds. 

Next we proceed as in the proof of Lemma \ref{L3.1} (iv). 
By (\ref{3.23}) and (\ref{3.55}) we have for  $ p \in [2, \infty) $ 
and any $ x _0 \in \R^N$,  
\beq
\label{3.58}
\begin{array}{l}
\|\nabla v_n - m( \nabla v_n, B( x_0, \frac 12 \tilde{R}_N )) \|
_{L^p(B( x_0, \frac 12 \tilde{R}_N ))} 
\\
\leq C \|\nabla v_n \|_{L^2(B( x_0,  \tilde{R}_N ))} ^{\frac 2p}
\|\nabla ^2 v_n \|_{L^N(B( x_0,  \tilde{R}_N ))} ^{1 - \frac 2p}
\leq C_1(p).
\end{array}
\eeq
Arguing as in (\ref{3.25}) we see that 
$\| m( \nabla v_n, B( x_0, \frac 12 \tilde{R}_N )) \|
_{L^p(B( x_0, \frac 12 \tilde{R}_N ))} $ is bounded independently on $n$ and hence
$$
\|\nabla v_n  \|_{L^p(B( x_0, \frac 12 \tilde{R}_N ))}  \leq C_2(p)
\qquad \mbox{ for any }  n \in \N^* \mbox{ and } x_0 \in \R^N. 
$$
Using this estimate for $ p = 2N$ together with the Morrey inequality 
(\ref{3.27}), 
we see that there exists $C_* >0$ such that 
for any $ x, y \in \R^N $
 with $|x-y| \leq \frac{\tilde{R}_N}{2} $  and any $ n \in \N^* $ we have
\beq
\label{3.59}
|v_n (x) - v_n (y)| \leq C_* |x - y |^{\frac 12} .
\eeq

\medskip

{\it Step 4. Conclusion. } 
Let $ \de _n = \| \, | v_n +1 | - 1 \|_{L^{\infty} ( \R^N)} $ 
and choose $ x_n \in \R^N$ such that 
$\big| \, | v_n (x_n) + 1| - 1 \big| \geq \frac{\de _n}{2}$. 
From (\ref{3.59}) it follows that 
$\big| \, | v_n (x) +1 | - 1 \big| \geq \frac{\de _n}{4}$
for all $ x \in B( x_n, r_n)$, where 
$$
r_n = \min \left( \frac{\tilde{R}_N}{2}, \left(\frac{\de _n}{4 C_*} \right)^2 \right).
$$
Then we have
\beq
\label{3.60} 
\begin{array}{l}
\ds \int_{B(x_n, 1)}
\left( \ph ^2(|1+ v_n (y)|) - 1 \right)^2 \, dy 
\geq 
\ds \int_{B(x_n, r_n)}
\left( \ph ^2(|1+ v_n (y)|) - 1 \right)^2 \, dy  
\\
\geq
{\ds \int_{B(x_n, r_n)} } \eta \left( \frac{\de_n}{4} \right) \, dy 
= \Lo ^N ( B(0, 1) ) \eta \left( \frac{\de_n}{4} \right) r_n ^N,  
\end{array}
\eeq
where $ \eta $ is as in (\ref{3.30}).

On the other hand, 
the function $ z\longmapsto \left( \ph ^2(| 1+ z|) - 1 \right)^2 $ is 
Lipschitz on $\C$. 
Using this fact, the Cauchy-Schwarz inequality, (\ref{3.2}) and assumption 
(a) we get
$$
\begin{array}{l}
\ds \int_{B(x, 1)} \Big\vert 
\left( \ph ^2(| 1+  v_n (y)|) - 1 \right)^2
- \left( \ph ^2(| 1+ u_n (y) |) - 1 \right)^2 \Big\vert  \, dy 
\\
\leq C \ds \int_{B(x, 1)} | v_n (y) - u_n (y) | \, dy 
\leq C' \| v_n - u_n \|_{L^2 ( B(x, 1) )} 
\leq C' \| v_n - u_n \|_{L^2 ( \R^N) } 
\leq C'' h_n ^{ \frac 2N} .
\end{array}
$$
Then using assumption (b) we infer that 
\beq
\label{3.61}
\ds \sup_{ x \in \R^N} \ds \int_{B(x, 1)}
\left( \ph ^2(| 1+ v_n (y)|) - 1 \right)^2 \, dy 
\lra 0 \qquad \mbox{ as } n \lra \infty.
\eeq

From (\ref{3.60}) and (\ref{3.61}) we get 
$ {\ds \lim_{n \ra \infty }} \eta \left( \frac{\de_n}{4} \right) r_n ^N = 0$
and this clearly implies $ \ds \lim_{n \ra \infty }  \de _n = 0$. 
Lemma \ref{vanishing} is thus proven.
\hfill
$\Box $

\medskip

The next result is based on Lemma \ref{L3.1} and will be very useful 
in the next sections to 
prove the "concentration" of minimizing sequences.
For $0 < R_1 < R_2 $ we denote 
$ \Om _{R_1, R_2} = B(0, R_2) \setminus \ov{B}(0, R_1)$. 

\medskip

\begin{Lemma}
\label{splitting}
Let $ A > A_3 > A_2 >1$. 
There exist $ \e _0 = \e _0 (  N, A, A_2, A_3)>0$ and 
$C_i = C_i (  N, A, A_2, A_3)>0$ such that for any $ R \geq 1$, 
$ \e \in (0, \e _0)$ and $ u \in \Xo $  verifying
$E_{GL}^{\Om _{R, AR}} (u) \leq \e, $
there exist two functions $ u_1, \, u_2 \in \Xo $ and a constant $ \theta _0 \in [0, 2 \pi)$ 
satisfying the following properties:

\smallskip

i) $\mbox{\rm supp}(u_1) \subset B(0, A_2 R) $ and $ 1+ u_1  = e^{-i \theta _0 } ( 1+ u)$ 
on $B(0, R)$, 

\smallskip

ii) $ u_2 = u $ on $ \R^N \setminus B(0, AR)$ and 
$ 1+  u _2 =  e^{ i \theta _0 }= constant$
on $B(0, A_3R)$,  

\medskip

iii) $ \ds \int_{\R^N} \Big\vert  \, 
 \Big\vert \frac{\p u}{\p x_j } \Big\vert ^2 - 
 \Big\vert \frac{\p u_1}{\p x_j } \Big\vert ^2 - 
 \Big\vert \frac{\p u_2}{\p x_j } \Big\vert ^2 \, \Big\vert \, dx 
 \leq C_1 \e  $ for $ j =1, \dots, N$, 
 
\medskip
 
iv) $ \ds \int_{\R^N} \Big\vert  
\left( \ph ^2(|1+ u|) - 1 \right)^2 
- \left( \ph ^2(|1+ u_1|) - 1 \right)^2 
- \left( \ph ^2(|1+  u_2|) - 1 \right)^2 
\Big\vert   \, dx \leq C_2 \e $, 

\medskip

v) $ | Q(u) - Q(u_1) - Q( u_2) | \leq C_3 \e $, 

\medskip

vi) If assumptions (A1) and  (A2) in the introduction hold, then
$$
 \ds \int_{\R^N}  \Big\vert V( | 1+ u |^2) - V( |1 + u _1|^2) - V( |1 + u_2 |^2) \Big\vert \, dx 
\leq C_4 \e + C_5 \sqrt{\e } \left( E_{GL}(u) \right)^{\frac{2^* -1}{2}}.
$$ 

\end{Lemma}

{\it Proof.  } 
Fix $ k > 0$,  $ A_1 $ and  $  A_4 $ such that 
$ 1 + 4k < A_1 < A_2 < A_3 < A_4 < A - 4k. $
Let $ h = 1 $ and $ \de = \frac{1}{2}$. 
We will prove that Lemma \ref{splitting} holds for 
$\e _0 = K( N, h=1, \de = \frac{1}{2}, k )$, 
where $K( N, h, \de, R)$ is as in Lemma \ref{L3.1} (iv). 

Fix two functions   $ \eta _1, \eta _2 \in C ^{\infty }(\R)$ satisfying  the following properties:
$$
\begin{array}{l}
\eta_1 = 1 \mbox{ on } (-\infty, A_1], \quad \eta_1 = 0 \mbox{ on } [A_2, \infty), 
\quad \eta_1 \mbox{ is nonincreasing, }
\\
\eta_2 = 0 \mbox{ on } (-\infty, A_3], \quad \eta_2 = 1 \mbox{ on } [A_4, \infty), 
\quad \eta_2 \mbox{ is nondecreasing. }
\end{array}
$$

Let $ \e < \e _0 $ and let  $u \in \Xo $ be such that 
$ E_{GL}^{\Om _{R, AR}}(u) \leq \e$. 
Let $ v_1 $ be a minimizer of $G_{1, \Om_{R, AR} }^u$  in the space 
$H_u ^1( \Om_{R, AR})$. 
The existence of $ v_1 $ is guaranteed by Lemma \ref{L3.1}. 
We have $ v_1 = u $ on $ \R^N \setminus \Om _{R, AR} $ and by Lemma \ref{L3.1} (iii) we 
know that $ v_1 \in W_{loc}^{2, p} (\Om_{R, AR} )$ for any  $ p \in [1, \infty)$.
Moreover, since $E_{GL}^{\Om _{R, AR}}(u) \leq K( N, 1, \frac{r_0}{2}, k)$, 
Lemma \ref{L3.1} (iv)  implies that
\beq
\label{3.62}
\frac{1}{2} < | 1+v_1 (x) | < \frac{3}{2} \qquad
\mbox{ if } R + 4k \leq |x| \leq AR - 4k .
\eeq
Since $ N\geq 3$, $ \Om_{A_1R, A_4R}$ is simply connected and it follows directly 
from Theorem 3 p. 38 in \cite{BBM} that there exist two real-valued functions 
$ \rho, \, \theta \in W^{2, p}(\Om_{A_1R, A_4 R})$, $1 \leq p < \infty$,  such that 
\beq
\label{3.63}
1+ v_1 (x) = \rho (x) e^{ i \theta (x)} \qquad \mbox{ in } 
\Om_{A_1R, A_4 R}.
\eeq
For  $ j =1, \dots, N$  we have
\beq
\label{3.64}
\frac{\p v_1}{\p x_j} = \left(  \frac{ \p \rho }{\p x_j }
+ i \rho \frac{ \p \theta }{\p x_j } \right)e^{i \theta} 
\quad \mbox{ and } \quad 
\Big\vert 
\frac{\p v_1}{\p x_j} 
\Big\vert  ^2 
= 
\Big\vert 
\frac{\p \rho}{\p x_j} 
\Big\vert  ^2  
+ \rho ^2 \Big\vert 
\frac{\p \theta}{\p x_j} 
\Big\vert  ^2 
\quad \mbox{ a.e. in } \Om_{A_1R, A_4 R}.
\eeq
Thus we get the following estimates: 
\beq
\label{3.65}
\int_{\Om_{A_1 R, \, A_4R}}
|\nabla \rho |^2 \, dx \leq 
\int_{\Om_{A_1 R, \, A_4R}}
|\nabla v_1 |^2 \, dx \leq \e, 
\eeq
\beq
\label{3.66}
 \frac 12  \int_{\Om_{A_1 R, \, A_4R}} 
\left( \rho ^2 - 1 \right)^2 \, dx \leq E_{GL}^{\Om_{A_1 R, \, A_4R}} (v_1)
\leq \e, 
\eeq
\beq
\label{3.67}
\int_{\Om_{A_1 R, \, A_4R}}
|\nabla \theta |^2 \, dx \leq 
4 \int_{\Om_{A_1 R, \, A_4R}}
|\nabla v_1 |^2 \, dx \leq 4 \e. 
\eeq
The Poincar\'e inequality and a scaling argument imply that 
\beq
\label{3.67bis}
\int_{\Om_{A_1 R, \, A_4R}} 
|f - m(f , \Om_{A_1 R, \, A_4R} )|^2 \, dx 
\leq C (N, A_1, A_4) R^2 \int_{\Om_{A_1 R, \, A_4R}} |\nabla f|^2 \, dx 
\eeq
for any $ f \in H^1( \Om_{A_1 R, \, A_4R}) $, where $C (N, A_1, A_4)$ 
does not depend on $R$. 
Let $ \theta _0 = m(\theta, \Om_{A_1 R, \, A_4R})$. We may assume that 
$\theta_0 \in [0, 2 \pi)$ (otherwise we replace $ \theta $ by 
$ \theta - 2 \pi \left[ \frac{\theta_0}{2 \pi} \right]$). 
Using (\ref{3.67}) and (\ref{3.67bis})  we get 
\beq
\label{3.68}
\int_{\Om_{A_1 R, \, A_4R}}
|\theta - \theta _0 |^2 \, dx \leq 
C( N, A_1, A_4) R^2 \int_{\Om_{A_1 R, \, A_4R}}
|\nabla v_1 |^2 \, dx \leq C( N, A_1, A_4) R^2 \e. 
\eeq

We define $ \tilde{u}_1$ and $ u_2 $ by 
\beq
\label{3.69}
 \tilde{u}_1 (x) = \left\{ \begin{array}{l} 
 u(x) \quad \mbox{ if } x \in \ov{B}(0, R), 
\\
 v_1 (x) \quad \mbox{ if } x \in {B}(0, A_1R) \setminus \ov{B}(0, R), 
\\
\left(1 + \eta_1 (\frac{|x|}{R} ) (\rho (x) - 1 )\right) 
e^{i \left(\theta _0 + \eta_1 (\frac{|x|}{R} ) (\theta (x) - \theta_0 )\right)}  -1
\\
\qquad \qquad \qquad \qquad  \qquad  \quad \mbox{ if } x \in {B}(0, A_4R) \setminus B(0, A_1R),
\\
 e^{i \theta _0} -1 \quad \mbox{ if } x \in \R^N \setminus {B}(0, A_4R ), 
\end{array}
\right.
\eeq
\beq
\label{3.70}
 u_2 (x) = \left\{ \begin{array}{l} 
 e^{i \theta _0}  -1 \quad \mbox{ if } x \in  \ov{B}(0, A_1R ), 
\\
\left(1 + \eta_2 (\frac{|x|}{R} ) (\rho (x) - 1 )\right) 
e^{i \left(\theta _0 + \eta_2 (\frac{|x|}{R} ) (\theta (x) - \theta_0 )\right)}  -1
\\
\qquad \qquad \qquad \qquad \quad \qquad  \mbox{ if } x \in {B}(0, A_4R) \setminus \ov{B}(0, A_1R),
\\
 v_1 (x) \quad \mbox{ if } x \in {B}(0, AR) \setminus {B}(0, A_4 R),
\\
 u(x)  \quad \mbox{ if } x \in \R^N \setminus {B}(0, AR ), 
\end{array}
\right.
\eeq
then we define $ u_1 $ in such a way that 
$ 1+ u_1 = e^{- i \theta _0} ( 1+ \tilde{u}_1)$.
Since $ u \in \Xo $ and $ u - v_1 \in H_0^1(\Om_{R, \, AR})$, it is clear  that 
$u_1 \in H^1(\R^N)$, $ u_2 \in \Xo $ and (i), (ii) hold.

Since $ \rho + 1 \geq \frac 32 $ on $ \Om _{A_1R, \, A_4 R }$, 
from (\ref{3.66})   we get 
\beq
\label{3.71}
\| \rho - 1 \|_{L^2 ( \Om _{A_1R, \, A_4 R } ) } ^2
\leq \frac{8}{ 9} \e.
\eeq
Obviously, 
$$
\nabla \left(1 + \eta_i (\frac{|x|}{R} ) (\rho (x) - 1 )\right) 
= \frac 1R \eta_i ' (\frac{|x|}{R} ) (\rho (x) - 1 ) \frac{x}{|x|}
+ \eta_i (\frac{|x|}{R} ) \nabla \rho
$$ 
and using (\ref{3.65}), (\ref{3.71})  and the fact that $R\geq 1$ we get 
\beq
\label{3.72}
\begin{array}{l}
\| \nabla  \left(1 + \eta_i (\frac{|x|}{R} ) (\rho (x) - 1 )\right) 
\|_{L^2 ( \Om _{A_1R, \, A_4 R } ) } 
\\
\leq
\frac 1R \sup | \eta _i ' |\cdot  \|\rho - 1 \|_{L^2 ( \Om _{A_1R, \, A_4 R } ) } 
+ \|  \eta_i (\frac{|\cdot |}{R} ) \nabla \rho \|_{L^2 ( \Om _{A_1R, \, A_4 R } ) } 
\leq C \sqrt{\e }.
\end{array}
\eeq
Similarly, using (\ref{3.67}) and (\ref{3.68}) we find
\beq
\label{3.73}
\begin{array}{l}
\| \nabla  \left(\theta _0 + \eta_i (\frac{|x|}{R} ) (\theta (x) - \theta_0 )\right) 
\|_{L^2 ( \Om _{A_1R, \, A_4 R } ) } 
\\
\leq
\frac 1R \sup | \eta _i ' |\cdot  \|\theta - \theta_0 \|_{L^2 ( \Om _{A_1R, \, A_4 R } ) } 
+ \|  \eta_i (\frac{|\cdot |}{R} ) \nabla \theta \|_{L^2 ( \Om _{A_1R, \, A_4 R } ) } 
\leq C \sqrt{\e }.
\end{array}
\eeq
From (\ref{3.72}), (\ref{3.73}) and the definition of $u_1$, $u_2$ it follows that 
$\| \nabla u_i \|
_{L^2 ( \Om _{A_1R, \, A_4 R } ) } 
\leq C \sqrt{\e } $ for  $ i = 1, 2$. 
Therefore
$$
\begin{array}{l}
 \ds \int_{\R^N} \Big\vert  \, 
 \Big\vert \frac{\p u}{\p x_j } \Big\vert ^2 - 
 \Big\vert \frac{\p u_1}{\p x_j } \Big\vert ^2 - 
 \Big\vert \frac{\p u_2}{\p x_j } \Big\vert ^2 \, \Big\vert \, dx 
=  \ds \int_{\Om _{R, AR} } \Big\vert  \, 
 \Big\vert \frac{\p u}{\p x_j } \Big\vert ^2 - 
 \Big\vert \frac{\p u_1}{\p x_j } \Big\vert ^2 - 
 \Big\vert \frac{\p u_2}{\p x_j } \Big\vert ^2 \, \Big\vert \, dx 
\\
\\
\leq 
 \ds \int_{\Om _{R, A_1 R} \cup \Om _{A_4R, A R} } 
 \Big\vert \frac{\p u}{\p x_j } \Big\vert ^2 + 
 \Big\vert \frac{\p v_1}{\p x_j } \Big\vert ^2  dx 
+ 
  \ds \int_{\Om _{A_1R, A_4R} } 
 \Big\vert \frac{\p u}{\p x_j } \Big\vert ^2 + 
 \Big\vert \frac{\p u_1}{\p x_j } \Big\vert ^2  +
 \Big\vert \frac{\p u_2}{\p x_j } \Big\vert ^2 \, dx 
\leq C_1 \e
\end{array}
$$
and (iii) is proven.

 On $\Om _{A_1R, A_4R} $ we have
$ \rho \in [\frac{1}{2}, \frac{3}{2}]$, hence
$\ph \left(1 + \eta_i (\frac{|x|}{R} ) (\rho (x) - 1 )\right) 
= 1 + \eta_i (\frac{|x|}{R} ) (\rho (x) - 1)$ and
\beq
\label{3.74} 
\begin{array}{l}
\left( \! \ph ^2 \! \left(1 \! + \eta_i (\frac{|x|}{R} ) (\rho (x) - 1 )\right) 
- 1 \right)^2 \! \! 
= (\rho (x)- 1 )^2 \eta_i ^2 (\frac{|x|}{R} ) 
\left(\! 2 \! + \eta_i (\frac{|x|}{R} )( \rho (x)- 1)\!  \right)^2
\\
\\
\leq
\frac{25}{4} ( \rho (x)- 1) ^2 \leq \frac{25}{8} | \rho (x)-1 |.
\end{array}
\eeq
From (\ref{3.69})$-$(\ref{3.71}) and (\ref{3.74}) it follows that 
$\|  \ph ^2( |1+  u_i |) - 1 \|_{L^2 ( \Om _{A_1R, \, A_4 R } ) } 
\leq C\sqrt{ \e}$.
As above, we get 
$$
\begin{array}{l}
\ds \int_{\R^N} \Big\vert  
\left( \ph ^2(|1+ u|) - 1 \right)^2 
- \left( \ph ^2(|1+ u_1|) - 1 \right)^2 
- \left( \ph ^2(|1+ u_2|) - 1 \right)^2 
\Big\vert  \, dx
\\
\\ 
\leq 
 \ds \int_{\Om _{R, A_1 R} \cup \Om _{A_4R, A R} } 
\left( \ph ^2(|1+ u|) - 1 \right)^2  + 
\left( \ph ^2(|1+  v_1|) - 1 \right)^2  \, dx 
\\
\\
+ 
 \ds \int_{\Om _{A_1R, A_4R} }  
\left( \ph ^2(|1+ u|) - 1 \right)^2 
+ \left( \ph ^2(|1+ u_1|) - 1 \right)^2 
+ \left( \ph ^2(|1+ u_2|) - 1 \right)^2  \, dx 
\leq C_2 \e.
\end{array}
$$
This proves (iv). 

\medskip

Next we  prove (v).
Since 
$\langle i\frac{\p \tilde{u}_1}{\p x _1} , \tilde{u}_1 \rangle $
has compact support, a simple computation gives 
\beq
\label{3.74bis}
\ds Q(u_1) = L\left( \langle i\frac{\p u_1}{\p x _1} , u_1 \rangle \right)
=  L \left( \langle ie ^{-i \theta _0} \frac{\p \tilde{u}_1}{\p x _1}, 
    e ^{-i \theta _0} -1 + e ^{-i \theta _0}\tilde{u}_1 \rangle \right)
= \ds \int_{\R^N} 
 \langle i\frac{\p \tilde{u}_1}{\p x _1} , \tilde{u}_1 \rangle \, dx.
\eeq
From the definition of $\tilde{u}_1$ and $ u_2$ and the fact that 
$ u = v_1 $ on $ \R^N \setminus \Om _{R, AR}$ we get 
$$
 \langle i\frac{\p v_1}{\p x _1} , v_1 \rangle 
- \langle i\frac{\p \tilde{u}_1}{\p x _1} , \tilde{u}_1 \rangle
- \langle i\frac{\p u_2}{\p x _1} , u _2 \rangle
= 0 \qquad \mbox{  a.e.  on } \R^N \setminus \Om _{A_1 R, \, A_4R}.
$$
Using this identity, Definition \ref{D2.4}, (\ref{3.74bis}), then 
(\ref{2.3}) and (\ref{3.69}), (\ref{3.70}) we obtain 
\beq
\label{3.77}
\begin{array}{l}
Q(v_1 ) - Q(u_1) - Q(u_2) = 
\ds \int _{\Om _{A_1 R, \, A_4R} } 
 \langle i\frac{\p v_1}{\p x _1} , v_1 \rangle 
- \langle i\frac{\p \tilde{u}_1}{\p x _1} , \tilde{u}_1 \rangle
- \langle i\frac{\p u_2}{\p x _1} , u _2 \rangle \, dx 
\\
\\
= 
\ds \int _{\Om _{A_1 R, \, A_4R} }  

\I \left( \frac{\p v_1}{\p x _1} -  \frac{\p \tilde{u}_1}{\p x _1} 
   - \frac{\p u_2}{\p x _1} \right)\, dx 
- 
\ds \int _{\Om _{A_1 R, \, A_4R} }   
(\rho^2 - 1) \frac{ \p \theta }{ \p x_1} \, dx 
\\
\\
+ \ds \int _{\Om _{A_1 R, \, A_4R} }  
\sum_{i=1}^2
\left( \! \left( \! 1 + \eta_i (\frac{|x|}{R} ) (\rho  - 1 )\right) ^2 \! 
 - 1 \! \right)
\frac{\p }{\p x_1} \!
\left(\! \theta _0 + \eta_i (\frac{|x|}{R} ) (\theta  - \theta_0 ) \! \right)  
 dx
\\
\\
 -  \ds \int _{\Om _{A_1 R, \, A_4R} }      \frac{ \p \theta }{ \p x_1} 
- \sum_{i=1}^2 \frac{\p }{\p x_1} 
\left(\theta _0 + \eta_i (\frac{|x|}{R} ) (\theta (x) - \theta_0 ) \right) \, dx.
\end{array}
\eeq
The functions $ v_1 - \tilde{u}_1 - u_2 $ and 
$ \theta^* =  \theta - \sum_{i=1}^2 
\left(\theta _0 + \eta_i (\frac{|x|}{R} ) (\theta (x) - \theta_0 ) \right)  $
belong to $C^1 ( \Om_ {R, AR})$ and 
 $ v_1 - \tilde{u}_1 - u_2 = 1- e^{ i \theta_0 }  = const.$, 
$ \theta ^* = - \theta _0 = const.$ on 
$\Om_ {R, AR} \setminus \Om_ {A_1R, A_4R}$.
Therefore
\beq
\label{3.78}
\int _{\Om _{A_1 R, \, A_4R} }  
\frac{\p }{\p x _1} \left( \I (v_1 -  \tilde{u_1} - u_2 ) \right) \, dx  =0 
\qquad \mbox{ and } \qquad
\int _{\Om _{A_1 R, \, A_4R} }   \frac{ \p \theta ^*}{ \p x_1} \, dx = 0.
\eeq
Using (\ref{3.66}), (\ref{3.67}) and the Cauchy-Schwarz inequality we have 
\beq
\label{3.79}
\Big\vert 
\ds \int _{\Om _{A_1 R, \, A_4R} }   
(\rho^2 - 1) \frac{ \p \theta }{ \p x_1} \, dx  
\Big\vert
\leq 2 \sqrt{2}  \e. 
\eeq
Similarly, from (\ref{3.71}), (\ref{3.73}), (\ref{3.74}) and the 
Cauchy-Schwarz inequality we  get
\beq
\label{3.80}
\Big\vert  
 \ds \int _{\Om _{A_1 R, \, A_4R} }   
\left( \! \left( \! 1 + \eta_i (\frac{|x|}{R} ) (\rho  - 1 )\right) ^2 \! 
 - 1 \! \right)
\frac{\p }{\p x_1} \!
\left(\! \theta _0 + \eta_i (\frac{|x|}{R} ) (\theta  - \theta_0 ) \! \right)  dx
\Big\vert  
\leq C\e. 
\eeq
From (\ref{3.77})$-$(\ref{3.80}) we obtain 
$
| Q(v_1) - Q(u_1) - Q(u_2) |  
\leq C \e 
$
and  (\ref{3.4}) gives
$
|Q(u) - Q(v_1)| 
 \leq C E_{GL}^{\Om _{R, AR} } (u ) \leq C\e.
$
These estimates clearly imply (v). 

\medskip

It remains to prove (vi). 
Assume that assumptions (A1) and (A2) in the introduction are satisfied and let 
$
W (s) = V(s)- V( \ph ^2 (\sqrt{s}))  , 
$
so that $ W(s) = 0 $ for $ s \in [0, 4 ]$. 
It is not hard to see that 
there exists $ C > 0$ 
such that
\beq
\label{i3}
|W(b^2) - W(a^2) | \leq C |b-a| 
\left( a^{2p _0 +1} \1_{\{ a > 2  \} } + b^{2p _0 +1} \1_{\{ b > 2  \} } 
\right)
\; \;  \mbox{ for any } a, b \geq 0.
\eeq
Using (\ref{i1}) and (\ref{i3}), then H\"older's inequality we obtain
\beq
\label{3.82}
\begin{array}{l}
\ds \int _{\R ^N  } 
\Big\vert V( | 1+ u|^2) - V( | 1+ v_1|^2) \Big\vert \, dx
\\
\\
\leq \! \! 
\ds \int _{\Om _{ R, \, A R} } \! 
\Big\vert V( \ph ^2(| 1\! + \! u|) ) \! - \! V( \ph ^2(| 1\! + \! v_1| )) \Big\vert 
+ \Big\vert W( | 1 \! + \!  u|^2) \! - \! W( | 1 \!  + \! v_1|^2) \Big\vert  dx 
\\
\\
\leq C \ds \int _{\Om _{ R, \, A R} }  
\left( \ph ^2(| 1+ u|) - 1 \right)^2 + \left( \ph ^2(| 1+ v_1|) - 1 \right)^2 \, dx 
 \\
 \\
 \qquad + C \ds \int _{\Om _{ R, \, A R} } 
 \Big\vert \, | 1+ u| - | 1+ v_1 | \, \Big\vert
 \left( | 1+ u| ^{2 p_0 + 1} \1 _{\{ | 1+ u| > 2  \} } \right.
 \\ 
 \qquad \qquad \qquad \qquad \left.
 + | 1+ v_1| ^{2 p_0 + 1} \1 _{\{ | 1+ v_1| > 2  \} } \right)\, dx
\\
\\
\leq \! C' \e + C'  \! \! \! \ds \int _{\Om _{ R, \, A R} }   \! \! \!  \! \! \! \! \!
|u \! - \! v_1| \left( | 1\!  + \! u| ^{2^* - 1} \1 _{\{ | 1+ u| > 2  \} }
 + | 1 \! + \! v_1| ^{2^* - 1} \1 _{\{ | 1+ v_1| > 2  \} } \right) \! dx
\\
\\
\leq  C' \e + C' \| u - v_1 \|_{L^{2^*}( \Om _{ R, \, A R})}
\left( \| \, | 1+ u| \1 _{\{ | 1+ u| > 2 \} } \|_{L^{2^*}( \Om _{ R, \, A R})}^{2 ^* -1} \right.
\\
\qquad \qquad \qquad \left. 
+ \| \, | 1+ v_1| \1 _{\{ | 1+ v_1| > 2  \} } \|_{L^{2^*}( \Om _{ R, \, A R})}^{2 ^* -1} \right).
\end{array}
\eeq
From the Sobolev embedding we have 
\beq
\label{3.83}
\begin{array}{l}
\| u - v_1 \|_{L^{2^*} (\R^N)} \leq C_S \|\nabla( u - v_1) \|_{L^{2} (\R^N)} 
\\
\leq C_S( \|\nabla u \|_{L^{2}( \Om _{ R, \, A R})} 
+ \|\nabla v_1 \|_{L^{2}( \Om _{ R, \, A R})}  ) \leq 2 C_S \sqrt{\e }.
\end{array}
\eeq
It is clear that $ | 1+ u| > 2 $ implies $|u| > 1$ and 
$ | 1+ u| < 2|u|$, hence
\beq
\label{3.84}
\begin{array}{l}
\| \, |1+ u| \1 _{\{ | 1+ u| > 2  \} } \|_{L^{2^*}( \Om _{ R, \, A R})}
\\
\leq 2 \|u\| _{L^{2 ^*} (\R^N)} 
\leq 2 C_S \|\nabla u\| _{L^{2 } (\R^N)}  
\leq 2 C_S \left( E_{GL}(u) \right)^{\frac 12}.
\end{array}
\eeq
Obviously, a similar estimate  holds for $ v_1$. 
Combining (\ref{3.82}), (\ref{3.83}) and (\ref{3.84}) we find
\beq
\label{3.85}
\ds \int _{\Om _{ R, \, A R} } 
\Big\vert V( |1+ u|^2) - V( | 1+ v_1|^2) \Big\vert \, dx
\leq C' \e + C'' \sqrt{ \e} \left( E_{GL}(u) \right)^{\frac {2^* -1}{2}}.
\eeq
From (\ref{3.69}) and (\ref{3.70}) it follows that 
$V( | 1+ v_1|^2 ) - V( | 1+ u_1|^2 ) - V( | 1+ u_2|^2 )= 0 $ on 
$ \R^N  \setminus \Om _{A_1R,\,  A_4R}$ and 
$|1+ v_1|, \, |1+ u_1|, \, |1+ u_2| 
\in \left[ \frac{1}{2}, \frac{3}{2} \right]$ on $\Om _{A_1R,\,  A_4R}$.
Then using (\ref{i1}), (\ref{3.66}), (\ref{3.74}) and (\ref{3.71}) we get 
\beq
\label{3.86}
\ds \int _{\Om _{ A_1R, \, A_4 R} } 
|V( |1+ v_1| ^2) | \, dx \leq C \int _{\Om _{ A_1R, \, A_4 R} } (\rho ^2 - 1)^2\, dx
\leq C \e, \qquad \mbox{ respectively } 
\eeq
\beq
\label{3.87}
\begin{array}{l}
{\ds \int _{\Om _{ A_1R, \, A_4 R} } } \! \! \! \! \! \! \! \! 
|V( |1+ u_i| ^2) |\, dx 
\leq C \! {\ds \int _{\Om _{ A_1R, \, A_4 R} } } \! 
\left( \! \left(1 + \eta_i (\frac{|x|}{R} ) (\rho  - 1 )\right) ^2 \! - 1 \! \right) ^2 dx
\leq C\e.
\end{array}
\eeq
Therefore
\beq
\label{3.88}
\begin{array}{l}
\ds \int _{\R^N  } 
\Big\vert  V( | 1+ v_1|^2 ) - V( | 1+ u_1|^2 ) - V( | 1+ u_2|^2 ) \Big\vert \, dx 
\\
\leq \ds \int _{\Om _{ A_1R, \, A_4 R} } 
|  V( | 1+ v_1|^2 ) |
+  |  V( | 1+ u_1|^2 ) |
+  |  V( | 1+ u_2|^2 ) |  \, dx
\leq C\e. 
\end{array}
\eeq
Then (iv) follows from (\ref{3.85}) and (\ref{3.88}) and Lemma \ref{splitting} is proven.
\hfill
$\Box $

\section{The variational framework }

The aim of this section is to study the properties of the  functionals $E_c$, $A$, $B_c$ and $P_c$
introduced  in  (\ref{Ec}), (\ref{A}), (\ref{Bc}) and (\ref{Pc}), respectively. 
We assume throughout that the assumptions (A1) and (A2)  in the introduction are satisfied.
Let
$$
\Co = \{ u \in \Xo \; | \; u \neq 0, P_c (u) = 0 \}.
$$

In particular, we will prove that $ \Co \neq \emptyset $ and 
$\ds \inf \{ E_c (u) \; | \; u \in \Co \} > 0$. 
This will be done in a sequence of lemmas.
In the next sections we show that $E_c$ admits a minimizer in 
$ \Co $ and this minimizer 
is a solution of (\ref{1.3}).

We begin by proving that the above mentioned functionals are well-defined on $ \Xo$. 
Since we have already seen in section 2 that $Q$ is well-defined on $\Xo$, 
all we have to do is to prove   that $V(|1+ u|^2) \in L^1(\R^N)$
for any $ u \in \Xo$. This will be done in the next lemma. 

\begin{Lemma}
\label{L4.1}
For any $ u \in \Xo $ we have $V(|1+ u|^2) \in L^1(\R^N)$.
Moreover, for any $ \de > 0$ there exist $C_1 (\de), \, C_2 ( \de ) >0$ such that 
for any $ u \in \Xo $ we have 
\beq
\label{4.1}
\begin{array}{l}
\ds \frac{1 - \de }{2}   \int_{\R^N}  
\left( \ph^2(|1+ u|) - 1 \right)^2 \, dx 
- C_1 ( \de ) \|\nabla u \|_{L^2(\R^N)}^{2^*} 
\\
\\
\leq 
\ds \int_{\R^N}  V(|1+ u |^2 ) \, dx 
\\
\\
\leq 
\ds \frac{1 + \de }{2}   \int_{\R^N}  
\left( \ph^2(| 1+ u|) - 1 \right)^2 \, dx 
+ C_2 ( \de ) \|\nabla u \|_{L^2(\R^N)}^{2^*}  .
\end{array}
\eeq
\end{Lemma}

{\it Proof. } 
Fix $ \de > 0$. 
Using (\ref{1.4}) we see that there exists $ \beta = \beta(\de ) \in ( 0, 1  ]$
 such that 
\beq
\label{4.2}
\frac{1 - \de}{2} (s - 1)^2 \leq V(s) \leq \frac{1 + \de}{2}  (s - 1)^2
\quad \mbox{ for any } s \in ((1  - \beta)^2 , (1 + \beta)^2). 
\eeq
Let $ u \in \Xo$. If $|u(x)| < \beta$ we have 
$ |1+ u(x)|^2 \in ((1  - \beta)^2 , (1 + \beta)^2) $ 
and it follows from (\ref{4.2}) that 
$V(|1 + u|^2)  \1_{\{ |u | < \beta \} } \in L^1( \R^N)$ and
\beq
\label{4.3}
\begin{array}{l}
\ds \frac{1- \de}{2}  \int_{ \{ |u| < \beta \} }
\left( \ph^2(|1+ u|) - 1 \right)^2 \, dx 
\leq \ds \int_{ \{ |u| < \beta \} } V(|1+ u|^2) \, dx 
\\
\\
\leq \ds \frac{1+ \de}{2}  \int_{ \{ |u| < \beta \} }
\left( \ph^2(|1+ u|) - 1 \right)^2 \, dx. 
\end{array}
\eeq
Assumption (A2) implies that there exists $ C_1 '( \de ) > 0 $ such that 
$$
\Big\vert V(| 1+ z |^2) - \frac{1 - \de}{2} (\ph ^2(|1+ z|) -1)^2 \Big\vert
\leq C_1'(\de) |z|^{2 p_0 + 2} 
\leq  C_1''(\de) |z|^{2 ^*}
$$
 for any $ z \in \C$ satisfying $|z| \geq \beta$. 
Using the Sobolev embedding we obtain
\beq
\label{4.4}
\begin{array}{l}
\ds \int_{ \{ |u| \geq \beta \} }
\Big\vert V(| 1+ u |^2) - \frac{1 - \de}{2}  (\ph ^2(|1+ u|) - 1)^2 \Big\vert
\, dx 
\\
\\
\leq
C_1''(\de) \ds \int_{ \{ |u| \geq \beta \} } |u |^{2^*} \, dx 
\leq C_1''(\de) \ds \int_{\R^N} |u |^{2^*} \, dx 
\leq C_1 ( \de ) \| \nabla u \|_{L^2(\R^N)}^{2^*}. 
\end{array}
\eeq
Consequently $V(|1+u|^2)  \1_{\{ |u | \geq \beta \} } \in L^1( \R^N)$ and
it follows from (\ref{4.3}) and (\ref{4.4})
that the first inequality in (\ref{4.1}) holds; 
the proof of the second inequality is similar. 
\hfill
$\Box $

\begin{Lemma}
\label{L4.2}
Let $ \de \in (0, 1)$ and let $ u \in \Xo $  such that 
$ 1 - \de \leq |1+ u| \leq 1 + \de $ a.e. on $\R^N$. Then 
$$
| Q(u) | \leq \frac{ 1}{\sqrt{2} (1 - \de )} E_{GL}(u).
$$
\end{Lemma}

{\it Proof. } 
From Lemma \ref{lifting} we know that there are two real-valued functions 
$ \rho, \, \theta$ such that $ \rho - 1 \in H^1( \R^N)$, $ \theta \in \DR $ and 
$ 1+ u = \rho e^{ i \theta } $ a.e. on $ \R^N$.
Moreover, from (\ref{2.3}) and Definition \ref{D2.4} we infer that 
$$
Q(u) = - \ds \int_{\R^N} ( \rho ^2 - 1) \theta _{x_1} \, dx.
$$
Using the Cauchy-Schwarz inequality we obtain
$$
\begin{array}{l}
\sqrt{2} ( 1 - \de ) | Q(u) | \leq \sqrt{2}  ( 1 - \de ) \| \theta _{x_1} \|_{L^2 (\R^N)}
\|  \rho ^2 - 1  \|_{L^2 (\R^N)} 
\\
\\
\leq 
(1 - \de )^2 \ds \int_{\R^N} | \theta _{x_1} | ^2 \, dx 
+ \frac 12 \int_{\R^N}  \left( \rho ^2 - 1 \right)^2 \, dx 
\\
\\
\leq
\ds \int_{\R^N} 
  \rho ^2 |\nabla \theta |^2 
 + \frac 12 \left( \rho ^2 - 1 \right)^2 \, dx  \leq E_{GL}(u).  
 \end{array}
$$
$\; $ 

\vspace{-26pt}
$\; $
\hfill 
$\Box $

\begin{Lemma}
\label{L4.3}
Assume that $ 0 \leq c < v_s $ and let $ \e \in (0, 1 - \frac{c}{v_s})$. 
There exists a constant $ K _1 = K_1( F, N, c, \e ) > 0$ 
such that for any $ u \in \Xo $ satisfying $E_{GL}(u) < K_1$ we have
$$
\ds \int_{\R^N} |\nabla u |^2 \, dx +  \int_{\R^N}  V(| 1+ u |^2 ) \, dx 
- c | Q(u) | \geq \e E_{GL}(u). 
$$
\end{Lemma}

{\it Proof. } 
Fix $ \e _1   $ such that $ \e < \e _1  < 1 - \frac{c}{v_s}$.
Then fix $ \de _1 \in (0, \e _1 - \e )$. 
By Lemma \ref{4.1}, there exists $C_1 ( \de _1 ) > 0$ such that 
for any $ u \in \Xo$ there holds
\beq
\label{4.5}
\ds \int_{\R^N} V(|1+ u|^2) \, dx 
\geq
\frac{1 - \de _1 }{2} \ds \int_{\R^N}  
\left( \ph^2(|1+ u|) - 1 \right)^2 \, dx 
- C_1 ( \de _1 ) \left(E_{GL}(u)\right) ^{\frac{2^*}{2}}.
\eeq

Using (\ref{3.4}) we see that there exists $A>0$ such that for any 
$ w \in \Xo $ with $E_{GL}(w) \leq 1$, for any 
$h \in (0, 1]$ and for any   minimizer $ v_h $ 
of $G_{h, \R^N}^w $ in $ H_w^1(\R^N)$
we have 
\beq
\label{4.6}
|Q(w) - Q( v_h) | \leq A h^{\frac 2N} E_{GL}(w).
\eeq
Choose $ h \in (0, 1] $ such that 
$
\e _1 - \de _1 - cA h^{\frac 2N} > \e 
$
(this choice is possible because $\e_1 - \de _1- \e > 0$). 
Then fix $ \de > 0 $ such that 
$
\frac{c}{ \sqrt{2} ( 1 - \de )} < 1 - \e _1 
$
(such  $\de $  exist  because 
$\e _1 < 1 - \frac{c}{v_s} = 1- \frac{c}{\sqrt{2} }  $).

Let $K = K( N, h , \de , 1)$ be as in Lemma \ref{L3.1} (iv). 

Consider $ u \in \Xo $ such that $E_{GL}(u) \leq \min (K, 1)$. 
Let $ v_h $ be a minimizer of $G_{h, \R^N} ^u$ in $H_u^1( \R^N)$. 
The existence of $ v_h $ follows from Lemma \ref{L3.1} (i).
By Lemma \ref{L3.1} (iv) we have 
$
1 - \de < |1+ v_h | < 1 + \de 
$
a.e. on $ \R^N $
and then Lemma \ref{L4.2} implies
\beq
\label{4.7}
c | Q( v_h)| \leq \frac{c}{ \sqrt{2} ( 1 - \de )} E_{GL}( v_h )
\leq ( 1 - \e _1 ) E_{GL}( v_h ) \leq ( 1 - \e _1 ) E_{GL}( u).
\eeq

\vspace{-6pt}

We have:
\beq
\label{4.8}
\begin{array}{l}
\ds \int_{\R^N} |\nabla u |^2 \, dx +  \int_{\R^N}  V(| 1+ u |^2 ) \, dx 
- c | Q(u) | 
\\
\\
\geq (1 - \de _1)  E_{GL}(u) - 
 C_1 ( \de _1 ) \left(E_{GL}(u)\right) ^{\frac{2^*}{2}}
- c|Q(u)| \qquad \mbox{ by (\ref{4.5}) }
\\
\\
\geq (1 - \de _1)  E_{GL}(u) - 
 C_1 ( \de _1 ) \left(E_{GL}(u)\right) ^{\frac{2^*}{2}} 
-c| Q(u) - Q( v_h)| - c | Q( v_h)|
\\
\\
\geq 
(1 - \de _1)  E_{GL}(u) 
-  C_1 ( \de _1 ) \left(E_{GL}(u)\right) ^{\frac{2^*}{2}}  
- c A h^{\frac 2N} E_{GL}(u ) - ( 1 - \e _1 ) E_{GL}( u)
\\
 \qquad \qquad \qquad \qquad \qquad \qquad \qquad \qquad \qquad \qquad \qquad \qquad
  \mbox{ by (\ref{4.6}) and  (\ref{4.7})}
\\
= \left( \e _1 - \de _1 - c A h^{\frac 2N} - C_1 ( \de _1) 
\left(E_{GL}(u)\right) ^{\frac{2^*}{2}-1} \right) E_{GL}(u).
\end{array}
\eeq
\vspace{-5pt}
\noindent
Note that (\ref{4.8}) holds for any $ u \in \Xo $ with 
$E_{GL}(u) \leq \min (K, 1)$. 
Since $ \e _1 - \de _1 - c A h^{\frac 2N} > \e$, 
it is obvious that 
$
 \e _1 - \de _1 - c A h^{\frac 2N} - C_1 ( \de _1) 
\left(E_{GL}(u)\right) ^{\frac{2^*}{2}-1} > \e 
$
if $E_{GL}(u) $ is sufficiently small and the conclusion of Lemma 
\ref{L4.3} follows.
\hfill
$\Box $

\medskip

An obvious consequence of  Lemma \ref{L4.3} is that $E_c(u) >0$ if 
$ u \in \Xo \setminus \{ 0 \}$ and $ E_{GL}(u)$ is sufficiently small. 
The next lemma implies that there are functions $ v \in \Xo $ such that  $E_c (v ) < 0$. 

\begin{Lemma}
\label{L4.4}
Let $ N \geq 2$.
There exists a continuous map from $[2, \infty) $ to $H^1( \R^N)$, 
$R \longmapsto v_R $ such that 
$v_R \in C_c ( \R^N)$ for any $ R \geq 2$ 
and the following estimates hold:

\medskip

i) $\ds \int_{\R^N} |\nabla v_R|^2 \, dx 
\leq C_1 R^{N-2} + C_2 R^{N-2} \ln R,  $

\medskip

ii) $ \Big\vert \ds \int_{\R^N}  V(|1+ v_R|^2 ) \, dx \Big\vert
\leq C_3  R^{N-2} ,  $

\medskip

iii) $ \Big\vert \ds \int_{\R^N}  
\left( \ph ^2 (|1+ v_R| ) - 1 \right)^2 \, dx \Big\vert
\leq C_4  R^{N-2} ,  $

\medskip

iv) $ - 2 \pi  \omega_{N-1} R^{N-1} \leq Q( v_R) 
\leq - 2 \pi  \omega_{N-1} (R- 2 ) ^{N-1}, $

\medskip

\noindent 
where the constants $C_1-C_4$ depend only on $N$ and 
$\omega_{N-1} = \Lo ^{N-1}(B_{\R^{N-1}}(0,1)) $.
\end{Lemma}

{\it Proof. } 
Let
\vspace{-7pt}
$$\begin{array}{l}
T_{ R }= \{ x = (x_1, x')  \in \R^N \; \big| \; 0 \leq |x'| \leq R \mbox{ and } 
 -R +|x'|  < x_1 <R- |x'| \}.
 \end{array}
 $$

 We define $ \theta _R : \R^N \lra \R$ in the following way: 
if $|x'| \geq R$ we put $\theta _R (x) = 0$ and if $|x'| < R$ we define
\vspace{-10pt}
\beq
\label{4.9}
\theta _R (x) =
\left\{
\begin{array}{l}
0 \qquad \mbox{ if } x_1 \leq  - R + |x'|  , 
\\
\\
\frac{ \pi }{R-|x'| } x_1 + \pi \qquad \mbox{ if } x \in T_{ R  }, 
\\
\\
2 \pi \qquad \mbox{ if } x_1 \geq   R- |x'| .
\end{array}
\right.
\eeq
It is easy to see that $ x \longmapsto e^{i \theta _R (x)}$ 
is continuous on 
$ \R^N \setminus \{ x \; | \; x_1 = 0 , \; |x'| = R \}$ and equals 1 on 
$\R^N \setminus  T_{ R  }$.

Fix $ \psi \in C^{\infty }(\R)$  such that 
$ \psi = 0 $ on $(-\infty, 1]$, $\psi = 1$ on $[2, \infty )$ and 
$ 0 \leq \psi ' \leq 2$. 
Let 
\vspace{-5pt}
\beq
\label{4.10}
\psi _R (x) = \psi \left( \sqrt{x_1 ^2 + ( | x'|  - R)^2} \right)
\quad 
\mbox{ and }
\quad
v_R (x) = \psi_R(x) e^{i \theta _R(x)} -1.
\eeq
It is obvious that $ v_R \in C_c(\R^N)$
(in fact, 
$ v_R $ is $ C ^{\infty }$ on 
$  \R^N \setminus B$, where 
$ B = \p T_{ R} \cup \{ (x_1, 0, \dots, 0) \; | \; x_1 \in [-R,R ] \} $).
On $  \R^N \setminus B$ we have
\beq
\label{4.11}
\frac{\p \theta _R }{\p x_1 } = \left\{
\begin{array}{l}
\frac{ \pi }{ R - | x'|} \quad \mbox{ if } x \in T_{R} ,
\\
0 \quad \mbox{ otherwise }, 
\end{array}
\right.
\quad
\frac{\p \theta _R }{\p x_j } = \left\{
\begin{array}{l}
\frac{ \pi   x_1}{( R - | x'|)^2} \frac{x_j}{|x'|}
\quad \mbox{ if } x \in T_{A,R} , 
\\
0 \quad \mbox{ otherwise,}
\end{array}
\right.
\eeq
\vspace{-7pt}
\beq
\label{4.12}
\frac{ \p \psi _R}{ \p x_1 } (x) 
=  \psi ' \left(  \sqrt{ x_1 ^2 + (|x'| - R )^2 } \right)
\frac{x_1}{\sqrt{ x_1 ^2 + (|x'| - R )^2 }}, 
\eeq
\vspace{-7pt}
\beq
\label{4.13}
\frac{ \p \psi _R }{ \p x_j } (x) 
=  \psi ' \left(  \sqrt{ x_1 ^2 + (|x'| - R )^2 } \right)
\frac{|x'| - R}{\sqrt{ x_1 ^2 + (|x'| - R )^2 }} \frac{x_j}{|x'|}
\quad \mbox{ for } j \geq 2 \mbox{ and } x' \neq 0.
\eeq
Then a simple computation gives 
$ \langle i \frac{\p  v_R }{\p x_1}, v_R \rangle 
 = - \psi _R ^2 \frac{ \p \theta _R}{\p x _1} + \frac{ \p }{\p x_1} \left( \I (v_R)  \right) $
on $ \R^N \setminus B$. 
Thus we have 
$$ 
Q(  v_R ) = -  \ds \int _{\R^N } 
\psi _R ^2 \frac{\p \theta _R }{\p x_1 }  \, dx .
$$
\vspace{-1pt}
It is obvious that 
\beq
\label{4.14}
{\ds \int_{- \infty }^{\infty } }
\frac{\p \theta _R }{\p x_1 } \, d x_1 = 0 \quad \mbox{ if }   |x'| > R 
\quad
\mbox{ and }
\quad 
{\ds \int_{- \infty }^{\infty } }
\frac{\p \theta _R }{\p x_1 } \, d x_1 =  2 \pi 
\quad \mbox{  if }
0 <  |x'| < R. 
\eeq
Since $ \frac{\p \theta _R }{\p x_1 } \geq 0$ a.e. on $ \R^N$ and 
$ 0 \leq \psi _R \leq 1$, we get 
$$
\ds \int_{\{ | R - |x'|\, | \geq 2  \} }
\frac{\p \theta _R }{\p x_1 } \, dx 
\leq
\int_{\R^N}  \psi _R^2 \frac{\p \theta _R }{\p x_1 } \, dx 
\leq 
\int_{\R^N}  \frac{\p \theta _R }{\p x_1 } \, dx , 
$$
and using Fubini's theorem and (\ref{4.14}) we obtain that $v_R$
satisfies (iv). 

\medskip

Using cylindrical coordinates $(x_1, r, \zeta )$ in $\R^N$, where 
$ r = |x'|$ and $ \zeta  =  \frac{ x'}{|x'|}  \in  S^{N-2}$, we get 
\beq
\label{4.15}
\begin{array}{l}
{\ds \int_{\R^N} } \! \! V(|1 + v_R |^2 ) \, dx
= |S^{N-2}| 
{\ds \int_{- \infty }^{\infty } 
\int_{0 }^{\infty } }
\! \! V\! \! \left( \!  \psi ^2 \left(\sqrt{ x_1 ^2 + (r - R )^2 }  \! \right) \! \right)
r^{ N-2} \, dr \, d x_1. 
\end{array}
\eeq
Next we use polar coordinates in the $(x_1, r)$ plane, that is we write 
$ x_1 = \tau \cos \al $, $ r = R +  \tau \sin  \al $ 
(thus $ \tau = \sqrt{ x_ 1 ^2 + ( R - r ) ^2}$).
Since $ V(  \psi ^2 (s) ) = 0 $ for $ s \geq 2$, we get 
\beq
\label{4.16}
{ \ds \int_{\! - \infty }^{\infty } \!
\int_{0 }^{\infty } } \! \! \! 
V  \!    \left( \!   \psi ^2 \! \left( \! \! \sqrt{ x_1 ^2 + \! (r \! - \! R )^2 } \right) \! \right) \!
r^{ N-2}  dr \, d x_1
= \! { \ds  \int_{0}^{ 2  } \! \!  \int_{0 }^{ 2 \pi } } \! \!
V(  \psi ^2 ( \tau )) (R + \! \tau \sin \al )^{N-2} \tau \, d \al \,  d \tau .
\eeq
It is obvious that 
$\Big\vert  \ds \int_{0 }^{ 2 \pi }  ( R + \tau \sin \al )^{N-2} d \al \Big\vert
\leq 2 \pi (R + 2  )^{N-2}$ for any $ \tau \in [0,2]$, 
and then using (\ref{4.15}) and (\ref{4.16}) we infer that 
$v_R$ satisfies (ii). 
The proof of (iii) is similar. 

\medskip

It is clear that on $ \R^N \setminus B$ we have  
\beq
\label{4.17}
| \nabla v_R | ^2 =  | \nabla \psi _R |^2 + 
 |\psi _R |^2 |\nabla \theta _R|^2 .
\eeq
From (\ref{4.12}) and  (\ref{4.13}) we see that 
$
| \nabla \psi _R (x)  |^2 
=  \Big\vert \psi ' \left (  \sqrt{ x_1 ^2 + (|x'| - R )^2 } \right) \Big\vert ^2 .
$
Proceeding as above and using cylindrical coordinates $(x_1, r , \zeta )$ in 
$\R^N$, then passing to polar coordinates 
$ x_1 = \tau \cos \al$, $ r =  R + \tau \sin \al $, we obtain 
\beq
\label{4.18}
\ds \int _{\R^N} 
\Big\vert 
\psi ' \left ( \sqrt{ x_1 ^2 + (|x'| - R )^2 } \right)
\Big\vert ^2  \, dx
\leq
2 \pi |S^{N-2} |  (R+ 2  )^{N-2} \int_0 ^2 s |\psi '(s) |^2  \, ds.
\eeq

It is easily seen from (\ref{4.11}) that 
$|\nabla \theta _R (x) | ^2 =  \frac{\pi ^2  }{ (R - |x'|)^2} 
\left( 1 + \frac{ x_1 ^2}{(R- |x'|)^2} \right) $
if $ x \in T_{ R}$, $|x'| \neq 0$, 
and $\nabla \theta _R (x)  = 0 $ a.e. on  $  \R^N \setminus \ov{T}_{ R}$. 
Moreover, if $(x_1, x') \in T_{ R}$ and 
$|x'| \geq R - \frac{ 1}{\sqrt{ 2}}$, we have 
$\psi _R ( x_1, x') = 0$. Therefore
\beq
\label{4.19} 
\begin{array}{l}
\ds \int_{\R^N} 
|\psi _R |^2 |\nabla \theta _R |^2 \, dx 
\leq
\ds \int_{ T_{ R}\cap \{ |x'| < R - \frac{1 }{\sqrt{2}} \} } 
|\nabla \theta _R |^2 \, dx 
\\
\\
= \ds \int_{ \{ |x'| < R - \frac{ 1 }{\sqrt{ 2}} \} }
\int _{- R + |x'|  } ^{ R - |x'|}
|\nabla \theta _R |^2 \, dx_1 \, d x'
\\
\\
= \ds \int_{ \{ |x'| < R - \frac{1 }{\sqrt{2}} \} } 
\frac{ 2 \pi ^2 }{ R - |x'|}
+ \frac {2 \pi ^2 }{3}  \frac{1}{R - |x'|} \, d x'
= \ds \frac 83 \pi ^2  |S^{N-2}|
\ds \int_0^{ R - \frac{ 1 }{\sqrt{2}} } 
\frac{r^{N-2}}{R - r} \, dr
\\
\\
\ds = \frac 83 \pi ^2  |S^{N-2}| R^{N-2}
\left( - \sum_{k =1}^{N-2} 
\frac 1k \left( \! 1 - \frac{1 }{R \sqrt{ 2}}\right)^k \! \!
+ \ln \left( R \sqrt{2} \right) \! \right)\! .
\end{array}
\eeq

From (\ref{4.17}), (\ref{4.18}) and (\ref{4.19}) it follows that  $ v _R$
satisfies (i). 
It is not hard to see that the mapping $R \longmapsto v_R $ 
is continuous from $ [2, \infty)$ to $ H^1( \R^N)$ and Lemma \ref{L4.4} is proven.
\hfill
$\Box $

\begin{Lemma}
\label{L4.5}
For any $ k > 0$, the functional $Q$ is bounded on the set 
$$
\{ u \in \Xo \; | \; E_{GL}(u) \leq k \}.
$$ 
\end{Lemma}

{\it Proof. } 
Let $ c \in (0, v_s)$ and let $ \e \in (0, 1 - \frac{c}{v_s} )$. 
From Lemmas \ref{L4.1} and \ref{L4.3} it follows that there exist
two positive constants $ C_2( \frac{\e }{2})$ and  $K_1 $ such that for 
any $ u \in \Xo $ satisfying $ E_{GL}(u) < K_1 $ we have
$$
\begin{array}{l}
(1 + \frac{ \e}{2} ) E_{GL}(u) + C_2( \frac{\e }{2})  \left( E_{GL}(u) \right)^{\frac{2^*}{2}}
- c | Q(u)|
\\
\\
\geq 
\ds \int_{\R^N} |\nabla u |^2 \, dx +  \int_{\R^N}  V(| 1+ u |^2 ) \, dx 
- c | Q(u) | \geq \e E_{GL}(u). 
\end{array}
$$
This inequality implies that there exists $ K_2 \leq K_1 $ such that 
for any $ u \in \Xo $ satisfying $ E_{GL}(u) \leq K_2 $ we have
\beq
\label{4.20}
c | Q(u)| \leq E_{GL}(u). 
\eeq
Hence Lemma \ref{L4.5} is proven if $ k \leq K_2$. 

Now let $ u \in \Xo $ be such that $E_{GL}(u) > K_2$. 
Using the notation (\ref{1.7}), it is clear that  for $ \si > 0 $ we have
$Q( u_{\si, \si }) = \si ^{N-1} Q(u) $ (see (\ref{2.13})) and 
$$
E_{GL}( u_{\si, \si }) = \si ^{N-2} \ds \int_{\R^N} |\nabla u |^2 \, dx 
+ \frac{\si ^N }{2} \ds \int_{\R^N} \left( \ph ^2(| 1+ u|) - 1 \right)^2 \, dx .
$$
Let $ \si _0 = \left( \frac{K_2}{E_{GL}(u)  }\right)^{\frac{1}{N-2}}$.
Then $ \si _0 \in (0, 1)$ and we have 
$E_{GL}(u_{\si _0, \si _0} ) \leq \si_0 ^{N-2} E_{GL}(u) = K_2$. 
Using (\ref{4.20}) we infer that 
$c|Q( u_{\si _0, \si _0})| \leq E_{GL}(u_{\si _0, \si _0} )$, 
and this implies 
$ c \si _0 ^{N-1} |Q(u)| \leq \si_0 ^{N-2} E_{GL}(u) $, or equivalently
\beq
\label{4.21}
| Q(u)| \leq \frac{1}{ c \si _0 } E_{GL}(u) 
= \frac 1c K_2^{- \frac{1}{N-2}} \left( E_{GL}(u)  \right)^{\frac{N-1}{N-2}}.
\eeq
Since (\ref{4.21}) holds for any $ u \in \Xo$ with $E_{GL}(u) > K_2$, 
Lemma \ref{L4.5} is proven.
\hfill
$\Box $

\medskip

From Lemma \ref{L4.1} and Lemma \ref{L4.5} it follows that for any $ k >0$, 
the functional $E_c$ 
is bounded on the set $ \{ u \in \Xo \; | \; E_{GL}(u) = k \}$.  
For $ k > 0$ we define 
$$
E_{c,\, min}(k) = \inf \{ E_c (u) \; | \; u \in \Xo, \; E_{GL}(u) = k \}.
$$
Clearly, the function $ E_{c,\, min} $ is bounded on any bounded interval of $\R _+$.
The next result gives some basic properties of $E_{c, min} $ which will be important for our variational argument. 

\begin{Lemma}
\label{L4.6}
Assume that $ N \geq 3 $ and $ 0 < c < v_s$. 
The function $E_{c,\,  min}$ has the following properties: 

i) There exists $k_0 > 0$ such that $E_{c,\,  min}(k) > 0$
for any $ k \in (0, k_0)$.

ii) We have $ \ds \lim _{k \ra \infty} E_{c,\,  min}(k)  = -\infty$. 

iii) For any $ k > 0$ we have $E_{c,\,  min}(k) < k$. 
\end{Lemma}

{\it Proof. } 
(i) is an easy consequence of Lemma \ref{L4.3}.

\smallskip

(ii) It is obvious that $H^1(\R^N) \subset \Xo $ and the functionals $ E_{GL}$, 
$E_c$ and $Q$ are continuous on $H^1(\R^N)$. For $ \e = 1$ and $ R > 2$, 
consider the functions $ v_R $ constructed in Lemma \ref{L4.4}.
Clearly, $R \longmapsto v_R $ is a continuous curve in $ H^1( \R^N)$. 
Lemma \ref{L4.4} implies $E_c ( v_R ) \lra - \infty $ as $ R \lra \infty $. 
From Lemma \ref{L4.5} we infer that 
$E_{GL}(v_R ) \lra \infty $ as $ R \lra \infty $ 
and then it is not hard to see that 
(ii) holds. 

\smallskip

(iii) Fix $ k > 0$. 
Let $v_R$ be as above and let $ u = v_R $ 
for some $ R$ sufficiently large, so that 
$$
E_{GL}(u) > k, \quad Q(u) < 0 \quad \mbox{  and } E_c(u) < 0.
$$
In particular, we have
$$
E_c(u) - E_{GL}(u) = c Q(u) + \ds \int_{\R^N} V(|1+ u |^2) 
- \frac 12  \left( \ph^2(| 1+ u|^2) - 1 \right)^2 \, dx 
< 0.
$$
It is obvious that $ E_{GL}( u_{\si, \si }) \lra 0 $ as $ \si \lra 0$, 
hence there exists $ \si _0 \in (0, 1)$ such that $E_{GL}(u_{\si _0, \si _0} )=k$.
Using the fact that $E_{GL}(u) - E_c (u) < 0$ and $Q(u) < 0$ we get 
$$
\begin{array}{l}
E_c( u_{\si _0, \si _0 }) - E_{GL}(u_{\si _0, \si _0} )
\\
= \si _0 ^{N-1} c Q(u) 
+ \si _0 ^N \ds \int_{\R^N} V(| 1 + u |^2) 
- \frac 12   \left( \ph^2(| 1 +  u|^2) - 1 \right)^2 \, dx 
\\
= ( \si _0 ^{N-1} - \si _0 ^N) c Q(u) + \si _0 ^N (E_c(u) - E_{GL}(u))
< 0 .
\end{array}
$$
Thus $E_c( u_{\si _0, \si _0 }) < E_{GL}(u_{\si _0, \si _0} )$. 
Since $E_{GL}(u_{\si _0, \si _0} ) =k$, we have necessarily
$E_{c,\,  min}(k)  \leq E_c( u_{\si _0, \si _0 }) < k$. 
\hfill
$\Box $

\medskip

 From Lemma \ref{L4.6} (i) and (ii) it follows that 
\beq
\label{4.22}
0 < S_c : = \sup \{ E_{c, min }(k) \; | \; k > 0 \} < \infty .
\eeq

\begin{Lemma}
\label{L4.7}
The set 
$
\Co = \{ u \in \Xo \; | \; u \neq 0,\;  P_c (u) = 0 \}
$
is not empty and we have 
$$
T_c := \inf \{ E_c (u) \; | \; u  \in \Co  \} \geq S_c > 0.
$$ 
\end{Lemma}

{\it Proof. } 
Let $ w \in \Xo \setminus \{0 \}$ be such that $ E_c(w) < 0$
(we have  seen in the proof of Lemma \ref{L4.6} that such functions $w$ exist). 
It is obvious that  $ A(w) > 0 $ and 
$ \ds \int_{\R^N} \Big\vert \frac{ \p w }{ \p x_1 } \Big\vert ^2  \, dx > 0$; 
therefore $ B_c (w) = E_c (w) - A(w) < 0 $ and $ P_c (w) = E_c (w) - \frac{2}{N-1} A(w) < 0.$
Clearly, 
\beq
\label{4.22bis}
P_c( w_{\si, 1}) = 
\frac{1}{\si } \int_{\R^N} \Big\vert \frac{ \p w }{ \p x_1 } \Big\vert ^2  \, dx 
+ \frac{N-3}{N-1} \si  A(w)
+ c Q(w) + \si \int_{\R^N} V(| 1+ w|^2) \, dx.
\eeq
Since $P_c( w_{1, 1}) = P_c( w) < 0$ and 
${\ds \lim _{ \si \ra 0 }} P_c( w_{\si, 1}) = \infty $, 
there exists $ \si _0 \in (0, 1) $ such that $P_c( w_{\si _0, 1})= 0$, that is 
$  w_{\si _0, 1} \in \Co $.  Thus $ \Co \neq \emptyset$.

To prove the second part of Lemma \ref{L4.7}, 
consider first the case  $ N \geq 4$. Let  $ u \in \Co $. 
It is clear that $ A(u) > 0$, $B_c (u) = - \frac{N-3}{N-1} A(u) < 0$ 
and for any $ \si > 0 $ we   have 
$ E_c( u _{1, \si } ) = A(u_{1, \si } ) + B_c (u_{1, \si })  
= \si ^{N-3} A(u) + \si ^{N-1} B_c(u)$ , hence
$$
\frac{d }{d \si } (  E_c( u _{1, \si } ) )
= (N-3 ) \si ^{N-4} A(u) + (N-1) \si ^{N-2} B_c(u)
$$
is positive on $ ( 0, 1)$ and negative on $(1, \infty)$. 
Consequently the function $ \si \longmapsto  E_c( u _{1, \si } ) $ achieves its
maximum at $ \si = 1$. 

On the other hand, we have 
$$
 E_{GL}( u _{1, \si } ) =  \si ^{N-3} A(u) + \si ^{N-1} 
 \left( \int_{\R^N} \Big\vert \frac{ \p u }{ \p x_1 } \Big\vert ^2  
 + \frac 12  \left( \ph ^2(| 1+  u |) - 1 \right)^2 \, dx \right).
$$
It is easy to see that the mapping $ \si \longmapsto  E_{GL}( u _{1, \si } ) $ 
is strictly increasing and one-to-one from $(0, \infty )$ to  $(0, \infty )$.
Hence for any $ k > 0$, there is a unique $ \si ( k, u ) > 0$ such that 
$E_{GL}( u _{1, \si (k, u)} ) = k$. 
Then we have 
$$
E_{c, min }(k) \leq E_c (  u _{1, \si (k, u)} ) \leq E_c(  u _{1, 1}) = E_c(  u). 
$$
Since this is true for any $ k > 0 $ and any $ u \in \Co$, the conclusion follows. 

\medskip

Next we consider the case $ N =3$. Let  $ u \in \Co $. 
We have $ P_c(u) =  B_c(u) = 0$ and $ E_c(u) = A(u) > 0$. 
For  $ \si > 0 $ we get 
$$
\begin{array}{c}
E_c( u_{1, \si }) = A(u ) + \si ^2 B_c(u ) = A(u) \qquad \mbox{ and } 
\\
\\
E_{GL}( u_{1, \si }) = A(u) + \si ^2  
\ds \left( \int_{\R^3} \Big\vert \frac{ \p u }{ \p x_1 } \Big\vert ^2  
 + \frac 12 \left( \ph ^2(| 1+ u |) - 1  \right)^2 \, dx \right).
\end{array}
$$
Clearly, $ \si \longmapsto E_{GL}( u_{1, \si }) $ is increasing on $(0, \infty )$ and 
is one-to-one from $(0, \infty )$ to $(A(u), \infty )$.

Fix $ \e > 0$. Consider  $ k_ {\e } > 0 $  such that $ E_{c, min} ( k_{\e }) > S_c - \e$. 
If $A(u) \geq k_{\e }$, from Lemma \ref{L4.6} (iii) we have $E_{c, min}(k_{\e}) < k_{\e}$, hence
$$
E_c(u) = A(u) \geq k_{\e } > E_{c, min} ( k_{\e }) > S_c - \e.
$$
If $A(u) < k_{\e }$, there exists $ \si (k_{\e }, u ) > 0$ such that 
$E_{GL}(u_{1,  \si (k_{\e }, u ) }) = k_{\e } $. 
Then we get 
$$
E_c (u) = A(u) = E_c (u_{1,  \si (k_{\e }, u ) }) 
\geq E_{c, min}(k_{\e}) > S_c - \e. 
$$
So far we have proved that for any $ u \in \Co $ and any $ \e > 0$ we have 
$E_c(u) > S_c - \e$. 
The conclusion follows letting $ \e \lra 0$, then taking the infimum for $ u \in \Co$.
\hfill
$\Box $

\medskip

We do not know whether $ T_c = S_c$ in Lemma \ref{L4.7}.

\begin{Lemma}
\label{L4.8}
Let $T_c$ be as in Lemma \ref{L4.7}. The following assertions hold. 

\medskip
 
i) For any $ u \in \Xo $ with $P_c(u) < 0$ we have $ A(u) > \frac{N-1}{2} T_c$. 

\medskip
 
ii) Let $(u_n)_{n \geq 1} \subset \Xo $ be a sequence such that 
$\left(E_{GL}(u_n)\right)_{n \geq 1}$ is bounded and $ 
\ds \lim_{n \ra \infty } P_c ( u_n ) = \mu < 0$. 
Then ${\ds \liminf_{n \ra \infty } } A( u_n) > \frac{N-1}{2} T_c$. 
\end{Lemma}

{\it Proof. } i) Since $P_c(u) < 0$, it is clear that $ u \neq 0$ and thus
$ \ds \int_{\R^N} \Big\vert \frac{ \p u}{\p x_1} \Big\vert ^2 \, dx > 0$. 
As in the proof of Lemma \ref{L4.7}, 
we have $ P_c(u_{1, 1}) = P_c(u) < 0$ and  (\ref{4.22bis}) implies that 
$\ds \lim_{\si \ra 0 } P_c( u_{\si, 1 }) = \infty$, hence there exists 
$ \si _0 \in (0, 1)$ such that $P_c( u_{\si _0, 1}) = 0$. 
From Lemma \ref{L4.7} we get $E_c( u_{\si _0, 1})  \geq T_c $ and this implies 
$E_c( u_{\si _0, 1}) - P_c( u_{\si _0, 1})  \geq T_c $, that is 
$\frac{2}{N-1} A ( u_{\si _0, 1})  \geq T_c$. From the last inequality we find
\beq
\label{4.23}
A(u) \geq \frac{N-1}{2} \frac{ 1}{ \si _0 } T_c > \frac{N-1}{2} T_c.
\eeq

ii) For  $n$ sufficiently large (so that $P_c( u_n) < 0$) we have $ u_n \neq 0$ and
$ \ds \int_{\R^N} \Big\vert \frac{ \p u _n }{\p x_1} \Big\vert ^2 \, dx > 0$.
As in the proof of part (i), using (\ref{4.22bis}) we see that 
for each $n$ sufficiently big there exists $ \si _n \in (0, 1)$
such that 
\beq
\label{4.24} 
P_c( (u_n)_{\si _n, 1 }) = 0 
\eeq
and we infer that 
$A(u_n ) \geq \frac{N-1}{2} \frac{ 1}{ \si _n } T_c$. 
We claim that 
\beq
\label{4.25} 
\limsup_{n \ra \infty }\;  \si _n < 1.
\eeq
Notice that if (\ref{4.25}) holds,  we have 
$$
{\ds \liminf_{n \ra \infty } } A(u_n ) 
\geq \frac{N-1}{2} \frac{ 1}{ \limsup_{n \ra \infty } \si _n } T_c
> \frac{N-1}{2}  T_c
$$
and Lemma \ref{L4.8} is proven. 

To prove (\ref{4.25}) we argue by contradition and assume that there is a subsequence
$ (\si_{n_k} )_{k\geq 1}$ such that $ \si_{n_k} \lra 1 $ as $ k \lra \infty $. 
Since $\left(E_{GL}(u_n)\right)_{n \geq 1}$ is bounded, using
Lemmas \ref{L4.1} and \ref{L4.5} we infer that 
$\left( \ds \int_{\R^N} \Big\vert \frac{ \p u_n }{\p x_1} \Big\vert ^2 \, dx \right)_{n \geq 1}$, 
$\left( \ds \int_{\R^N}  V(|1 + u_n |^2) \, dx \right)_{n \geq 1}$,
$ ( A(u_n) )_{n \geq 1}$, and  $ ( Q(u_n) )_{n \geq 1}$  are bounded. 
Consequently there is a subsequence $( n_{k_{\ell}})_{\ell \geq 1}$ 
and there are  $ \al _1, \, \al _2, \, \beta, \, \g \in \R$  such that 
$$
\begin{array}{ll}
\ds \int_{\R^N} \Big\vert \frac{ \p u_{n_{k_{\ell }}} }{\p x_1}  \Big\vert ^2 \, dx 
\lra \al _1, 
& 
\ds \int_{\R^N}  V(| 1+   u_{n_{k_{\ell }}}  |^2) \, dx  \lra \g
\\
\\
A( u_{n_{k_{\ell }}} ) \lra \al _2,
&
Q( u_{n_{k_{\ell }}} ) \lra \beta 
 \quad \mbox{ as }
\ell \lra \infty.
\end{array}
$$
Writing  (\ref{4.24}) and (\ref{4.22bis}) (with  $(u_{n_{k_{\ell }}})_{\si _{n_{k_{\ell }}} , 1}$ 
instead of $ (u_n)_{\si _n , 1}$ and $w_{\si , 1}$, respectively)
then passing to the limit as $ \ell \lra \infty $ and using the fact that $ \si_{n_k} \lra 1 $
we find
$
\al _1 + \frac{N-3}{N-1} \al _2 + c \beta + \g = 0. 
$
On the other hand we have $\ds \lim_{ \ell \ra \infty } P_c( u_{n_{k_{\ell }}}) = \mu <0$ 
and this gives $ \al _1 + \frac{N-3}{N-1} \al _2 + c \beta + \g = \mu < 0$.
This contradiction  proves that (\ref{4.25}) holds 
and the proof of Lemma \ref{L4.8} is complete. 
\hfill
$\Box $

\section{The case $N \geq 4$}

Throughout this section we assume that $ N \geq 4$,  $ 0 < c < v_s$
and the assumptions (A1) and (A2) in the introduction are satisfied.  
Most of the results below do {\it not} hold for $ c > v_s$. 
Some of them may not hold for $ c =0$ and some particular nonlinearities $F$.

\begin{Lemma}
\label{L5.1}
Let $ (u_n)_{n \geq 1} \subset \Xo $ be a sequence such that 
$(E_c(u_n) )_{n \geq 1}$ is bounded and 
$P_c (u_n) \lra 0 $ as $ n \lra \infty $. 

Then $(E_{GL}(u_n))_{n \geq 1} $ is bounded. 
\end{Lemma}

{\it Proof. } 
We have $ \frac{2}{N-1} A(u_n) = E_c (u_n) - P_c( u_n ) $, hence 
$(A(u_n) )_{n \geq 1}$ is bounded. 
It remains to prove that 
$\ds \int_{\R^N} \Big\vert \frac{ \p u_n } {\p x _1} \Big\vert + 
\frac 12 \left( \ph ^2 ( | 1+ u_n | ) - 1 \right)^2 \, dx $ is bounded. 
We argue by contradiction and we assume that there is a subsequence, still denoted
$ (u_n)_{n \geq 1} $, such that 
\beq
\label{5.1}
{\ds \int_{\R^N}} \Big\vert \frac{ \p u_n } {\p x _1} \Big\vert + 
\frac 12 \left( \ph ^2 ( | 1+ u_n | ) - 1 \right)^2 \, dx 
\lra \infty \qquad \mbox{ as } n \lra \infty .
\eeq
Fix $ k_0 > 0$ such that $E_{c, min}( k_0 ) > 0$. 
Arguing as in the proof of Lemma \ref{L4.7}, it is easy 
to see that there exists a sequence $(\si _n )_{n \geq 1}$ such that 
\beq
\label{5.2}
E_{GL}((u_n)_{1, \si _n}) = \si _n ^{N-3} A(u_n) + \si _n ^{N-1} \int_{\R^ N} 
\Big\vert \frac{ \p u_n } {\p x _1} \Big\vert  + 
\frac 12  \left( \ph ^2 ( | 1+ u_n | ) - 1 \right)^2 \, dx  = k_0 .
\eeq
From (\ref{5.1}) and (\ref{5.2}) we have $ \si _n \lra 0 $ as $ n \lra \infty $. 
Since $ B_c ( u_n ) = - \frac{ N-3}{N-1} A(u_n) + P_c (u_n)$, 
it is clear that $(B_c(u_n) )_{n \geq 1}$ is bounded and 
we obtain 
$$
E_{c}((u_n)_{1, \si _n}) = \si _n ^{N-3} A(u_n) +  \si _n^{N-1} B_c( u_n) 
\lra 0 
\qquad \mbox{ as } n \lra \infty. 
$$
But this contradicts the fact that $E_{c, min}( k_0 ) > 0$ and the proof of Lemma \ref{L5.1} 
is complete. 
\hfill
$\Box $

\begin{Lemma}
\label{L5.2} 
Let $(u_n )_{n \geq 1} \subset \Xo $ be a sequence satisfying the 
following properties:

a) There exist $C_1, \, C_2  > 0$ such that   $C_1 \leq E_{GL}(u_n)$ 
and $A(u_n) \leq C_2 $ for any $ n \geq 1$. 

b) $P_c( u_n ) \lra 0 $ as $ n \lra \infty $. 

Then $\ds \liminf_{n \ra \infty } E_c(u_n) \geq T_c$, 
where $T_c$ is as in Lemma \ref{L4.7}. 
\end{Lemma}

Notice that in Lemma \ref{L5.2} the assumption 
$E_{GL}(u_n) \geq C_1 > 0  $ is   necessary. 
To see this, consider a sequence $(u_n )_{n \geq 1} \subset  H^1( \R^N) $ such that 
$ u_n \neq 0 $ and $ u_n \lra 0 $ as $ n \lra \infty $. 
It is clear that $ P_c( u_n ) \lra 0 $ and $E_c( u_n ) \lra 0 $ as $ n \lra \infty $. 

\medskip

{\it Proof. } 
First we prove that 
\beq
\label{5.3}
C_3 := \liminf_{n \ra \infty } A(u_n) > 0.
\eeq
To see this, fix
 $ k_0 > 0$ such that $E_{c, min}( k_0 ) > 0$. 
Exactly as in the proof of Lemma \ref{L4.7}, it is easy 
to see that for each $ n$ there exists a unique $ \si _n > 0$ such that 
(\ref{5.2}) holds. 
Since $ k_0 = E_{GL}((u_n)_{1, \si _n})  
\geq \min ( \si _n ^{N-3}, \si _n ^{N-1}) E_{GL}((u_n)) 
\geq  \min ( \si _n ^{N-3}, \si _n ^{N-1}) C_1$, 
it follows that $ (\si _n )_{n \geq 1}$ is bounded. 
On the other hand, we have 
$E_c( (u_n)_{1, \si _n})  
= \si _n ^{N-3}A(u_n) + \si _n ^{N-1} B_c( u_n) 
 \geq E_{c, min }( k_0 ) > 0 $, that is 
\beq
\label{5.4}
\si _n ^{N-3}A(u_n) 
+ \si _n ^{N-1} \left( P_c (u_n) - \frac{N-3}{N-1} A(u_n) \right)
\geq E_{c, min }( k_0 ) > 0  .
\eeq
If there is a subsequence $ (u_{n _k})_{k \geq 1}$ such that 
$A(u_{n _k} ) \lra 0 $, putting $u_{n _k} $ in (\ref{5.4}) 
and letting $ k \lra \infty$
we would get $ 0 \geq E_{c, min }( k_0 ) > 0 $, a contradiction. 
Thus (\ref{5.3}) is proven.

We have $B_c(u_n) = P_c( u_n) -  \frac{N-3}{N-1} A(u_n) $ and using 
(b) and (\ref{5.3}) we obtain  
\beq
\label{5.5}
 \ds \limsup_{n \ra \infty } B_c( u_n ) = -  \frac{N-3}{N-1}  C_3 < 0.
\eeq 
Clearly, for any $ \si > 0$ we have
$$
P_c( (u_n)_{1, \si })  
= \si  ^{N-3}\frac{N-3}{N-1}A(u_n) + \si ^{N-1} B_c( u_n)  
= \si ^{N-3} \left( \frac{N-3}{N-1}A(u_n) + \si ^2 B_c( u_n)  \right).
$$
For $ n $ sufficiently large (so that $B_c(u_n) < 0$), let 
$ \tilde{\si } _n = \left( \frac{\frac{N-3}{N-1}A(u_n) }{-  B_c( u_n) }\right)^{\frac 12}$. 
Then  $P_c( (u_n)_{1, \tilde{\si } _n}) = 0$, or equivalently 
$ (u_n)_{1, \tilde{\si } _n} \in \Co $.
From Lemma \ref{L4.7} we obtain 
$$
E_c( (u_n)_{1, \tilde{\si } _n} ) = 
\tilde{\si }_n ^{N-3} A(u_n) + \tilde{ \si }_n ^{N-1} B_c( u_n)   
\geq T_c,
\qquad \mbox{  that is}
$$ 
\beq
\label{5.6}
E_c( u_n) + \left( \tilde{\si }_n  ^{N-3} - 1 \right) A(u_n) 
+ \left( \tilde{\si }_n  ^{N-1} - 1 \right)
\left( P_c ( u_n) - \frac{N-3}{N-1}  A(u_n) \right) \geq T_c.
\eeq
Clearly, $ \tilde{ \si }_n$ can be written as 
$ \tilde{\si } _n = \left( \frac{ P_c (u_n)}{ - B_c (u_n)} + 1 \right)^{\frac 12}$
and using  (b) and (\ref{5.5}) it follows that 
$\ds \lim_{ n \ra \infty } \tilde{ \si }_n = 1$. 
Then passing to the limit as $ n \lra \infty $ in (\ref{5.6}) 
and using the fact that $ ( A(u_n) )_{n \geq 1}$ and 
$ ( P_c(u_n) )_{n \geq 1}$ are bounded, we obtain 
$\ds \liminf_{n \ra \infty } E_c(u_n) \geq T_c$. 
\hfill
$\Box $ 

\medskip

We can now state the main result of this section.

\begin{Theorem}
\label{T5.3}
Let $(u_n)_{n \geq 1} \subset \Xo \setminus \{ 0 \}$ be a sequence 
such that 
$$
P_c( u_n ) \lra 0 \qquad \mbox{ and } \qquad E_c( u_n) \lra T_c 
\qquad \mbox{ as } n \lra \infty. 
$$
There exist a subsequence $(u_{n_k})_{k \geq 1}$, 
a sequence $(x_k)_{k \geq 1} \subset \R^N$  and $ u \in \Co $ such that 
$$
\nabla u_{n_k}(\cdot + x_k)  \lra \nabla u \quad \mbox{ and } \; \;  
|1 + u_{n_k}( \cdot + x_k)| - 1 \lra |1+ u| - 1
\; \;  \mbox{ in } L^2 ( \R^N) .
$$
Moreover, we have $E_c(u) = T_c$, that is $u$ minimizes $E_c$ in $\Co$. 
\end{Theorem}

{\it Proof. } 
From Lemma \ref{L5.1} we know that $E_{GL}(u_n)$ is bounded. 
We have $ \frac{2}{N-1} A(u_n) = E_c( u_n) - P_c( u_n) \lra T_c $ 
as $ n \lra \infty$. Therefore
\beq
\label{5.6bis}
\lim_{ n \ra \infty } A(u_n) = \frac{N-1}{2} T_c \; \mbox{ and } \;
\liminf_{ n \ra \infty } E_{GL}( u_n) 
\geq \lim_{ n \ra \infty } A(u_n) = \frac{N-1}{2} T_c.
\eeq
Passing to a subsequence if necessary, we may assume that there exists
$ \al _0 \geq \frac{N-1}{2} T_c $ such that 
\beq
\label{5.7}
E_{GL}(u_n) \lra \al _0 
\qquad \mbox{ as } n \lra \infty. 
\eeq

We will use the concentration-compactness principle (\cite{lions}).
We denote by $q_n(t)$ the concentration function of $E_{GL}(u_n)$, that is
\beq
\label{5.8}
q_n(t ) = \sup_{y \in \R^N} \int_{B(y, t)}
\left\{ |\nabla u_n |^2 + \frac 12 \left( \ph ^2( |1+ u_n |) - 1 \right)^2 \right\}  dx. 
\eeq
As in \cite{lions}, it follows that there exists a subsequence of 
$((u_n, q_n))_{n \geq 1}$, still denoted $((u_n, q_n))_{n \geq 1}$, 
 there exists a nondecreasing function  $q : [0, \infty ) \lra \R$
and there is $\al \in [0, \al _0] $ such that 
\beq
\label{5.9}
q_n (t) \lra q(t)\; \mbox{  a.e on }  [0, \infty ) \; \mbox{ as } n \lra \infty
\qquad \mbox{ and } \qquad
q(t) \lra \al \mbox{ as } t \lra \infty.
\eeq

We claim that 
\beq
\label{5.10}
\mbox{ there is a nondecreasing sequence } t_n \lra \infty 
\mbox{ such that } \lim_{n \ra \infty }q_n (t_n) = \al .
\eeq
To prove the claim, fix an increasing sequence $ x_k \lra \infty $ such that 
$ q_n ( x_k) \lra q(x_k)$ as $ n \lra \infty $ for any $ k$. 
Then there exists $ n_k \in \N$ such that $| q_n ( x_k) - q(x_k)| < \frac 1k$
for any $ n \geq n_k$; 
clearly, we may assume that $ n_k < n_{k+1}$ for all $k$. 
If $ n_k \leq n < n_{k+1}$, put $ t_n = x_k$. 
Then for $ n_k \leq n < n_{k+1}$ we have
$$
| q_n ( t_n ) - \al | = | q_n ( x_k) - \al | 
\leq |q_n ( x_k) - q( x_k) | + |q(x_k) - \al | 
\leq \frac 1k + |q(x_k) - \al |  \lra 0 
$$
as $ k \lra \infty $ and (\ref{5.10}) is proved. 

Next we claim that 
\beq
\label{5.11} 
q_n ( t_n) - q_n \left( \frac{t_n}{2} \right) \lra 0 \quad \mbox{ as } n \lra \infty.
\eeq
To see this, fix $ \e > 0$. 
Take $ y > 0$ such that $ q(y) > \al - \frac{ \e}{4}$
and $q_n (y) \lra q(y)$ as $ n \lra \infty$. 
There is some $ \tilde{n}\geq 1$ such that $q_n(y) > \al - \frac{\e}{2}$ for 
$ n \geq \tilde{n}$. 
Then we can find $ n_* \geq \tilde{n}$ such that $ t_n > 2 y$ for $ n \geq n_*$, 
and consequently we have 
$ q_n ( \frac{t_n}{2}) \geq q_n (y) > \al - \frac{\e}{2}$. 
Therefore ${\ds \limsup_{n \ra \infty }} \left( q_n ( t_n) - q_n ( \frac{t_n}{2}) \right) 
= {\ds \lim_{ n \ra \infty } } q_n ( t_n) 
- {\ds \liminf_{ n \ra \infty } } \, q_n ( \frac{t_n}{2}) <\e.$
Since $ \e $ was arbitrary, (\ref{5.11}) follows. 

\medskip

Our aim is to show that $ \al = \al _0 $ in (\ref{5.9}).
It follows from the next lemma that $ \al > 0$.

\begin{Lemma}
\label{L5.4}
Let $(u_n) _{n \geq 1} \subset \Xo $ be a sequence satisfying

a) $ M_1 \leq E_{GL}(u_n) \leq M_2 $ for some positive constants $M_1$, $M_2$. 

b) $\ds \lim_{n \ra \infty } P_c( u_n) = 0$. 

There exists $ k > 0$ such that 
$ \ds \sup_{ y \in \R^N } \int_{B(y, 1)} 
\left\{ |\nabla u_n |^2 + \frac 12 \left( \ph ^2( |1+ u_n |) - 1 \right)^2 \right\} dx \geq k$
for all sufficiently large $n$. 
\end{Lemma}

{\it Proof. } We argue by contradiction and we  suppose that the conclusion is false. 
Then there exists a subsequence (still denoted $(u_n) _{n \geq 1} $) such that 
\beq
\label{5.12}
\lim_ {n \ra \infty } \sup_{ y \in \R^N } \int_{B(y, 1)} 
\left\{ |\nabla u_n |^2 + \frac 12  \left( \ph ^2( |1+ u_n |) - 1 \right)^2 \right\}  dx  =0.
\eeq
In order to get a contradiction we proceed in four steps. 

\medskip

{\it Step 1.  We show that $|E(u_n) - E_{GL}(u_n)| \lra 0 $ as $ n \lra \infty$. 
More precisely, we  prove that  }
\beq
\label{5.13}  
\lim_{ n \ra \infty} \int_{\R^N} 
\Big\vert 
V( |1+ u_n |^2) - \frac 12 \left( \ph ^2( |1+ u_n |) - 1 \right)^2 
\Big\vert  \, dx = 0.
\eeq
Fix $ \e > 0$. 
Assumptions (A1) and (A2) imply that there exists $ \de (\e ) > 0$ such that 
\beq
\label{5.14}   
\Big\vert  V(| 1+ z|^2) - \frac 12 \left( \ph ^2( |1+ z|) - 1 \right)^2  \Big\vert 
\leq \frac{\e}{2}   \left( \ph ^2( |1+ z|) - 1 \right)^2 
\eeq
for any $ z \in \C$ satisfying $\big\vert \, | 1+ z | - 1 \big\vert \leq \de ( \e )$
(see (\ref{4.2})). 
Therefore
\beq
\label{5.15}
\begin{array}{l}   
\ds \int_{ \{ | \, | 1+ u_n | - 1 | \leq \de ( \e ) \} }
\Big\vert 
V( |1+ u_n |^2) - \frac 12 \left( \ph ^2( |1+ u_n |) - 1 \right)^2 
\Big\vert  \, dx 
\\
\\
\ds \leq \frac{\e}{2}  \int_{ \{ | \, | 1+ u_n | - 1 | \leq \de ( \e ) \} }
\left( \ph ^2( |1+ u_n |) - 1 \right)^2  \, dx \leq \e M_ 2. 
\end{array}
\eeq
Assumption (A2) implies that there exists $ C(\e) > 0$ such that 
\beq
\label{5.16} 
\Big\vert 
V( |1+ z |^2) - \frac 12 \left( \ph ^2( |1+ z|) - 1 \right)^2 
\Big\vert  
\leq C( \e) | \, | 1+ z| - 1 |^{2 p_0 +2 } 
\eeq
for any $ z \in \C$ verifying $\big\vert \, | 1+ z | - 1 \big\vert \geq \de ( \e )$.

Let $w_n = | \, |1+u_n | - 1 |$. It is clear that $|w_n| \leq | u_n|$. 
Using the inequality $| \nabla | v| \, | \leq |\nabla v|$ a.e. for  
$ v \in H_{loc}^1(\R^N)$, we infer that $ w_n \in \DR $ and 
\beq
\label{5.17} 
\int_{\R^N} |\nabla w _n |^2 \, dx 
\leq M_2  \qquad \mbox{ for any } n .
\eeq

Using (\ref{5.16}), H\"older's inequality,  the Sobolev embedding and (\ref{5.17}) we find
\beq
\label{5.18}  
\begin{array}{l}
\ds \int_{\{| 1+ u_n | - 1 | > \de ( \e ) \} }
\Big\vert 
V( |1+ u_n |^2) - \frac 12  \left( \ph ^2( |1+ u_n |) - 1 \right)^2 
\Big\vert  \, dx 
\\
\\
\leq C( \e ) \ds \int_{\{ w_n > \de (\e ) \} } | w_n |^{ 2 p _0 +2} \, dx 
\\
\\
\leq C( \e ) \left( \ds \int_{\{ w_n > \de (\e ) \} } | w_n |^{ 2 ^* } \, dx
     \right)^{\frac{ 2 p_0 +2 }{2^*}}
\left( \Lo ^ N (  \{ w_n > \de (\e ) \} ) \right) ^{ 1 - \frac{ 2 p_0 +2 }{2^*} }
\\
\\
\leq C( \e )  C_S^{ 2 p_0 +2} \| \nabla w_n \|_{ L^2 ( \R^N)} ^{ 2 p _ 0 +2} 
\left( \Lo ^ N (  \{ w_n > \de (\e ) \} ) \right) ^{ 1 - \frac{ 2 p_0  +2}{2^*} } 
\\
\\
\leq C( \e )  C_S^{ 2 p_0 +2}  M_2^{p_0 +1} 
\left( \Lo ^ N (  \{ w_n > \de (\e ) \} ) \right) ^{ 1 - \frac{ 2 p_0 +2 }{2^*} } .
\end{array}
\eeq

We claim that for any $ \de > 0$ we have
\beq
\label{5.19}
\lim_{ n \ra \infty } 
\Lo ^ N (  \{ w_n > \de  \} )  = 0.
\eeq
To prove the claim, we argue by contradiction and we assume that there exist $ \de _ 0 > 0$, 
a subsequence $(w_{n _k})_{k \geq 1 } $ and $ \gamma > 0$ such that 
$ 
\Lo ^ N \left(  \{ w_{n _k} > \de _0 \} \right)  \geq \gamma > 0
$ for any $ k \geq 1$. 
Since $\| \nabla w_n \|_{L^2( \R^N)} $ is bounded, using Lieb's lemma 
(see Lemma 6 p. 447 in \cite{lieb} or Lemma 2.2 p. 101 in \cite{brezis-lieb}), we infer that there exists 
$ \beta > 0$ and  
$ y_k \in \R^N$ such that 
$
\Lo ^ N \left(  \{ w_{n _k} > \frac{\de _0 }{2} \} \cap B( y_k , 1)  \right)  \geq  \beta.
$
Let $ \eta $ be as in (\ref{3.30}). 
Then $ w_{n_k}( x) \geq \frac{ \de _0 }{2} $ implies 
$\left( \ph ^2( |1+ u_{n _k} (x)|) - 1 \right)^2 
\geq \eta\left( \frac{ \de _0}{2}  \right) >0$.
Therefore
$$
\ds \int_{B(y_k, 1)} \left( \ph ^2( | 1+ u_{n _k} (x)|) - 1 \right)^2  \, dx
\geq \eta\left( \frac{\de _0 }{2}  \right)  \beta > 0
$$
for any $k \geq 1$, and this clearly contradicts (\ref{5.12}). 
Thus we have proved that (\ref{5.19}) holds.

From (\ref{5.15}), (\ref{5.18}) and (\ref{5.19}) it follows that 
$$
\int_{\R^N} 
\Big\vert 
V( |1+ u_n |^2) - \frac 12  \left( \ph ^2( |1+ u_n |) - 1 \right)^2 
\Big\vert  \, dx
\leq 2 \e M_2
$$
for all sufficiently large $n$. Thus (\ref{5.13}) holds and the proof of step 1 is complete.

\medskip

{\it Step 2.  We find a convenient scaling of $u_n$. }
From Lemma \ref{L5.2} we know that $\ds \liminf_{ n \ra \infty } E_{c}(u_n) \geq T_c$. 
Combined with (b), this implies ${ \ds \liminf_{ n \ra \infty }} \frac{2}{N-1} A(u_n) \geq T_c$. 
Let $ \si _0 = \sqrt{ \frac{ 2(N-1)}{N-3}} $ and let $ \tilde{u}_n = (u_n)_{1, \si _ 0}$. 
It is obvious that 
\beq
\label{5.20}
\liminf_{ n \ra \infty } A(\tilde{u}_n) = 
\si _0 ^{N-3} \liminf_{ n \ra \infty } A(u_n)
\geq \frac{N-1}{2} \si _0 ^{ N-3} T_c.
\eeq
Using assumption (a), (\ref{5.12}) and (\ref{5.13}) it is easy to see that 
\beq
\label{5.21}
\mbox{ there exist } \tilde{M}_1, \, \tilde{M}_ 2 > 0 \mbox{ such that }
\tilde{M}_1 \leq E_{GL} ( \tilde{u} _n ) \leq \tilde{M}_2 
\mbox{ for any } n , 
\eeq
\beq
\label{5.22} 
\lim_ {n \ra \infty } \sup_{ y \in \R^N } \int_{B(y, 1)} 
\left\{ |\nabla \tilde{u}_n |^2 + \frac 12 \left( \ph ^2( |1+ \tilde{u}_n |) - 1 \right)^2 \right\} dx  =0 
\quad \mbox{ and } 
\eeq
\beq
\label{5.23}  
\lim_{ n \ra \infty} \int_{\R^N} 
\Big\vert 
V( |1+ \tilde{u}_n |^2) - \frac 12 \left( \ph ^2( | 1+ \tilde{u}_n |) - 1 \right)^2 
\Big\vert  \, dx = 0.
\eeq  
It is clear that $ P_c ( u_n) = \frac{N-3}{N-1} \si _0 ^{ 3 -N} A( \tilde{u}_n) + 
\si _0 ^{1-N} B_c( \tilde{u}_n) $ and then assumption (b) implies
\beq
\label{5.24}  
\lim_{ n \ra \infty} 
\left( \frac{N-3}{N-1} \si _0 ^ 2 A( \tilde{u}_n) + B_c( \tilde{u}_n)  \right)
= \lim_{ n \ra \infty}  \left( A( \tilde{u}_n) + E_c( \tilde{u}_n)  \right) = 0. 
\eeq

{\it Step 3.  Regularization of $\tilde{u}_n$. }
Using (\ref{5.21}), (\ref{5.22}) and Lemma \ref{vanishing} we infer that there is 
a sequence $ h_n \lra 0$ and for each $n$ there exists a minimizer 
$v_n$ of $G_{h_n, \R^N}^{\tilde{u}_n }$ in $H_{\tilde{u}_n} ^1( \R^N)$ such that 
$\de _n : =  \| \, |1 + v_n  | - 1 \|_{L^{\infty }(\R^N)} \lra 0 $ as $ n \lra \infty$. 
Then using Lemma \ref{L4.2} and the fact that $ |c| < v_s = \sqrt{2} $ we obtain 
\beq
\label{5.25}   
E_{GL}( v_n ) + c Q( v_n) \geq 0 
\qquad \mbox{ for all  sufficiently large $n$. }
\eeq

From (\ref{5.21}) and   (\ref{3.4}) we get
\beq
\label{5.27}    
| Q(\tilde{u}_n) - Q( v_n) | 
\leq C \left( h_n ^2 + h_n ^{\frac 4N} \tilde{M}_2^{\frac 2N} \right)^{\frac 12}
\tilde{M}_2 \lra 0 
\qquad \mbox{ as } n \lra \infty . 
\eeq

{\it Step 4.  Conclusion. }
Since $ E_{GL}(v_n) \leq E_{GL}(\tilde{u}_n)$, it is clear that 
$$
\begin{array}{rcl}
E_c( \tilde{u}_n ) & = & 
E_{GL}(\tilde{u}_n ) + cQ( \tilde{u}_n ) + 
\ds \int_{\R^N} \left\{ V(| 1+ \tilde{u}_n|^2) 
- \frac 12 \left( \ph ^2( |1+ \tilde{u}_n |) - 1  \right)^2  \right\} dx 
\\
\\
 &  \geq &
E_{GL}(v_n) + cQ(v_n) + c (Q(\tilde{u}_n) - Q( v_n) )
\\
& &  \qquad \qquad \qquad 
- \ds \int_{\R^N} \Big\vert V(|1+ \tilde{u}_n|^2) 
- \frac 12  \left( \ph ^2( |1+ \tilde{u}_n |) - 1 \right)^2 \Big\vert \, dx 
\end{array}
$$
Using the last inequality and (\ref{5.23}), (\ref{5.25}), (\ref{5.27}) 
we infer that $\ds \liminf_{n \ra \infty } E_c( \tilde{u}_n ) \geq 0$. 
Combined with (\ref{5.24}), this gives $\ds \limsup_{n \ra \infty } A( \tilde{u}_n ) \leq 0$, 
which clearly contradicts (\ref{5.20}). 
This completes the proof of Lemma \ref{L5.4}.
\hfill
$\Box $

\medskip

Next we prove that we cannot have $ \al \in (0, \al _0)$. 
To do this we argue again by contradiction and we assume that 
$ 0 < \al < \al _0$. 
Let $ t_n$ be as in (\ref{5.10}) and let $ R_n= \frac{t_n}{2}$. 
For each $n \geq 1$, fix $ y_n \in \R^N$ such that 
$E_{GL}^{B(y_n, R_n)} ( u_n) \geq q_n ( R_n) - \frac 1n$. 
Using (\ref{5.11}),  we have 
\beq
\label{5.28}
\begin{array}{l}
\e _n : = \ds \int_{B(y_n, 2R_n) \setminus B(y_n, R_n)} |\nabla u_n |^2 
+ \frac 12  \left( \ph ^2( |1+ u_n |) - 1 \right)^2 \, dx 
\\
\\
\leq q_n ( 2R_n) - \left( q_n( R_n) - \frac 1n \right) \lra 0 
\mbox{ as } n \lra \infty . 
\end{array}
\eeq
After a translation, we may assume that $ y_n = 0$. 
Using Lemma \ref{splitting} with $ A = 2$, $ R = R_n$, $ \e = \e_n$, 
we infer that for all $n$  sufficiently large there exist 
two functions $ u_{n, 1}$, $ u_{n, 2}$ having the properties (i)-(vi) 
in Lemma \ref{splitting}. 

From Lemma \ref{splitting} (iii) and (iv) we get 
$| E_{GL}(u_n) -  E_{GL}(u_{n, 1}) -  E_{GL}(u_{n, 2}) | \leq C \e _n$, 
while Lemma \ref{splitting} (i) and (ii) implies 
$ E_{GL}(u_{n, 1}) \geq E_{GL}^{ B(0, R_n)}(u_{n}) > q_n ( R_n) - \frac 1n$, 
respectively
$ E_{GL}(u_{n, 2}) \geq E_{GL}^{\R^N \setminus  B(0, 2R_n) }(u_{n}) \geq E_{GL}(u_n)
- q_n (2 R_n) $. Taking into account (\ref{5.7}), (\ref{5.10}), (\ref{5.11}) and (\ref{5.28}), we infer that 
\beq
\label{5.29}
E_{GL}(u_{n, 1}) \lra \al  \qquad \mbox{ and } \qquad 
E_{GL}(u_{n, 2}) \lra \al _0 - \al  \qquad \mbox{ as } n \lra \infty.
\eeq

By (\ref{5.28}) and Lemma \ref{splitting} (iii)$-$(vi) we obtain 
\beq
\label{5.30} 
| A(u_n) - A(u_{n,1}) -  A(u_{n,2}) | \lra 0, 
\eeq
\beq
\label{5.31} 
| E_c (u_n) - E_c(u_{n,1}) -  E_c(u_{n,2}) | \lra 0, \qquad \mbox{ and }
\eeq
\beq
\label{5.32} 
| P_c (u_n) - P_c(u_{n,1}) -  P_c(u_{n,2}) | \lra 0 \qquad \mbox{ as } n \lra \infty. 
\eeq

From (\ref{5.32}) and the fact that $ P_c( u_n) \lra 0 $ we infer that 
$P_c(u_{n,1}) + P_c(u_{n,2}) \lra 0 $ as $ n \lra \infty$. 
Moreover, Lemmas \ref{L4.1},  \ref{L4.5} and \ref{L5.1} imply that the sequences 
$(P_c(u_{n, i}))_{ n \geq 1}$ and $(E_c(u_{n, i}))_{ n \geq 1}$  are bounded for  $i= 1, 2$. 
Passing again to a subsequence (still denoted $(u_n)_{n \geq 1}$), we may assume that 
$\ds \lim_{ n \ra \infty } P_c(u_{n,1}) = p_1 $ and 
$\ds \lim_{ n \ra \infty } P_c(u_{n,2}) = p_2 $,   where  $ p_1, \, p_2 \in \R $ and  
$ p_1 + p_2 = 0.$
There are only two possibilities: either $ p_1 = p_2 = 0$, or 
one element of $\{ p_1, \, p_2 \}$ is negative. 

If $ p_1 = p_2 = 0$, then (\ref{5.29}) and Lemma \ref{L5.2} imply that 
$ \ds \liminf _{ n \ra \infty } E_c(u_{n, i}) \geq T_c$ for  $ i =1, \, 2$. 
Using (\ref{5.31}), we obtain $\ds \liminf _{ n \ra \infty } E_c(u_n) \geq 2 T_c$
and this clearly contradicts the assumption $E_{c} ( u_n) \lra T_c$ in 
Theorem \ref{T5.3}.

If $ p_i < 0$, it follows from (\ref{5.29}) and  Lemma \ref{L4.8} (ii) that 
 $ {\ds \liminf _{ n \ra \infty } } A(u_{n , i }) > \frac{N-1}{2}T_c$. 
Using (\ref{5.30}) and the fact that $ A \geq 0$, we obtain 
${ \ds \liminf _{ n \ra \infty } } A(u_n  ) > \frac{N-1}{2}T_c$, 
which is in  contradiction with (\ref{5.6bis}).

We conclude that we cannot have $ \al \in (0, \al _0)$.

\medskip

So far we have proved that $\ds \lim_{t \ra \infty } q(t) = \al _0$. 
Proceeding as in \cite{lions}, it follows that for each $n\geq 1$ there exists 
$ x_n \in \R^N$ such that for any $ \e > 0$ there are $ R_{\e} > 0$ and $ n_{\e } \in \N$ satisfying
\beq
\label{5.35}
E_{GL} ^{B(x_n, R_{\e})} ( u_n) > \al _0 - \e \qquad \mbox{ for any } n \geq n_{\e}. 
\eeq

Let $ \tilde{u}_n = u_n ( \cdot + x_n)$, so that $\tilde{u}_n $ satisfies (\ref{5.35}) with 
$B(0, R_{\e})$ instead of $ B(x_n, R_{\e})$. 
Let $ \chi \in C_c^{\infty } (\C, \R)$ be as in Lemma \ref{L2.2}
and denote $\tilde{u}_{n, 1} = \chi ( \tilde{u}_n ) \tilde{u}_n$, 
$\tilde{u}_{n, 1} =( 1-  \chi ( \tilde{u}_n )) \tilde{u}_n$.
Since $E_{GL}(  \tilde{u}_n  ) = E_{GL}( u_n) $ is bounded, 
we infer from Lemma \ref{L2.2}
that $(\tilde{u}_{n, 1} )_{n\geq 1}$ is bounded in $ \DR$, 
$(\tilde{u}_{n, 2} )_{n\geq 1}$ is bounded in $ H^1( \R^N)$ 
and $(E_{GL}( \tilde{u}_{n, i}))_{n \geq 1}$ is bounded for $ i =1, \, 2$.

Using Lemma \ref{lifting} we may write $ 1+ \tilde{u}_{n, 1}  = \rho _n e^{ i \theta _n}$, 
where $ \frac 12  \leq \rho _n \leq \frac 32  $ and $ \theta _n \in \DR$. 
From (\ref{2.4}) and (\ref{2.7}) we find that 
$(\rho _n - 1 )_{n \geq 1}$ is bounded in $ H^1( \R^N)$  and
$(\theta_n)_{n \geq 1}$ is bounded in $\DR$.

We infer that there exists a subsequence $ (n_k)_{k \geq 1}$ 
and there are functions $ u_1 \in \DR$, $ u_2 \in H^1( \R^N)$, 
$ \te \in \DR$, $ \rho \in 1 + H^1( \R^N)$   such that 
$$
 \tilde{u}_{n_k, 1}  \rightharpoonup u_1 \qquad \mbox{ and  } \qquad
 \theta_{n_k} \rightharpoonup \theta \qquad \mbox{ weakly in } \DR, 
$$
$$
\tilde{u}_{n_k, 2} \rightharpoonup u_2 \qquad \mbox{ and } \qquad
\rho_{n_k} - 1  \rightharpoonup \rho - 1  \qquad \mbox{ weakly in } H^1( \R^N),  
$$
$$
\tilde{u}_{n_k, 1} \lra u_1, \qquad 
\tilde{u}_{n_k, 2} \lra u_2, \qquad 
\theta_{n_k} \lra \theta , \qquad
\rho_{n_k} - 1  \lra \rho - 1
$$
strongly in $L^p(K)$, $ 1 \leq p < 2^*$ for any compact set $ K \subset \R^N$
and almost everywhere on $ \R^N$. 
Since 
$ \tilde{u} _{n _k, 1} =  \rho_{n_k} e^{ i \theta_{n_k}} -1 \lra  \rho e^{ i \theta} -1$
a.e., we have $  u_1 = \rho e^{ i \theta} -1 $ a.e. on $ \R^N$.

Denoting $ u = u_1 + u_2$, we see that $ \tilde{u} _{n _k} \rightharpoonup  u$ 
weakly in $\DR$, $ \tilde{u} _{n _k} \lra u $ a.e. on $ \R^N$ and strongly 
in $L^p(K)$, $ 1 \leq p < 2^*$ for any compact set $ K \subset \R^N$.

The weak convergence $ \tilde{u} _{n _k} \rightharpoonup  u$  in $\DR$ implies 
\beq
\label{5.37}  
\ds \int_{\R^N} \Big\vert \frac{\p u}{\p x_j } \Big\vert ^ 2 \, dx 
\leq
\liminf_{ k \ra \infty }
\ds \int_{\R^N} \Big\vert \frac{\p \tilde{u}_{n_k}}{\p x_j } \Big\vert ^ 2 \, dx 
< \infty \qquad \mbox{ for } j=1, \dots, N.
\eeq
Using the a.e. convergence $ \tilde{u} _{n _k} \lra  u$ and Fatou's lemma we obtain 
\beq
\label{5.38}  
\ds \int_{\R^N}  \left(\ph ^2 ( | 1+ u |)  - 1 \right)^2 \, dx 
\leq 
\liminf_{ k \ra \infty } 
\ds \int_{\R^N} \left(\ph ^2 ( | 1+ \tilde{u} _{n _k} |)  - 1 \right)^2 \, dx 
\eeq
From (\ref{5.37}) and(\ref{5.38}) it follows that $ u \in \Xo $ and 
$E_{GL}(u) \leq \ds \liminf_{ k \ra \infty }  E_{GL}( \tilde{u} _{n _k} )$.

We will prove that 
\beq
\label{5.39}   
\lim_{ k \ra \infty } \int_ {\R^N} V(|1+ \tilde{u}_{n_k}|^2) \, dx 
= \int_ {\R^N} V(| 1 +  u|^2) \, dx , 
\eeq
\beq
\label{5.39bis}
\lim_{ k \ra \infty } \| \, |1 + \tilde{u}_{n _k}| - |1 + u | \, \| _{L^2( \R^N)} = 0 \qquad \mbox{ and}
\eeq 
\beq
\label{5.40}   
\lim_{ k \ra \infty }  
Q(\tilde{u}_{n_k} ) = Q(u).
\eeq

Fix $ \e > 0$. Let $ R_{\e}$ be as in (\ref{5.35}).
 Since $ E_{GL}( \tilde{u}_{n_k}  ) \lra \al _0 $ as $ k \lra \infty$, 
it follows from (\ref{5.35})  that there exists $ k_{\e } \geq 1$ such that 
\beq
\label{5.41}   
E_{GL}^{ \R^N \setminus B(0, R_{\e})}( \tilde{u}_{n_k}  )  < 2 \e
\qquad \mbox{ for any } k \geq k_{\e }.
\eeq
As in (\ref{5.37})$-$(\ref{5.38}), the weak convergence 
$ \nabla \tilde{u}_{n_k} \rightharpoonup \nabla u $ in $L^2( \R^N \setminus B(0, R_{\e}))$
implies 
$$ \ds \int_{\R^N  \setminus B(0, R_{\e}) } |\nabla u | ^ 2 \, dx 
\leq
\liminf_{ k \ra \infty }
\ds \int_{\R^N  \setminus B(0, R_{\e}) } |\nabla \tilde{u}_{n_k} | ^ 2 \, dx ,
$$ 
while the fact that $ \tilde{u}_{n_k} \lra u $ a.e. on $ \R^N$ and Fatou's lemma imply
$$
  \ds \int_{\R^N  \setminus B(0, R_{\e}) }
 \left(\ph ^2 ( | 1+ u |)  - 1 \right)^2 \, dx  
 \leq 
\liminf_{ k \ra \infty }
\ds \int_{\R^N  \setminus B(0, R_{\e}) }  
\left(\ph ^2 ( | 1+ \tilde{u} _{n _k} |)  - 1 \right)^2 \,  dx. 
$$
Therefore 
\beq
\label{5.42}   
E_{GL}^{ \R^N \setminus B(0, R_{\e})}(u ) 
\leq 
\liminf_{ k \ra \infty }
E_{GL}^{ \R^N \setminus B(0, R_{\e})}( \tilde{u}_{n_k}  )  \leq 2 \e. 
\eeq

Let $ v \in \Xo $ be a function satisfying $ E_{GL}^{ \R^N \setminus B(0, R_{\e})}(v  ) \leq 2\e$.
Since $ \ph (| 1 + v |) = |1 + v |$ and $ \big| \, |1+ v | - 1 \big| ^2 \leq \left( \ph ^2 (| 1+ v|) - 1 \right)^2 $
 if $ |1 + v | \leq 2$, using (\ref{i1}) we find 
\beq
\label{5.43}    
\begin{array}{l}
\ds \int_{ \{ | 1 + v | \leq 2 \}   \setminus B(0, R_{\e}) } 
| V( | 1+ v|^2 ) | \, dx 
\leq C_1    {\ds \int_{ \{ | 1 + v | \leq 2 \} \setminus B(0, R_{\e}) } } 
\left(\ph ^2 ( | 1+ v  |)  - 1 \right)^2 \,  dx 
\\
\\
\leq  2 C_1  E_{GL}^{ \R^N \setminus B(0, R_{\e})}(v )  
\leq 4 C_1 \e \qquad \mbox{ and}
\end{array}
\eeq
\beq 
\label{5.43bis}
\int_{ \{ | 1 + v | \leq 2 \}   \setminus B(0, R_{\e}) } ( | 1 + v | - 1 )^2 \, dx 
\leq 
\int_{ \{ | 1 + v | \leq 2 \}   \setminus B(0, R_{\e}) } \left( \ph ^2 (| 1+ v|) - 1 \right)^2 \, dx
\leq 4 \e.
\eeq

On the other hand, $|1+ v(x) |  > 2  $ implies 
$ \left(\ph ^2 ( | 1+ v (x)  |)  - 1 \right)^2 > 9  $, consequently
$$
9  \Lo ^N (\{ x \in \R^N \setminus B(0, R_{\e}) \; | \; |1+ v(x) |  > 2  \})
\leq 
{\ds \int_{\R^N  \setminus B(0, R_{\e}) } } 
\left(\ph ^2 ( | 1+ v  |)  - 1 \right)^2 \,  dx  
\leq 4 \e.
$$
Using the fact that $ \big| V( |1+ s|^2 ) \big| \leq C\left( | 1+ s |^2 - 1 \right)^{ p_0 + 1} \leq C_2 |s |^{ 2 p_0 + 2} $ 
if $| 1 + s | \geq 2$,   H\"older's inequality,  the above estimate and the Sobolev embedding we find
\beq
\label{5.44}    
\begin{array}{l}
\ds \int_{ \{ | 1 + v | > 2 \}  \setminus B(0, R_{\e}) } |V ( |1 + v|^2) |\, dx 
\leq 
C _2 \ds \int_{ \{ | 1 + v | > 2 \}  \setminus B(0, R_{\e}) }
| v|^{ 2 p_0 + 2 } \, dx 
\\
\\
\leq C \left( \ds \int_
{\R^N} \! 
| v|^{ 2^*}  dx \! \right)^{\frac{2 p_0 + 2}{2 ^*}} 
\left( \Lo ^N (\{ x \in \R^N \setminus B(0, R_{\e}) \; \big| \; |1+ v(x) |  > 2  \}) 
     \right)^{ 1 - \frac{2 p_0 + 2}{2 ^*}} 
\\
\\
\leq C_3 \| \nabla v\|_{L^{2}(\R^N)} ^{ 2 p_0 + 2 } \e^{ 1 - \frac{2 p_0 + 2}{2 ^*}}
\leq C_3 \left( E_{GL}(v) \right) ^{ p_0 + 1 } \e^{ 1 - \frac{2 p_0 + 2}{2 ^*}}.
\end{array}
\eeq
Similarly we get 
\beq
\label{5.44bis}    
\begin{array}{l}
\ds \int_{ \{ | 1 + v | > 2 \}  \setminus B(0, R_{\e}) }   ( | 1 + v | - 1 )^2 \, dx 
\leq  \int_{ \{ | 1 + v | > 2 \}  \setminus B(0, R_{\e}) }  |v|^2 \, dx
\\
\\
\leq \left( \Lo ^N (\{ x \in \R^N \setminus B(0, R_{\e}) \; \big| \; |1+ v(x) |  > 2  \}) 
     \right)^{ 1 - \frac{2 }{2 ^*}} \| v \| _{L^{2^*} (\R^N)} ^2 
\\
\\
 \leq C \e^{ 1 - \frac{2 }{2 ^*}}  \| \nabla v \| _{L^{2} (\R^N)} ^2 
 \leq C E_{GL}(v)  \e^{ 1 - \frac{2 }{2 ^*}} .
\end{array}
\eeq

It is obvious that $ u$ and $ \tilde{u} _{n _k} $ (with  $ k \geq k_{\e}$) satisfy 
(\ref{5.43}) and (\ref{5.44}). 
If $M > 0$ is such that $E_{GL}(u_n) \leq M$ for all $n$, from (\ref{5.43}) and (\ref{5.44})  we infer that 
\beq
\label{5.45}  
\begin{array}{l}  
\ds \int_{\R^N  \setminus B(0, R_{\e}) } 
\big| V( | 1+ \tilde{u} _{n _k} |^2 ) - V( |1+ u |^2) \big| \, dx 
\\
\\
\leq 
\ds \int_{\R^N  \setminus B(0, R_{\e}) } 
\big| V( |1+ \tilde{u} _{n _k} |^2 ) \big| + \big| V( |1+ u |^2) \big| \, dx 
\leq 
C \e + C M^{p_0 + 1} \e ^{ 1 - \frac{2 p_0 + 2}{2 ^*}} , 
\end{array}
\eeq
while (\ref{5.43bis}) and (\ref{5.44bis})  give
\beq
\label{5.45bis}  
\| \, | 1+ \tilde{u}_{ n_k} | - 1 \|_{ L^2( \R^N \setminus B( 0, R_{\e}) )} ^2 \leq 4 \e +  C M \e^{ 1 - \frac{2 }{2 ^*}} .
\eeq
Of course,  a similar estimate is valid for $u$.

The mapping  $ z \longmapsto V(|1+ z|^2) $ is obviously $ C^1$. Since 
$| V(|1+z|^2)  | \leq C( 1 + |z |^{2 p_0 + 2} ) $ and 
$ \tilde{u} _{n _k}  \lra u $ in $L^2 \cap L^{2 p_0 + 2} ( B(0, R_{\e }))$
and almost everywhere, 
 it follows that $ | 1+ \tilde{u}_{ n_k} |  \lra |1 + u |$ in $L^2 ( B(0, R_{\e }))$  and
$  V( |1+ \tilde{u} _{n _k} |^2 ) \lra  V( |1+u |^2) $ in $L^1 ( B(0, R_{\e }))$  
(see, e.g., Theorem A2 p. 133 in \cite{willem}).
Hence
\beq
\label{5.46}   
\ds \int _{ B(0, R_{\e })} \big| V( | 1+ \tilde{u} _{n _k} |^2 )  - V( | 1+u |^2)  \big| \, dx 
\leq \e 
\qquad \mbox{ and}
\eeq
\beq
\label{5.46bis}   
\| \, | 1+ \tilde{u}_{ n_k} |  - |1 + u| \, \| _{L^2( B(0, R_{\e} ))} \leq \e  \qquad \mbox{ for all  $k$  sufficiently large.} 
\eeq
Since $ \e > 0$ is arbitrary,  (\ref{5.39}) follows from (\ref{5.45}) and (\ref{5.46}), 
while  (\ref{5.39bis}) is a consequence of (\ref{5.45bis}) and (\ref{5.46bis}).

Next we prove (\ref{5.40}). Fix $ \e > 0 $ and let $R_{\e}$ and $ k_{\e }$ be as in (\ref{5.35}) and (\ref{5.41}), respectively.
From (\ref{2.6}) we obtain 
$$
\|( 1 - \chi ^2 ( \tilde{u}_n) ) \tilde{u}_n \|_{L^2( \R^N)} 
\leq C\|\nabla \tilde{u}_n \| _{L^2( \R^N)}  ^{ \frac{2^*}{2}} 
\leq C \left( E_{GL}( u_n) \right)^{\frac{2^*}{4}}.
$$ 
Using the Cauchy-Schwarz inequality and (\ref{5.41}) we get 
\beq
\label{5.47}  
\begin{array}{l}
\ds \int _{ \R^N \setminus B(0, R_{\e })} 
\Big\vert 
( 1 - \chi ^2 ( \tilde{u} _{n _k} ) )
\langle i \frac{\p \tilde{u} _{n _k}}{\p x_1 }, \tilde{u} _{n _k} \rangle
\Big\vert  dx 
\\
\\
\leq
\|( 1 - \chi ^2 ( \tilde{u}_{n_k}) ) \tilde{u}_{n_k} \|_{L^2( \R^N)}  
\big\| \frac{\p \tilde{u} _{n _k}}{\p x_1  } \big\|_{L^2 (\R^N \setminus B(0, R_{\e }))}
\leq C M^{\frac{2^*}{4}} \sqrt{ \e } \quad \mbox{ for any } k \geq k_{\e }.
\end{array}
\eeq

From (\ref{2.7}) we infer that 
$$
 \|\rho_n ^2 - 1 \|_{L^2( \R^N)}   
\leq C \left(  E_{GL}( \tilde{u}_n) + \| \nabla \tilde{u}_n \|_{L^2( \R^N)}   ^{2^*} \right)^{\frac 12} 
\leq C\left( M + M  ^{\frac{2^*}{2}} \right)^{\frac 12}. 
$$
Using (\ref{2.4}) and (\ref{2.5}) we obtain 
$
\big\vert \frac{ \p \theta _n }{\p x_ 1 } \big\vert 
\leq 2 \big\vert \frac{ \p ( \chi (\tilde{u}_n) \tilde{u}_n ) }{\p x_ 1 } \big\vert  
\leq C \big\vert \frac{ \p  \tilde{u}_n  }{\p x_ 1 } \big\vert  
$
a.e. on $ \R^N$ and then (\ref{5.41}) implies 
$\big\| \frac{ \p \theta _{n _k}}{ \p x_ 1} \big\|_{L^2 (\R^N \setminus B(0, R_{\e }))} 
\leq C \sqrt{ \e}$ for any $ k \geq k_{\e}$ . 
Using again the Cauchy-Schwarz inequality we find
\beq
\label{5.48}  
\ds \int _{ \R^N \setminus B(0, R_{\e })} 
\Big\vert  
\left( \rho_{n_k} ^ 2 - 1 \right)  \frac{ \p \theta _{n _k}}{ \p x_ 1}
\Big\vert  \, dx 
\leq \|\rho_{n _k} ^2 - 1 \|_{L^2( \R^N)}    
\Big\|
 \frac{ \p \theta _{n _k}}{ \p x_ 1} 
\Big\| _{L^2 (\R^N \setminus B(0, R_{\e }))}  
\leq C \left( M \right) \sqrt{ \e } 
\eeq
for any $ k \geq k_{\e }.$
It is obvious that the estimates (\ref{5.47}) and (\ref{5.48}) also hold with $u$, $\rho $ and $ \theta $ instead of 
$\tilde{u} _{n _k} $, $\rho_{n_k}$ and $ \theta_{n_k}$, respectively. 

Using the fact that $ \tilde{u} _{n _k}  \lra u $ 
and $\rho_{n _k} - 1  \lra \rho - 1 $ 
in $L^2( B(0, R_{\e}))$ and a.e. 
and the dominated convergence theorem we infer that 
$$
( 1 - \chi ^2 ( \tilde{u} _{n _k} )) \tilde{u} _{n _k}  \lra ( 1 - \chi ^2 ( u)) u 
\qquad \mbox{ and } \qquad
\rho _{n_k}^2 - 1 \lra \rho ^2 - 1
\qquad \mbox{ in } L^2( B(0, R_{\e})).
$$
This information and the fact that 
$ \frac{ \p \tilde{u} _{n _k}}{ \p x_1} \rightharpoonup \frac {\p u }{\p x_1 }$ 
and 
$ \frac{ \p \theta _{n _k}}{ \p x_1} \rightharpoonup \frac {\p \theta }{\p x_1 }$ 
weakly in $L^2( B(0, R_{\e})) $ imply 
\beq
\label{5.49}
\ds \int_{ B(0, R_{\e}) } 
\langle i \frac{ \p \tilde{u} _{n _k}}{ \p x_1} , 
( 1 - \chi ^2 ( \tilde{u} _{n _k} )) \tilde{u} _{n _k} \rangle \, dx 
\lra 
\ds \int_{ B(0, R_{\e}) } 
\langle i \frac{ \p u}{ \p x_1} , 
( 1 - \chi ^2 ( u )) u \rangle \, dx 
\quad \mbox{ and } 
\eeq
\beq
\label{5.50}
\ds \int_{ B(0, R_{\e}) }  
\left( \rho _{n_k}^2 - 1  \right) \frac{ \p \theta _{n _k}}{ \p x_1} \, dx 
\lra 
\ds \int_{ B(0, R_{\e}) }  
\left( \rho ^2 - 1  \right) \frac{ \p \theta }{ \p x_1} \, dx  .
\eeq
Using (\ref{5.47})$-$(\ref{5.50}) 
and the representation formula (\ref{2.12}) we infer that 
 there is some $ k_1 ( \e ) \geq k_{\e } $ such that for any 
$ k \geq k_1 ( \e )  $ we have 
$$
| Q( \tilde{u}_{n_k} ) - Q(u) | \leq C    \sqrt{ \e } , 
$$
where $ C$ does not depend on $k \geq k_1 ( \e)$ and $ \e$.
Since $ \e > 0$ is arbitrary,  (\ref{5.40}) is proven.

Notice that the proofs of (\ref{5.39})-(\ref{5.40}) above are also valid if $N=3$. 

\medskip

It is obvious that 
$$
\begin{array}{l}
 -c Q(\tilde{u}_{n_k} )  - \ds \int_{\R^N} V( | 1+  \tilde{u}_{n_k} |^2) \, dx 
 \\ 
 \\
 = \ds \frac{N-3}{N-1} A(\tilde{u}_{n_k}) 
 + \int_{\R^N} \Big\vert \frac{ \p \tilde{u}_{n_k} }{\p x_1 } \Big\vert  ^2 dx 
 - P_c(\tilde{u}_{n_k}) 
\geq  \frac{N-3}{N-1} A(\tilde{u}_{n_k})  - P_c(\tilde{u}_{n_k}) .
\end{array}
$$
Passing to the limit as $ k \lra \infty $ 
in this inequality and using (\ref{5.39}), (\ref{5.40})
and the fact that $ A( \tilde{u}_n) \lra \frac{N-1}{2} T_c$, $P_c( \tilde{u}_n) \lra 0 $ as $ n \lra \infty $ 
we find
\beq
\label{5.51} 
- cQ(u) - \ds \int_{\R^N} V( | 1 +  u |^2) \, dx
\geq \frac{N-3}{2} T_c > 0.
\eeq
In particular, (\ref{5.51}) implies that $ u \neq 0$. 

From (\ref{5.37}) we get 
\beq
\label{5.52}  
A(u) \leq \liminf_ {k \ra \infty }  A(\tilde{u}_{n_k} ) = \frac{N-1}{2} T_c.
\eeq

Using (\ref{5.37}), (\ref{5.39}) and (\ref{5.40}) we find 
\beq
\label{5.53}
P_c ( u) \leq \liminf_ {k \ra \infty } P_c ( \tilde{u}_{n_k}) = 0.
\eeq
If $P_c(u) < 0$, from Lemma \ref{L4.8} (i) we get $A(u) > \frac{N-1}{2} T_c$, 
contradicting (\ref{5.52}). Thus necessarily $ P_c(u) =0$, that is $ u \in \Co$.
Since $ A(v) \geq \frac{N-1}{2} T_c $ for any $v \in \Co$, we infer from (\ref{5.52}) that 
$A(u) = \frac{N-1}{2} T_c$, therefore $E_c(u) = T_c$ and $u $ is a minimizer of $ E_c$ in $ \Co$.

It follows from the above that 
\beq
\label{5.54}
A(u) = \frac{N-1}{2} T_c = \lim_ {k \ra \infty }  A(\tilde{u}_{n_k} ).
\eeq
Since $ P_c(u) = 0$, $\ds \lim_{ k \ra \infty } P_c ( \tilde{u}_{n_k}) = 0$ and 
(\ref{5.39}), (\ref{5.40}) and (\ref{5.54})  hold, it is obvious that 
\beq
\label{5.55}
\int_{\R^N} 
\Big\vert \frac{ \p u }{\p x_1 } \Big\vert  ^2 dx 
= \lim_ {k \ra \infty }  
\int_{\R^N} \Big\vert \frac{ \p \tilde{u}_{n_k} }{\p x_1 } \Big\vert  ^2 dx .
\eeq
Now (\ref{5.54}) and (\ref{5.55}) imply 
$\ds  \lim_ {k \ra \infty }  \| \nabla \tilde{u}_{n_k} \|_{L^2(\R^N)} ^2 = 
\| \nabla u \|_{L^2(\R^N)} ^2$. 
Together with the fact that  $ \nabla \tilde{u}_{n_k} \rightharpoonup \nabla u $  weakly in $L^2(\R^N)$, 
this implies  $ \nabla \tilde{u}_{n_k} \lra \nabla u  $ strongly in $L^2(\R^N)$, that is 
$   \tilde{u}_{n_k} \lra  u  $ in $ \DR$ and  
the proof of Theorem \ref{T5.3} is complete. 
\hfill
$\Box $

\medskip

In order to prove that the minimizers provided by Theorem \ref{T5.3} solve equation (\ref{1.3}), 
we need the following regularity result. 

\begin{Lemma}
\label{L5.5}
Let $ N \geq 3$. 
Assume that the conditions (A1) and (A2) in the Introduction hold and  that 
$ u \in \Xo $ satisfies (\ref{1.3}) in $ \Do ' ( \R^N)$. 
Then $ u \in W_{loc }^{2, p } ( \R^N)$ for any $ p \in [1, \infty)$, 
$ \nabla u \in W^{1, p }(\R^N)$ for $ p \in [2, \infty)$, 
$ u \in C^{1, \al }(\R^N)$ for $ \al \in [0, 1)$ and $ u(x) \lra 0 $ as $ | x | \lra \infty$. 
\end{Lemma}

{\it Proof. } 
First we prove that for any $R>0$ and $ p \in [2, \infty)$ there exists $C(R, p) >0$
(depending on $u$, but not on $ x \in \R^N$) such that 
\beq
\label{reg1}
\| u \|_{W^{2, p}(B(x, R))} \leq C(R, p) 
\qquad \mbox{ for any } x \in \R^N. 
\eeq
We  write $ u = u_1 + u_2$, where $ u_1 $ and $ u_2$ are as in Lemma \ref{L2.2}. 
Then $ |u_1 | \leq \frac{ 1}{2}$, $ \nabla u_1 \in L^2 ( \R^N)$ and $ u_2 \in H^1( \R^N)$, hence 
for any  $R>0$ there exists $C(R) >0$ such that 
\beq
\label{reg2}
 \| u \|_{H^1(B(x, R))} \leq C(R) 
\qquad \mbox{ for any } x \in \R^N. 
\eeq
Let $ \phi (x) = e^{-\frac{i c x_1}{2}} (1 + u (x))$. 
It is easy to see that  $ \phi $  satisfies 
\beq
\label{reg3}
\Delta \phi + \left( F(|\phi |^2) + \frac{c^2}{4} \right) \phi = 0 
\qquad \mbox{ in } \Do '(\R^N).
\eeq
Moreover, (\ref{reg2}) holds for $ \phi $ instead of $u$. 
From (\ref{reg2}), (\ref{reg3}), (\ref{3.18}) and a standard bootstrap argument we infer that  
$ \phi $ satisfies (\ref{reg1}). 
(Note that  assumption (A2) is needed for this bootstrap argument.)
It is then clear that (\ref{reg1}) also holds for $u$. 

From (\ref{reg1}), the Sobolev embeddings and Morrey's inequality (\ref{3.27})
we find that $u$ and $ \nabla u$ are continuous and bounded on $ \R^N$ and 
$ u \in C^{1, \al }( \R^N)$ for $ \al \in [0, 1)$. 
In particular, $u$ is Lipschitz; since $ u \in L^{2^*}(\R^N)$, we  have necessarily
$ u(x) \lra 0 $ as $|x| \lra \infty$. 

The boundedness of $u$ implies that there is some $C>0$ such that 
$\big\vert  F( |1+u|^2) (1+ u) \big\vert \leq C \big\vert \ph ^2( |1+ u| ) -1 \big\vert$ on $ \R^N$. 
Therefore $  F( |1+ u|^2) (1+ u) \in L^2 \cap L^{\infty }(\R^N)$. 
Since $ \nabla u \in L^2 ( \R^N)$, from (\ref{1.3}) we find $ \Delta u \in L^2( \R^N)$. 
It is well known that $ \Delta u \in L^p( \R^N)$ with $1<p< \infty $ implies 
$\frac{ \p^2 u}{\p x_i \p x_j } \in L^p( \R^N)$ for any $ i, \, j$
(see, e.g., Theorem 3 p. 96 in \cite{stein}). 
Thus we get $ \nabla u \in W^{1, 2} ( \R^N)$. 
Then the Sobolev embedding implies $ \nabla u \in L^p( \R^N)$ for $ p \in [2, 2^*]$. 
Repeating the previous argument, after an easy induction we find 
$ \nabla u \in W^{1, p}(\R^N) $ for any $p\in [2, \infty)$.
\hfill
$\Box$

\begin{Proposition} 
\label{P5.6}
Assume that the conditions (A1) and (A2) in the introduction are satisfied. 
Let $ u \in \Co $ be a minimizer of $ E_c$ in $ \Co$. 
Then $ u \in W_{loc}^{2, p }(\R^N)$ for any $ p \in [1, \infty)$, 
$ \nabla u \in W^{1, p }(\R^N)$ for  $ p \in [2, \infty)$ 
and $u$ is a solution of  (\ref{1.3}).
\end{Proposition} 

{\it Proof. } 
It is standard to prove that for any $R>0$, 
$J_u(v) = \ds \int_{\R^N}  V(| 1+  u + v |^2) \, dx $ is a $C^1$ 
functional on $H_0 ^1(B(0, R))$ 
and $J_u '(v) . w = - 2  \ds \int_{\R^N}  F(| 1+ u + v |^2) \langle 1+ u + v, w \rangle\, dx $
(see, e.g., Lemma 17.1 p. 64 in \cite{kavian} or the appendix A in \cite{willem}). 
It follows easily that for any $R> 0$, the functionals 
$ \tilde{P}_c(v) = P_c( u+ v)$ and $ \tilde{E}_c(v) = E_c( u+ v)$ are 
$C^1$  on $H_0 ^1(B(0, R))$. The differentiability of $Q$ follows, for instance, from (\ref{2.18}).
We divide the proof of Proposition \ref{P5.6} into several steps.

\medskip

{\it Step 1. }  There exists a function $ w \in C_c ^1( \R^N)$ 
 such that $ \tilde{P}_c'(0) . w \neq 0$. 

\smallskip
 
To prove this, 
we argue by contradiction and we assume that the above statement is false. 
Then  $u$ satisfies
\beq
\label{5.57}
- \frac{ \p ^2 u}{ \p x_ 1 ^2 } - \frac{N-3}{N-1} 
\left( \sum_{ k = 2}^N \frac{ \p ^2 u}{ \p x_ k ^2 } \right)
+ i c u_{ x_ 1} - F(| 1+ u  |^2) ( 1+ u) = 0 
\qquad \mbox{ in } \Do ' (\R^N ) . 
\eeq
Let $ \si = \sqrt{ \frac{ N-1}{N-3}}$. 
It is not hard  to see that $ u_{1, \si }$ satisfies (\ref{1.3}) in $\Do ' (\R^N ) . $
Hence the conclusion of Lemma \ref{L5.5} holds for $ u_{1, \si }$ (and thus for $u$). 
This regularity is enough to prove that $u_{1, \si} $ satisfies  the Pohozaev identity
\beq
\label{5.58} 
\int_{\R^N} \! \Big\vert \frac{ \p u_{1, \si } }{ \p x_1} \Big\vert ^2 dx 
+ \frac{N-3}{N-1} \int_{\R^N}  \sum_{ k = 2}^N 
  \Big\vert \frac{ \p u_{1, \si} }{ \p x_k} \Big\vert ^2 dx 
+ cQ(u_{1, \si }) 
+ \int_{\R^N} \! \!  V(| 1+ u_{1, \si }  |^2) \, dx = 0.
\eeq
To prove (\ref{5.58}), we multiply (\ref{1.3}) by 
$\sum_{ k = 2}^N \tilde{\chi} ( \frac xn) \frac{ \p u_{1, \si}}{\p x_k}$,
where $ \tilde{\chi}\in C_c^{ \infty } ( \R^N)$ is a cut-off function such that
$ \tilde{\chi} = 1 $ on $B(0,1)$ and $\mbox{supp}(\tilde{\chi} ) \subset B(0,2)$, 
we integrate by parts, then we let $ n \lra \infty$; see the proof of Proposition 4.1 
and equation (4.13) 
in \cite{M8}, p.1094 for details.

Since $ \si = \sqrt{ \frac{ N-1}{N-3}}$, (\ref{5.58}) is equivalent to 
$ \left( \frac{ N-3}{N-1} \right)^2 A(u) + B_c (u) = 0$. 
On the other hand we have $ P_c (u) = \frac{ N-3}{N-1}  A(u) + B_c (u) = 0$
and we infer that $A(u) = 0$. 
But this contradicts the fact that $ A(u) = \frac{N-1}{2} T_c > 0$ and the proof of step 1 is complete. 

\medskip

{\it Step 2. } Existence of a Lagrange multiplier. 

\smallskip

Let $ w $ be as above and let $v\in H^1( \R^N)$ be a function with compact support
such that $ \tilde{P}_c'(0) .v =0$. 
For $ s, \, t \in \R$, put $ \Phi ( t, s ) = P_c( u + tv + s w) = \tilde{P}_c(tv + s w)$, 
so that $ \Phi (0,0) = 0$, $\frac { \p \Phi }{ \p t } (0, 0) = \tilde{P}_c'(0) .v =0$ and 
$\frac { \p \Phi }{ \p s } (0, 0) = \tilde{P}_c'(0) .w \neq 0$. 
The implicit function theorem implies that there exist $ \de > 0$ and  a $C^1 $ function 
$ \eta : (- \de , \de) \lra \R$ such that $ \eta (0) = 0$, $ \eta '(0) = 0$
and $P_c( u + tv + \eta (t) w ) = P_c ( u ) = 0$ for $ t \in (-\de, \de )$. 
Since $u$ is a minimizer of $A$ in $ \Co $, the function $ t \longmapsto A( u + tv + \eta (t) w )$
achieves a minimum at $ t =0$. 
Differentiating at $ t =0$ we get $ A'(u) . v = 0$. 

Hence $ A'(u) . v = 0$ for any $v\in H^1( \R^N)$ with compact support   
satisfying $ \tilde{P}_c'(0) .v =0$.  
Taking $ \al = \frac{ A'(u) . w}{ \tilde{P}_c'(0) . w}$ (where $w$ is as in step 1), 
we see that 
\beq
\label{5.59}  
A'(u) . v = \al P_c '(u) .v 
\quad \mbox{ for any } v\in H^1( \R^N) \mbox{ with compact support. } 
\eeq


{\it Step 3. } We have $ \al < 0$. 

\smallskip

To see this, we argue by contradiction. 
Suppose  that $ \al > 0$. Let $w$ be as in step 1. We may assume that $ P_c '(u) .w > 0$.
From (\ref{5.59}) we obtain $ A'(u) .w > 0$. 
Since $ A'(u) .w = {\ds \lim_{ t \ra 0 }} \frac{A( u + t w)- A(u) }{t} $ and
$ P_c '(u) .w = {\ds \lim_{ t \ra 0 }} \frac{P_c( u + t w) - P_c(u)}{t} $, 
we see that for $ t <0$, $t$ sufficiently close to $0$ we have 
$ u + tw\neq 0$, $P_c( u + t w) < P_c(u) = 0  $ and $A( u + t w) < A(u) = \frac{ N-1}{2} T_c$. 
But this contradicts Lemma \ref{L4.8} (i). 
Therefore $ \al \leq 0$. 

Assume that $ \al = 0$. Then (\ref{5.59}) implies 
\beq
\label{5.60}  
\int_{\R^N} \sum_{ k = 2}^N 
\langle \frac { \p u }{ \p x_j} , \frac { \p v }{ \p x_j} \rangle \, dx = 0 
\quad \mbox{ for any } v\in H^1( \R^N) \mbox{ with compact support. } 
\eeq
Let $ \tilde{\chi } \in C_c^{\infty } (\R^N)$ be such that $ \chi = 1$ on $ B(0,1)$ and 
$\mbox{ supp}(\tilde{\chi} ) \subset B(0,2)$. 
Put $ v_n (x) = \chi ( \frac xn) u(x)$, so that 
$ \nabla v_n (x) = \frac 1n \nabla \tilde{\chi} (\frac xn) u +\tilde{\chi} (\frac xn) \nabla u$. 
It is easy to see that $ \tilde{\chi} (\frac{ \cdot }{n}) \nabla u \lra \nabla u$ 
in $L^2( \R^N)$ and 
$ \frac 1n \nabla \tilde{\chi} (\frac{\cdot }{n}) u  \rightharpoonup 0 $ weakly in 
$L^2( \R^N)$. 
Replacing $v$ by $ v_n $ in (\ref{5.60}) and passing to the limit as $ n \lra \infty $ 
we get $A(u) = 0 $, which contradicts the fact that $ A(u) = \frac{N-1}{2} T_c$. 
Hence we cannot have $ \al = 0$. Thus necessarily $ \al < 0$. 

\medskip

{\it Step 4. } Conclusion. 

\smallskip

Since $ \al < 0$, it follows from (\ref{5.59}) that $u$ satisfies 
\beq
\label{5.61}   
- \frac{ \p ^2 u}{\p x_ 1 ^2 } - \left( \frac{N-3}{N-1} - \frac { 1}{ \al } \right) 
\sum_{ k = 2}^N  \frac { \p ^2 u}{\p x_ k ^2 }
+ i c u_{ x_1 } - F(|1+ u|^2) ( 1+ u) = 0 
\qquad
\mbox{ in } \Do '( \R^N). 
\eeq
Let $ \si _0 = \left( \frac{N-3}{N-1} - \frac { 1}{ \al } \right) ^{- \frac 12}$.
It is easy to see that $ u_{1, \si _0 } $ satisfies ({\ref{1.3}) in $\Do '( \R^N)$. 
Therefore the conclusion of Lemma \ref{L5.5} holds for $ u_{1, \si _0 } $ (and consequently for $u$). 
Then Proposition 4.1 in \cite{M8} implies that $ u_{1, \si _0 } $  satisfies the Pohozaev identity
$
\frac{N-3}{N-1} A( u_{1, \si _0}  ) + B_c ( u_{1, \si _0}  ) = 0, 
$
or equivalently $ \frac{N-3}{N-1} \si _0 ^{ N-3} A(u) + \si _0 ^{ N-1} B_c (u) = 0$, 
which implies
$$
\frac{N-3}{N-1} \left( \frac{N-3}{N-1} - \frac { 1}{ \al } \right) A(u) + B_c (u) = 0.
$$
On the other hand we have $ P_c(u) = \frac{N-3}{N-1} A(u)+ B_c (u) = 0$. 
Since $A(u) > 0$, we get $\frac{N-3}{N-1} - \frac { 1}{ \al } = 1 $.
Then coming back to (\ref{5.61}) we see that $u$ satisfies (\ref{1.3}). 
\hfill
$\Box $

\section{The case $ N =3$ }

This section is devoted to  the proof of Theorem \ref{T1.1} in space dimension $ N =3$. 
We only indicate the differences with respect to  the case $ N \geq 4$. 
Clearly, if $ N =3$ we have $ P_c = B_c$. For $ v \in \Xo $ we denote 
$$
D(v) = \int_{\R^3} \Big\vert \frac{ \p v }{\p x_1 } \Big\vert ^2 \, dx 
+  \frac 12 \int_{\R^3} \left( \ph^2(|1+ v |) - 1 \right)^2 \, dx .
$$
For any $ v \in \Xo $ and $ \si > 0$ we have 
\beq
\label{6.1}
A(v_{1, \si }) = A(v), \qquad 
B_c (v_{1, \si }) = \si^2 B_c (v) \qquad
\mbox{ and } \qquad
D (v_{1, \si }) = \si^2 D (v) .
\eeq

If $ N =3$ we cannot have  a result similar to Lemma \ref{L5.1}.
To see this consider $ u \in \Co$, so that $ B_c(u) = 0$. 
Using (\ref{6.1}) we see that $ u_{1, \si } \in \Co $ for any $ \si > 0$ and we have 
$E_c( u_{1, \si }) = A(u) + \si ^2 B_c (u) = A(u)$, 
while $E_{GL}(u_{1, \si } ) = A(u) + \si ^2 D(u) \lra \infty $ as $ \si \lra \infty$. 

However, for any $ u \in \Co $ there exists $ \si > 0 $ such that
$D(u_{1, \si }) = 1 $ (and obviously $ u_{1, \si } \in \Co$, $E_c( u_{1, \si }) = E_c(u)$).
Since $ \Co \neq \emptyset $ and $ T_c = \ds \inf \{ E_c(u) \; | \; u \in \Co \}$, 
we see that there exists a sequence $ (u_n)_{n \geq 1} \subset \Co $ such that 
\beq
\label{6.2}
D(u_n) = 1 \qquad \mbox{ and } \qquad E_c( u_n) = A( u_n) \lra T_c
\quad \mbox{ as } n \lra \infty. 
\eeq
In particular, (\ref{6.2}) implies $ E_{GL}(u_n) \lra T_c + 1 $ as $ n \lra \infty $. 

The following result is the equivalent  of Lemma \ref{L5.2} in the case $N=3$. 

\begin{Lemma}
\label{L6.1}
Let $ N =3$ and let $(u_n)_{n \geq 1} \subset \Xo $ be a sequence satisfying

\medskip

a) There exists $ C > 0$ such that $ D(u_n) \geq C$ for any $n$, and

\medskip

b) $B_c (u_n) \lra 0 $ as $ n \lra \infty$. 

\medskip

\noindent
Then $\ds \liminf_{n \ra \infty } E_c( u_n) = \liminf_{n \ra \infty } A( u_n) \geq S_c$, 
where $S_c$ is given by (\ref{4.22}).
\end{Lemma}

{\it Proof. }
It suffices to prove that for any $ k > 0$ there holds
\beq
\label{6.3}
\liminf_{n \ra \infty } A( u_n) \geq E_{c, min} (k). 
\eeq
Fix $ k > 0$. Let $ n \geq 1$. 
If $A( u_n) \geq k$, by Lemma \ref{L4.6} (iii) we have $A( u_n) \geq  k > E_{c, min} (k)$. 
If $ A(u_n) < k$, 
since $E_{GL}((u_n)_{1, \si }) = A( u_n) + \si ^2 D(u_n)$ we see that there exists $ \si _n >0$
such that $E_{GL}((u_n)_{1, \si _n})  = k$.
Obviously, we have $ \si _n ^2 D(u_n) < k$, hence $ \si _n ^2 \leq \frac kC$ by (a). 
It is clear that $E_c((u_n)_{1, \si _n}) = A(u_n) + \si _n ^2 B_c( u_n) \geq E_{c, min} (k) $, 
therefore
$ A( u_n) \geq E_{c, min} (k) - \si _n ^2 | B_c ( u_n) |  \geq E_{c, min} (k) - 
\frac kC | B_c ( u_n) |$.
Passing to the limit as $ n \lra \infty $ we obtain (\ref{6.3}). 
Since $k>0$ is arbitrary, Lemma \ref{6.1} is proven.
\hfill
$\Box $

\medskip

Let 
$$
\begin{array}{rcl}
\Lambda _c & = &
 \{ \la \in \R \; | \; \mbox{ there exists a sequence } (u_n)_{n \geq 1}  \subset \Xo 
\mbox{ such that } 
\\
& & 
D(u_n) \geq 1, \;  B_c ( u_n) \lra 0 \mbox{ and } A(u_n) \lra \la 
\mbox{ as } n \lra \infty \}. 
\end{array}
$$
Using a scaling argument, we see that 
$$
\begin{array}{rcl}
\Lambda _c & = &
 \{ \la \in \R \; | \; \mbox{ there exist a sequence } (u_n)_{n \geq 1} \subset \Xo \mbox{ and } C > 0
\mbox{ such that } 
\\
& & 
D(u_n) \geq C, \;  B_c ( u_n) \lra 0 \mbox{ and } A(u_n) \lra \la 
\mbox{ as } n \lra \infty \}. 
\end{array}
$$ 
Let $ \la _c = \inf \Lambda _c$. 
From (\ref{6.2}) it follows that $ T_c \in \Lambda _c$. 
It is standard to prove that $\Lambda _c  $ is closed in $ \R$, hence $ \la _c \in \Lambda _c $.
From Lemma \ref{6.1} we obtain 
\beq
\label{6.4}
S_c \leq \la _c \leq T_c .
\eeq
The main result of this section is as follows. 

\begin{Theorem} 
\label{T6.2}
Let $ N=3$ and let $(u_n)_{n \geq 1} \subset \Xo  $ be a sequence such that
\beq
\label{6.5} 
D(u_n) \lra 1, \quad B_c( u_n) \lra 0 \quad \mbox{ and } \quad
A( u_n) \lra \la _c \quad \mbox{ as } n \lra \infty . 
\eeq
There exist a subsequence $(u_{n_k})_{k \geq 1}$, 
a sequence $(x_k)_{k \geq 1} \subset \R^3$  and $ u \in \Co $ such that 
$$
\nabla u_{n_k}(\cdot + x_k)  \lra \nabla u \quad \mbox{ and } \; \;  
|1+ u_{n_k}( \cdot + x_k)|  -1 \lra |1+ u|  - 1
\; \;  \mbox{ in } L^2 ( \R^3) .
$$
Moreover, we have $E_c(u) = A(u) = T_c =  \la _c  $  
 and $u$ minimizes $E_c$ in $\Co$. 
\end{Theorem}

{\it Proof. } 
By (\ref{6.5}) we have $E_{GL}(u_n) = A( u_n) + D( u_n) \lra \la _c + 1 $ as $ n \lra \infty$. 
Let $ q_n(t) $ be the concentration function of $E_{GL}(u_n)$, as in (\ref{5.8}). 
Proceeding as in the proof of Theorem \ref{T5.3}, we infer that there exist a subsequence 
of $(u_n, q_n)_{n \geq 1}$, still denoted  $(u_n, q_n)_{n \geq 1}$, a nondecreasing function 
$ q :[0, \infty ) \lra [0, \infty )$ and $ \al \in [0, \la _c + 1 ]$ 
such that (\ref{5.9}) holds. 
We see also that there exists a sequence $t_n \lra \infty $ satisfying (\ref{5.10}) and (\ref{5.11}). 

Clearly, our aim is to prove that $ \al = \la _c +1$. 
The next result implies that $ \al >0$. 

\begin{Lemma}
\label{L6.3}
Assume that $ N =3$, $ 0 \leq c < v_s $ and 
 $(u_n)_{ n\geq 1} \subset \Xo $ is a sequence satisfying
$ D( u_n) \lra 1$, $B_c ( u_n) \lra 0 $ as $ n \lra \infty $ and 
$ \ds \sup_{n \geq 1} E_{GL}(u_n) = M < \infty$. 

There exists $ k > 0$ 
 such that 
$$ \ds \sup_{ y \in \R^3 } \int_{B(y, 1)} 
\left\{ |\nabla u_n |^2 + \frac 12 \left( \ph ^2( |1+ u_n |) - 1 \right)^2 \right\} \, dx \geq k
\qquad \mbox{ for all sufficiently large } n.
$$ 
\end{Lemma}

{\it Proof. } We argue by contradiction and assume that 
the conclusion of Lemma \ref{L6.3} is false. 
Then there exists a subsequence, still denoted $(u_n)_{ n\geq 1}$, such that 
\beq
\label{6.6}
\sup_{ y \in \R^3} E_{GL} ^{B(y, 1)} (u_n) \lra 0  \qquad \mbox{ as } n \lra \infty.
\eeq
Exactly as in Lemma \ref{L5.4} we prove that (\ref{5.13}) holds, that is 
\beq
\label{6.7}  
\lim_{ n \ra \infty} \int_{\R^3} 
\Big\vert 
V( |1+ u_n |^2) - \frac 12 \left( \ph ^2( |1+ u_n |) - 1 \right)^2 
\Big\vert  \, dx = 0.
\eeq
Using  (\ref{6.7}) and the assumptions of Lemma \ref{L6.3} we find
\beq
\label{6.8}
cQ(u_n) = B_c(u_n) - D( u_n) - 
\int_{\R^3} 
\left\{ V( |1+ u_n |^2) - \frac 12 \left( \ph ^2( |1+ u_n |) - 1 \right)^2 \right\} \, dx 
\lra -1 
\eeq
 as $n \lra \infty.$
If $ c =0$, (\ref{6.8}) gives a contradiction and Lemma \ref{L6.3} is proven. 
From now on we assume that $ 0< c< v_s $. 
 
Fix $ c_1 \in (c, v_s)$,  then fix $ \si > 0 $ such that 
\beq
\label{6.9} 
\si ^2 > \frac{ Mc}{ c_1 - c} .
\eeq
A simple change of variables shows that  
 $ \tilde{M} := \ds \sup_{ n \geq 1} E_{GL}( (u_n)_{1, \si }) < \infty $
and (\ref{6.7}) holds with  $(u_n)_{1, \si } $  instead of  $ u_n$. 
It is  easy to see that  $((u_n)_{1, \si }) _{ n\geq 1}$ also satisfies (\ref{6.6}).
Using Lemma \ref{vanishing} we infer that there exists a sequence $ h_n \lra 0$ and 
for each $n$ there exists 
a minimizer $ v_n $ of $ G_{h_n, \R^3} ^{ (u_n)_{1, \si } }$ in $ H_{(u_n)_{1, \si }}^1( \R^3)$ 
such that 
\beq
\label{6.10}  
\| \; |1+ v_n  | - 1 \| _{ L^{ \infty }( \R^3)} \lra 0 \qquad \mbox{ as } n \lra \infty .
\eeq
From   (\ref{3.4}) we obtain 
\beq
\label{6.11}    
| Q( (u_n)_{1, \si }) - Q( v_n) | 
\leq C \left( h_n ^2 + h_n ^{\frac 43} \tilde{M} ^{\frac 23} \right)^{\frac 12}
\tilde{M} \lra 0 
\qquad \mbox{ as } n \lra \infty . 
\eeq
Using (\ref{6.10}), the fact that $ 0 < c_1 < v_s $ and Lemma \ref{L4.2} we infer that for 
all sufficiently large $n$ there holds 
\beq
\label{6.12}   
E_{GL}( v_n) + c_1 Q(v_n ) \geq 0.
\eeq
Since $E_{GL}( v_n) \leq E_{GL}( (u_n)_{1, \si })$,    for large $n$ we have 
\beq
\label{6.13}
\begin{array}{l}
0 \leq E_{GL}( v_n) + c_1 Q(v_n ) 
\\
\\
\leq E_{GL}( (u_n)_{1, \si }) + c_1 Q( (u_n)_{1, \si }) + c_1| Q( (u_n)_{1, \si }) - Q( v_n) | 
\\
\\
= A( u_n ) + B_c ( (u_n)_{1, \si }) + (c_1 - c) Q( (u_n)_{1, \si })
+ c_1| Q( (u_n)_{1, \si }) - Q( v_n) |  
\\
\qquad + \ds \int_{\R^3}  \left\{ \frac 12  \left( \ph ^2( |1+ (u_n)_{1, \si} |) - 1 \right)^2 - 
V( |1+ (u_n)_{1, \si} |^2) \right\} \, dx 
\\
\\
= A( u_n) + \si ^2 B_c ( u_n) + \si ^2 ( c_1 - c) Q( u_n)  
+ a_n 
\\
\\
\leq M + \si ^2 B_c ( u_n) + \si ^2 ( c_1 - c) Q( u_n)  +  a_n, 
\end{array}
\eeq
where 
$$
a_n = c_1| Q( (u_n)_{1, \si }) - Q( v_n) | + 
\ds \int_{\R^3} \left\{ \frac 12  \left( \ph ^2( |1+ (u_n)_{1, \si} |) - 1 \right)^2 - 
V( |1+ (u_n)_{1, \si} |^2) \right\} dx .
$$
From (\ref{6.7}) and (\ref{6.11}) we infer that $\ds \lim_{ n \ra \infty } a_n = 0$. 
Then passing to the limit as $ n \lra \infty $ in (\ref{6.13}),  
using (\ref{6.8}) and the fact that $\ds \lim_{ n \ra \infty } B_c ( u_n) = 0$ we find
$
0 \leq M - \si ^2 \frac{ c_1 - c}{c} .
$
The last inequality clearly contradicts the choice of $ \si $ in (\ref{6.9}).
This contradiction shows that (\ref{6.6}) cannot hold and 
Lemma \ref{L6.3} is   proven. 
\hfill
$\Box $

\medskip

Next we show that we cannot have $ \al \in (0, \la _c + 1)$. 
We argue again by contradiction and we assume that $ \al \in ( 0 , \la _c + 1)$. 
Proceeding exactly as in the proof of Theorem \ref{T5.3} and using Lemma \ref{splitting}, 
we infer that for each $n$ sufficiently large  there exist two functions $ u_{n,1}$, $ u_{n,2}$
having the following properties: 
\beq
\label{6.14} 
E_{GL}(u_{n, 1}) \lra \al, 
\qquad E_{GL}(u_{n, 2}) \lra \la _c + 1 - \al , 
\eeq
\beq
\label{6.15}  
| A( u_n) - A( u_{n,1}) - A( u_{n,2}) | \lra 0 , 
\eeq
\beq
\label{6.16}  
| B_c( u_n) - B_c( u_{n,1}) - B_c( u_{n,2}) | \lra 0 , 
\eeq
\beq
\label{6.17}  
| D( u_n) - D( u_{n,1}) - D( u_{n,2}) | \lra 0  \qquad \mbox{ as } n \lra \infty.  
\eeq
Since $( E_{GL}(u_{n, i}))_{n \geq 1}$ 
are bounded,  from Lemmas \ref{L4.1} and \ref{L4.5} we see that 
$B_c (u_{n, i}))_{n \geq 1}$ 
are bounded. 
Moreover, by (\ref{6.16}) we have 
$ \ds \lim_{ n \ra \infty } \left ( B_c (u_{n, 1}) + B_c (u_{n, 2}) \right)
= \ds \lim_{ n \ra \infty }  B_c (u_{n}) = 0 . $ 
Similarly,  
$(D (u_{n, i}))_{n \geq 1}$ 
are bounded and  
$ \ds \lim_{ n \ra \infty } \left ( D (u_{n, 1}) + D (u_{n, 2}) \right)
= \ds \lim_{ n \ra \infty }  D(u_{n}) = 1 . $ 
Passing again to a subsequence (still denoted $(u_n)_{n \geq 1}$), 
we may assume that 
\beq
\label{6.18}
\ds \lim_{ n \ra \infty }  B_c (u_{n, 1}) = b_1, \qquad
\ds \lim_{ n \ra \infty }  B_c (u_{n, 2}) = b_2, \qquad 
\mbox{ where } b_i \in \R, \; b_1 + b_2 = 0, 
\eeq
\beq
\label{6.19}
\ds \lim_{ n \ra \infty }  D (u_{n, 1}) = d_1, \qquad
\ds \lim_{ n \ra \infty }  D (u_{n, 2}) = d_2, \qquad 
\mbox{ where } d_i \geq 0,  \; d_1 + d_2 = 1.
\eeq 
From (\ref{6.18}) it follows that either $ b_1 = b_2 = 0$, or one of $b_1$ or $b_2$ is negative. 

\medskip

{\it Case 1. } If $ b_1 = b_2 = 0$, we distinguish two subcases: 

\smallskip

{\it Subcase 1a. } We have $d_1 > 0$ and $ d_2 > 0$. 
Let $ \si _i = \frac {2}{ \sqrt{ d_i}}$, $ i=1, 2$. 
Then $ D((u_{n,i})_{1, \si _i }) = \si _i ^2 D(u_{n,i} ) \lra 4 $ and 
$ B_c((u_{n,i})_{1, \si _i }) = \si _i ^2 B_c (u_{n,i} ) \lra 0$ as $ n \lra \infty $. 
From (\ref{6.1}) and the definition of $ \la _c $ it follows that 
$ \ds \liminf_{n \ra \infty } A(u_{n,i}) =\ds \liminf_{n \ra \infty } A((u_{n,i})_{1, \si _i })  
\geq \la _c $, $i=1,2$. Then (\ref{6.15}) implies
$$
\ds \liminf_{n \ra \infty } A(u_{n }) \geq 
 \liminf_{n \ra \infty } A(u_{n,1}) +  \liminf_{n \ra \infty } A(u_{n,2})
\geq 2 \la _c
$$
and  this is a contradiction because by (\ref{6.5}) we have 
$\ds \lim_{n \ra \infty } A(u_{n })= \la _c$.

\smallskip

{\it Subcase 1b. }
One of the $d_i$'s is zero, say $ d_1 = 0$. 
Then necessarily $ d_2 = 1$, that is $\ds \lim_{n \ra \infty}  D( u_{n, 2}) = 1$. 
Since $E_{GL} ( u_{n,2}) = A( u_{n, 2}) + D( u_{n, 2}) \lra 1 + \la _c - \al $ as $ n \lra \infty $, 
we infer that $\ds \lim_{n \ra \infty } A(u_{n, 2 })= \la _c   - \al $. 
Hence $D(u_{n,2}) \lra 1 $, $B_c  (u_{n,2}) \lra 0 $ and $ A(u_{n,2}) \lra \la _c - \al $ 
as $ n \lra \infty$, which implies $ \la _c - \al \in \Lambda _c$. 
Since $\al >0$, this contradicts the definition of $ \la _c$. 

\medskip

{\it Case 2. } One of $b_i$'s is negative, say $ b_1 < 0$. 
From Lemma \ref{L4.8} (ii) we get $\ds \liminf_{n \ra \infty } A(u_{n, 1 }) > T_c \geq \la _c$
and then using (\ref{6.15}) we find $\ds \liminf_{n \ra \infty } A(u_{n }) > \la _c$, 
in contradiction with (\ref{6.5}). 

\medskip

Consequently in all cases we get a contradiction and this proves that we cannot have 
$ \al \in (0, \la _c + 1)$. 

Up to now we have proved that $ \ds \lim_{ t \ra \infty } q(t) = \la _c + 1$, 
that is "concentration" occurs. 

Proceeding as in the case $ N \geq 4$, we see that there exist a subsequence 
$ (u_{n_k})_{k \geq 1}$, a sequence of points $ (x_k ) _{k \geq 1} \subset \R^3$ 
and $ u \in \Xo $ such that, denoting $ \tilde{u}_{n_k} (x) = u_{n_k}( x + x_k)$, 
we have:
\beq
\label{6.20}
 \nabla \tilde{u}_{n_k} \rightharpoonup \nabla u  \mbox{ in } L^2( \R^3) \quad \mbox{ and } \quad 
\tilde{u}_{n_k}  \lra u  \mbox{ in } L_{loc}^p( \R^3) \mbox{ for } 1 \leq p < 6 
\mbox{ and a.e. on } \R^3 ,
\eeq
\beq
\label{6.21}
| 1 + \tilde{u}_{n_k} | - 1 \lra | 1 + u | -1 \qquad \mbox{ in } L^2( \R^3 ), 
\eeq
\beq
\label{6.22} 
\int_{ \R^3} V( | 1+  \tilde{u}_{n_k} |^2) \, dx \lra 
\int_{ \R^3} V( | 1+ u  |^2) \, dx , 
\eeq
\beq
\label{6.24}  
Q( \tilde{u}_{n_k} ) \lra Q(u) \qquad \mbox{ as } k \lra \infty. 
\eeq

Since $\big| \left( \ph ^2 ( s) -1  \right) ^2 - \left( \ph ^2 ( t) -1  \right) ^2 \big| 
\leq 24 |s-t|\left( | s-1| + |t-1| \right)$, 
from (\ref{6.21}) we get 
\beq
\label{6.23} 
\int_{ \R^3} 
\left( \ph ^2 ( | 1+ \tilde{u}_{n_k} | ) - 1 \right)^2 \, dx 
\lra 
\int_{ \R^3} 
\left( \ph ^2 ( | 1+ u | ) - 1 \right)^2 \, dx .
\eeq

Passing to the limit as $ k \lra \infty $ in the identity 
$$
\int_{ \R^3} \left\{ V( | 1+ \tilde{u}_{n_k} |^2) - \frac 12
\left( \ph ^2 ( | 1+ \tilde{u}_{n_k} | ) - 1 \right)^2 \right\} \, dx 
+ c Q( \tilde{u}_{n_k} )  
= B_c ( \tilde{u}_{n_k} )  - D ( \tilde{u}_{n_k} ) , 
$$
using (\ref{6.22})$-$(\ref{6.23}) 
and the fact that $  B_c ( \tilde{u}_{n_k} ) \lra 0$, $ D ( \tilde{u}_{n_k} ) \lra 1$ 
we obtain
$$
\int_{ \R^3} \left\{ V( | 1+ u |^2) - \frac 12
\left( \ph ^2 ( | 1+ u | ) - 1 \right)^2 \right\} \, dx  + cQ(u) = -1.
$$
Thus $ u \neq 0$. 

From the weak convergence $ \nabla \tilde{u}_{n_k} \rightharpoonup \nabla u  $ in $ L^2( \R^3)$ 
we get 
\beq
\label{6.25} 
\int_{ \R^3}  
\Big\vert \frac{ \p u}{ \p x_j } \Big\vert ^2 \, dx \leq 
\liminf_{ k \ra \infty } 
\int_{ \R^3}  
\Big\vert \frac{ \p \tilde{u}_{n_k}}{ \p x_j } \Big\vert ^2 \, dx 
\qquad \mbox{ for } j=1, 2, 3.
\eeq
In particular, we have 
\beq
\label{6.26}  
A(u) \leq \lim_{ k \ra \infty } A(\tilde{u}_{n_k} ) = \la _c.
\eeq
From (\ref{6.22}), (\ref{6.24}) and (\ref{6.25}) we obtain 
\beq
\label{6.27} 
B_c(u) \leq \lim_{ k \ra \infty }  B_c( \tilde{u}_{n_k} ) = 0.
\eeq
Since $ u \neq 0$, (\ref{6.27}) and Lemma \ref{L4.8} (i) imply 
$A(u) \geq T_c$. 
Then using (\ref{6.26}) and the fact that $ \la _c \leq T_c$, we infer that 
necessarily
\beq
\label{6.28}
A(u) = T_c = \la _c = \lim_{ k \ra \infty } A(\tilde{u}_{n_k} ).
\eeq

The fact that $B_c(\tilde{u}_{n_k} ) \lra 0 $, (\ref{6.22}) and (\ref{6.24}) imply that 
$\left( \ds \int_{ \R^3}  \Big\vert \frac{ \p \tilde{u}_{n_k}}{ \p x_1 } \Big\vert ^2 \, dx 
\right) _{k \geq 1}$
converges. 
If $\ds \int_{ \R^3}  
\Big\vert \frac{ \p u}{ \p x_1 } \Big\vert ^2 \, dx 
< \lim_{ k \ra \infty } 
\int_{ \R^3}  
\Big\vert \frac{ \p \tilde{u}_{n_k}}{ \p x_1 } \Big\vert ^2 \, dx  $,
we  get  $B_c(u) < \ds \lim_{ k \ra \infty }  B_c( \tilde{u}_{n_k} ) = 0$ in (\ref{6.27}) 
and then Lemma \ref{L4.8} (i)  implies $A(u) > T_c$, a contradiction. 
Taking (\ref{6.25}) into account, we see that necessarily
\beq
\label{6.29} 
\int_{ \R^3}  
\Big\vert \frac{ \p u}{ \p x_1 } \Big\vert ^2 \, dx 
= \lim_{ k \ra \infty } 
\int_{ \R^3}  
\Big\vert \frac{ \p \tilde{u}_{n_k}}{ \p x_1 } \Big\vert ^2 \, dx  
\qquad \mbox{ and } \qquad
B_c(u) = 0.
\eeq

Thus we have proved that $ u \in \Co $ and 
$\|\nabla u\|_{L^2(\R^3) } = \ds \lim_{ k \ra \infty } \| \nabla \tilde{u}_{n_k} \|_{L^2(\R^3) }.$
Combined with the weak convergence $  \nabla \tilde{u}_{n_k} \rightharpoonup \nabla u $ in $L^2( \R^3)$,  
this implies the strong convergence $  \nabla \tilde{u}_{n_k} \lra \nabla u $
in $L^2(\R^3) $ and the proof of Theorem \ref{T6.2} is complete. 
\hfill
$\Box $

\medskip

To prove that any minimizer provided by Theorem \ref{T6.2} satisfies an Euler-Lagrange 
equation, we will need the next lemma. 
It is clear that for any $ v \in \Xo $ and any $ R> 0$, the functional 
 $ \tilde{B}_c^v ( w): = B_c ( v + w)$ is $C^1$ on $H_0^1(B(0, R))$. 
We denote by $( \tilde{B}_c ^v )' ( 0).w 
= {\ds \lim_{t \ra 0}} \frac{B_c( v + tw) - B_c (v) }{t}$ its  derivative at the origin. 

\begin{Lemma}
\label{L6.4}
Assume that $ N \geq 3$ and the conditions (A1) and (A2) are satisfied. 
Let $ v \in \Xo$ be such that $ ( \tilde{B}_c ^v )' ( 0) .w = 0$ for any $ w \in C_c^1( \R^N)$. 
Then $ v = 0 $ almost everywhere in $ \R^N$. 
\end{Lemma}

{\it Proof. } 
We denote by  $ v^*$  the precise representative of $v$, that is 
$ v^*(x) = \ds \lim_{r \ra 0} m(v, B(x, r)) $ if this limit exists, and
$0$  otherwise. 
Since $ v \in L_{loc}^1( \R^N)$, it is well-known that $ v = v^*$ almost everywhere
 on $\R^N$ 
(see, e.g., Corollary 1 p. 44 in \cite{EG}). 
Throughout the proof of Lemma \ref{L6.4} we replace $v$ by $ v^*$. 
We proceed in three steps. 

\medskip

{\it Step 1. } There exists a set $ S \subset \R^{N-1}$ 
such that 
 $\Lo ^{N-1}(S) = 0$ and for any  $ x' \in \R^{N-1} \setminus S$ the function 
$ v_{x'} := v ( \cdot , x')$ belongs to $ C^2 ( \R)$ and solves the differential equation 
\beq 
\label{n1}
- (v_{x'})'' (s)  + i c (v_{x'}) ' (s) - F(|1+ v_{x'}(s)|^2) ( 1+ v_{x'}(s)) = 0 \qquad \mbox{ for any } s \in \R. 
\eeq  
Moreover, we have $ | v_{x'}(s) | \lra 0 $ as $ s \lra \pm \infty $ and $v_{x'}$ satisfies the following  properties:
\beq
\label{n3}
v_{x'} \in L^{2^*}( \R) ,  \qquad 
\ph ^2( |1+ v_{x'}|) - 1 \in L^2( \R) \quad \mbox{ and } \quad
\left( v_{x'} \right)' = \frac{\p v}{\p x_1} ( \cdot, x') \in L^2( \R), 
\eeq
\beq
\label{n4}
F(|1+ v_{x'}|^2) (1+ v_{x'} ) \in L^2(\R) + L^{\frac{ 2^*}{2 p_0 + 1}}( \R).
\eeq

It is easy to see that 
$ F(|1+ v|^2 )( 1+ v) \in L^2 (\R^N) + L^{\frac{ 2^*}{2 p_0 + 1}}( \R^N)$. 
Since $ v \in H_{loc }^1 ( \R^3)$, 
using Theorem 2 p. 164 in \cite{EG} and Fubini's Theorem, respectively, we see that 
there exists a set $ \tilde{S} \subset \R^{N-1}$ such that $ \Lo ^{N-1}(\tilde{S}) = 0 $ and 
for any $ x' \in \R^{N-1} \setminus \tilde{S}$ the function $ v_{x'}$ is  absolutely continuous, 
$v_{x'} \in H_{loc}^1( \R) $ and 
(\ref{n3})$-$(\ref{n4}) hold. 

Given  $ \phi \in C_c^1( \R)$,  we denote 
$$
 \Lambda _{\phi} ( x_1, x') = \langle \frac{ \p v}{\p x_1} (x_1, x'), \phi '(x_1) \rangle
+ c \langle i \frac{ \p v}{\p x_1} (x_1, x'), \phi (x_1) \rangle  
- \langle F( |1+ v|^2) ( 1+ v)(x_1, x') , \phi (x_1) \rangle,  
$$
where $\langle \cdot , \cdot \rangle $ denotes the scalar product of two complex numbers. 
From (\ref{n3}) and (\ref{n4}) it follows that $\Lambda_{\phi }(\cdot, x') \in L^1( \R)$ for 
$ x' \in \R^{N-1} \setminus \tilde{S}$. We define 
$ \la _{\phi }(x') = \ds \int_{\R} \Lambda _{\phi} ( x_1, x') d x_1$ if  $ x' \in \R^{N-1} \setminus \tilde{S}$ and 
 $ \la _{\phi }( x') = 0 $ if $ x' \in \tilde{S}$. 
Let $ \psi \in C_c^1( \R^{N-1})$. 
It is obvious that the function $ ( x_1, x') \longmapsto  \Lambda _{\phi} ( x_1, x') \psi ( x')$ belongs to 
$L^1( \R^N)$ and using Fubini's Theorem we get 
$$
 \ds \int_{\R^N} \Lambda _{\phi} ( x_1, x') \psi ( x')\,   d x = 
\ds \int_{\R^{N-1}} \la _{\phi }(x') \psi ( x')\, d x'. 
$$ 
On the other hand, using the assumption of Lemma \ref{L6.4} we obtain 
$$
 2 \ds \int_{\R^N} \Lambda _{\phi} ( x_1, x') \psi ( x')\,   d x = 
\left(\tilde{B}_c ^v \right)' ( 0) . (\phi( x_1) \psi ( x') ) = 0.
$$ 
Hence we have $  \ds \int_{\R^{N-1}} \la _{\phi }(x') \psi ( x')\, d x' = 0 $ for any $ \psi \in C_c^1( \R^{N-1})$ 
and this implies that there exists a set $ S_{\phi} \subset \R^{N-1} \setminus \tilde{S}$ such that 
$ \Lo^{N-1} ( S_{\phi }) = 0$ and $ \la _{\phi } = 0 $ on $\R^{N-1} \setminus   S_{\phi} $. 

Denote $ q _0 = \frac{ 2^*}{ 2 p_0 + 1} \in (1, \infty)$. 
There exists a countable set $ \{ \phi_n \in C_c^1( \R) \; | \; n \in \N \}$ which is dense in 
$ H^1 ( \R) \cap L^{ q_0 ' } ( \R)$. 
For each $n$ consider the set $ S_{\phi _n } \subset \R^{N-1}$ as   above.
Let $ S = \tilde{S} \cup \ds \bigcup_{n \in \N}  S_{\phi _n } $. It is clear that $\Lo ^{N-1}(S) = 0$. 

Let $ x' \in \R^{N-1} \setminus S$.
Fix $ \phi \in C_c^1( \R)$. There is a sequence $ (\phi_{n_k})_{k \geq 1} $ such that 
$\phi_{n_k} \lra \phi $ in $ H^1 ( \R)$ and in $L^{ q_0 ' } ( \R)$. 
 Then $ \la_{\phi _{n_k}} ( x') = 0$ for each $k$ and (\ref{n3})$-$(\ref{n4}) imply that 
$ \la_{\phi _{n_k}} ( x') \lra \la _{\phi} ( x')$. 
Consequently $ \la _{\phi }(x') = 0 $ for any $ \phi \in C_c^1( \R)$ and this implies that $ v_{x'}$ satisfies 
the equation  (\ref{n1}) in $ \Do ' ( \R)$. 
Using (\ref{n1}) we infer that $\left(v_{x'}\right)'' $ (the {\it weak } second derivative of  $v_{x'}$)
belongs to $ L_{loc}^1( \R)$ and  then
it follows that $\left(v_{x'}\right)' $ is continuous on $\R$ (see, e.g.,  Lemma VIII.2 p. 123 in \cite{brezis}).
In particular, we have $ v_{x'} \in C^1( \R)$. 
Coming back to (\ref{n1}) we see that $\left(v_{x'}\right)'' $ is continuous, hence $v_{x'} \in C^2( \R)$
and (\ref{n1}) holds at each point of $\R$. 
Finally, we have $|v_{x'}(s_2) - v_{x'}(s_1)| \leq |s_2 - s_1|^{\frac 12} \| \left(v_{x'}\right)' \|_{L^2}$; 
this estimate and the fact that $ v_{x'} \in L^{2^*}(\R)$ imply that $ v_{x'}(s) \lra 0 $ as $ s \lra \pm \infty$. 

\medskip

{\it Step 2. } There exist two positive constants $ k_1, \; k_2$ (depending only on $F$ and $c$) 
such that for any $ x' \in \R^{N-1} \setminus S$ we have either $v_{x'} = 0$ on $\R$ or there exists 
an interval $I_{x'} \subset \R$ with $\Lo ^1(I_{x'}) \geq k_1 $ 
and $\big|\, | 1+ v_{x'}| - 1 \big|  \geq k_2 $ on $I_{x'}$. 

\smallskip

To see this, fix $ x' \in \R^{N-1} \setminus S$ and denote $ g = | 1+ v_{x'}|^2 - 1$. 
Then $ g \in C^2(\R, \R)$  and $g$ tends to zero at $ \pm \infty$. 
Proceeding exactly as in \cite{M8}, p. 1100-1101 we integrate (\ref{n1}) and we see that $g$ satisfies
\beq
\label{n5}
(g')^2 (s) + c^2 g^2 (s) - 4 ( g (s)+ 1) V( g(s) + 1) = 0 \qquad \mbox{ in } \R. 
\eeq
Using (\ref{1.4}) we have $ c^2 t^2 - 4( t+ 1) V( t + 1) = t^2 ( c^2 - v_s ^2 + \e _1 (t))$, 
where $ \e_1 (t) \lra 0 $ as $ t \lra 0$. 
In particular, there exists $ k_0 > 0$ such that 
\beq
\label{n6}
c^2 t^2 - 4( t+ 1) V( t + 1 )  < 0 \qquad \mbox{ for } t \in [-2 k_0, 0) \cup (0, 2 k_0].
\eeq

If $ g =0$ on $ \R$ then $ |1+ v_{x'} | = 1 $ and consequently there exists a lifting 
$1+ v_{x'} (s) =  e^{i \theta (s)}$ with $ \theta \in C^2( \R, \R)$. 
Using equation (\ref{n1}) and proceeding as in \cite{M8} p. 1101 we see that either 
$1+ v_{x'}(s)  =  e^{i \theta _0}$ or $1+ v_{x'}(s)  =  e^{i c s + \theta _0}$, 
where $ \theta _0 \in \R$ is a constant. 
Since $ v_{x'} \in L^{2^*}(\R)$, we must have  $ v_{x'} = 0$. 

If $ g \not\equiv 0$, the function $g$ achieves a negative minimum or a positive maximum at some $ s_0 \in \R$. 
Then $ g'( s_0) = 0 $ and using (\ref{n5}) and (\ref{n6}) we infer that $|g(s_0)| > 2 k_0$. 
Let 
$$
 s_2 = \ds \inf \{ s < s_0 \; \big| \; |g(s)| \geq 2 k_0 \},  \qquad 
 s_1 = \sup \{ s < s_2 \; \big| \; |g(s )|  \leq k_0 \}, 
$$
 so that $ s_1 < s_2$, 
$ |g( s_1)| = k_0 $, $ |g(s_2)| = 2 k_0$ and $ k_0 \leq | g(s) | \leq 2 k_0$ for $ s \in [s_1, s _2]$. 
Denote $ M = \ds \sup \{ 4( t+ 1) V( t + 1) - c^2 t^2 \; | \; t \in [ -2 k_0, 2 k_0] \}$. 
From (\ref{n5}) we obtain $|g'(s)| \leq \sqrt{M}$ if $ g(s) \in [-2 k_0, 2k_0 ]$ and we infer that 
$$
k_0 = | g( s_2)| - | g( s_1)| \leq \Big\vert \int_{s_1}^{s_2} g'(s) \, ds \Big\vert \leq \sqrt{M} ( s_2 - s_1), 
$$
hence $ s_2 - s _1 \geq \frac {k_0}{\sqrt{M}}$. 
Obviously,  there exists $ k_2 > 0$ such that $ | \, | 1+ z |^2 - 1| \geq k_0 $ implies 
$|\, | 1+z| - 1 | \geq k_2$. 
Taking  $ k_1 = \frac{k_0}{\sqrt{M}}$ and $ I_{x'}= [s_1, s_2]$,   the proof of step 2 is complete. 

\medskip

{\it Step 3. } Conclusion. 

\smallskip

Let $ K = \{ x' \in \R^{N-1} \setminus S \; | \; v_{x'} \not\equiv 0 \}$. 
It is standard to prove that $K$ is $\Lo ^{ N-1}-$measurable. 
The conclusion of Lemma \ref{L6.4} follows if we prove that $ \Lo ^{N-1}(K) = 0$. 
We argue by contradiction and we assume that $ \Lo ^{N-1}(K) > 0$. 

If $ x' \in K$, it follows from step 2 that there exists an interval $I_{x'}$ of length at least $k_1$ such that 
$ \left( \ph ^2 ( |1+ v_{x'}| ) - 1 \right)^2 \geq \eta ( k_2 ) $ on $I_{x'}$, where $ \eta$ 
is as in (\ref{3.30}). 
This implies 
$ \ds \int_{\R}   \left( \ph ^2 ( |1+ v ( x_1, x')| ) - 1 \right)^2 \, d x_1 \geq k_1 \eta ( k_2)$ 
and using Fubini's theorem we get 
$$
\begin{array}{l}
 \ds \int_{\R^N}   \left( \ph ^2 ( |1+ v ( x)| ) - 1 \right)^2 \, d x
= \int_K \left(  \ds \int_{\R}   \left( \ph ^2 ( |1+ v ( x_1, x')| ) - 1 \right)^2 \, d x_1  \right) d x'
\\
\geq k_1 \eta ( k_2) \Lo ^{N-1}(K).
\end{array}
$$
Since $ v \in \Xo$, we infer that $ \Lo ^{N-1}(K) $ is finite. 

It is obvious that there exist $ x_1' \in K$ and $ x_2 ' \in \R^{N-1} \setminus (K\cup S)$ 
arbitrarily close to each other. 
Then $|v_{x_1 '}| \geq k_2$ on an interval $I_{x_1 '}$ of length $ k_1$, while   $ v_{x_2 '} \equiv 0$. 
If we knew that $v$ is uniformly continuous, this would lead to a contradiction. 
However, the equation (\ref{n1}) satisfied by $v$ involves only derivatives with respect to $ x_1$ and 
does not imply any regularity properties of $v$ with respect to the transverse variables 
(notice that if $v$ is a solution of (\ref{n1}),  
then $ v( x_1 + \de( x'), x')$ is also a solution, even if $\de$ is discontinuous). 
For instance, for the Gross-Pitaevskii nonlinearity $F(s) = 1-s$ it is possible to construct bounded, $C^{\infty }$ functions 
$v$ such that $ v \in L^{2^*}(\R^N)$,  (\ref{n1}) is satisfied for a.e. $x'$, and the set $K$ 
constructed as above is a nontrivial ball in $\R^{N-1}$
(of course, these functions do not tend uniformly to zero at infinity,  
are not uniformly continuous and their gradient is not in $L^2(\R^N)$). 

We  use the fact that one transverse derivative of $v$ (for instance, $\frac{\p v}{\p x_2}$)
 is in $ L^2( \R^N)$ to get a contradiction.

For $ x' = ( x_2, x_3, \dots, x_N) \in \R^{N-1}$, we denote $ x'' = ( x_3, \dots, x_N)$. 
Since $ v \in H_{loc}^1(\R^N)$, from Theorem 2 p. 164 in \cite{EG} it follows that there exists 
$J \subset \R^{N-1}$ such that $ \Lo^{N-1} (J) =0$ and $ u( x_1, \cdot, x'') \in H_{loc}^1( \R^N)$ 
for any $ (x_1, x'') \in \R^{N-1} \setminus J$.
Given $ x'' \in \R^{N-2}$, we denote 
$$
\begin{array}{l}
K_{x''} = \{ x_2 \in \R \; | \; ( x_2, x'' ) \in K \}, \\
S_{x''} = \{ x_2 \in \R \; | \; ( x_2, x'' ) \in S \}, \\
J_{x''} = \{ x_1 \in \R \; | \; ( x_1, x'' ) \in J \}.
\end{array}
$$ 
Fubini's Theorem implies that for almost all $ x'' \in \R^{N-2}$,  the sets 
$ K_{x''}$, $ S_{x''}$, $ J_{x''}$ are $\Lo ^1-$measurable,  $\Lo ^1(  K_{x''} ) < \infty $ and 
$ \Lo ^1(S_{x''} ) =  \Lo ^1(J_{x''} ) = 0$.  
Let 
\beq
\label{n7}
\begin{array}{ccl}
G & = & \{ x'' \in \R^{N-2} \; | \; K_{x''}, S_{x''}, J_{x''} \mbox{ are } \Lo ^1 \mbox{ measurable, }
\\
& & \Lo ^1(S_{x''} ) =  \Lo ^1(J_{x''} ) = 0 \mbox{ and } 0< \Lo ^1(K_{x''}) < \infty \}.
\end{array}
\eeq
Clearly, $G$ is $ \Lo^{N-2}-$measurable and $\ds \int_G \Lo ^1(K_{x''}) \, d x'' = \Lo^{N-1}(K) > 0$, 
thus $ \Lo^{N-2}(G) >0$. We claim that 
\beq
\label{n8}
\int_{\R^2} \Big\vert \frac{ \p v}{\p x_2} ( x_1, x_2, x'' )\Big\vert ^2 dx_1 \, dx_2 = \infty 
\qquad \mbox{ for any } x'' \in G.
\eeq

Indeed, let $ x'' \in G$. Fix $ \e > 0$. 
Using  (\ref{n7}) we infer that there exist $s_1, \,  s_2 \in \R$ such that 
$(s_1, x'' ) \in \R^{N-1} \setminus (K \cup S)$,  $ (s_2, x'' ) \in  K$ and $|s_2 - s_1 | < \e$. 
Then $v( t, s_1, x'' ) = 0$ for any $ t \in \R$. 
From step 2 it follows that there exists an interval $I$ with $ \Lo ^1( I) \geq k_1$ such that 
$| v(  t, s_2, x'' )| \geq | \, | 1+  v(  t, s_2, x'' )| - 1 | \geq k_2$ for $ t \in I$. 
Assume $ s_1 < s_2$.
If $ t \in I \setminus J_{x''} $ we have $ v(t, \cdot, x'') \in H_{loc}^1( \R)$, hence 
$$
\begin{array}{l}
k_2 \leq |v(  t, s_2, x'' ) - v(  t, s_1, x'' )|
= \ds \Big\vert { \int_{s_1}^{s_2}} \frac{ \p v }{\p x_2} ( t, \tau, x'') \, d \tau \Big\vert
\\
\ds \leq ( s_2 - s_1)^{\frac 12} \left( \int_{s_1}^{s_2} \Big\vert \frac{ \p v }{\p x_2} ( t, \tau, x'') \Big\vert ^2 d \tau \right)^{\frac 12}.
\end{array}
$$
Clearly, this implies 
$\ds  \int_{s_1}^{s_2}  \Big\vert \frac{ \p v }{\p x_2} ( t, \tau, x'') \Big\vert ^2 d \tau \geq \frac{ k_2 ^2}{\e }$. 
Consequently
$$
\int_{\R^2} \Big\vert \frac{ \p v}{\p x_2} ( x_1, x_2, x'' )\Big\vert ^2 dx_1 \, dx_2 
\geq \int_I  \int_{s_1}^{s_2} \Big\vert \frac{ \p v }{\p x_2} ( t, \tau, x'') \Big\vert ^2 d \tau  \, dt 
\geq \frac{ k_1 k_2 ^2}{\e }.
$$
Since the last inequality holds for any $ \e > 0$, (\ref{n8}) is proven. 
Using  (\ref{n8}), the fact that $ \Lo^{N-2}(G) >0$ and Fubini's Theorem we get 
$ \ds \int_{\R^N} \Big\vert \frac{ \p v }{\p x_2} \Big\vert ^2 \, dx = \infty $, 
contradicting the fact that $ v \in \Xo$. 
Thus necessarily $ \Lo ^{N-1}(K) = 0 $ and the proof of Lemma \ref{L6.4} is complete.
\hfill
$\Box $

\begin{Proposition} 
\label{P6.4}
Assume that $ N =3$  and the conditions (A1)  and (A2)  are satisfied. 
Let $ u \in \Co $ be a minimizer of $ E_c$ in $ \Co$. 
Then $ u \in W_{loc}^{2, p }(\R^3)$ for any $ p \in [1, \infty)$,  
$ \nabla u \in W^{1, p }(\R^3)$ for  $ p \in [2, \infty)$ 
and there exists $ \si > 0$ such that $u_{1, \si} $ is a solution of  (\ref{1.3}).
\end{Proposition} 

{\it Proof. } The proof is very similar to the proof of Proposition \ref{P5.6}. 
It is clear that $A(u) = E_c(u) = T_c$ and $u$ is a minimizer of $ A$ in $ \Co$. 
For any $ R>0$, the functionals $\tilde{B}_c^u$ and $\tilde{A}(v): = A(u + v)$
are $C^1$ on $H_0^1( B(0, R))$. 
We proceed in four steps. 

\medskip

{\it Step 1. } There exists $ w \in C_c^1( \R^3 )$ such that 
$(\tilde{B}_c ^u)'(0). w \neq 0$. 
This follows from Lemma \ref{L6.4}.

\medskip

{\it Step 2. } There exists a Lagrange multiplier $ \al \in \R$ such that 
\beq
\label{6.32}
\tilde{A} ' (0). v = \al (\tilde{B}_c ^u) ' (0) . v 
\qquad \mbox{ for any } v \in  H^1( \R^3 ) , \; v \mbox{ with compact support.}
\eeq


{\it Step 3. } We have $ \al < 0$. 

\smallskip

\noindent
The proof of steps 2 and 3
is the same as 
the proof of steps 2 and 3 in Proposition \ref{P5.6}.

\medskip 

{\it Step 4. } Conclusion. \\
Let $ \beta = - \frac{1}{\al }$. 
Then (\ref{6.32}) implies that $u$ satisfies 
$$
- \frac{ \p^2 u }{ \p x_1 ^2 } 
- \beta \left( \frac{ \p^2 u }{ \p x_2 ^2 }  + \frac{ \p^2 u }{ \p x_3 ^2 } \right) 
+  i c u_{x_1} - F(|1+  u|^2) ( 1+ u ) = 0 
\mbox{  in } \Do ' ( \R^3) . 
$$
For $ \si ^2 = \frac{ 1}{\beta }$ we see that $ u_{1, \si} $ satisfies (\ref{1.3}). 
It is clear that $ u_{1, \si} \in \Co$ and $u_{1, \si} $ minimizes $A$ (respectively $E_c$)
in $ \Co$. 
Finally, the regularity of $u_{1, \si} $ (thus the regularity of $u$) 
follows from Lemma \ref{L5.5}.
\hfill
$\Box $

\section{Further properties of traveling waves } 

By Propositions \ref{P5.6} and \ref{P6.4} we already know that the solutions of 
(\ref{1.3}) found there are in $W_{loc}^{2, p} (\R^N)$ for any $ p \in [1, \infty)$ and in $ C^{1, \al } (\R^N) $ for any $ \al \in (0, 1)$. 
 In general, a straightforward bootstrap argument shows that 
the finite energy traveling waves of (\ref{1.1}) have the best regularity allowed by 
the nonlinearity $F$. 
For instance, if $F \in C^k ( [0, \infty ))$ for some $ k \in \N^*$, 
it can  be proved 
that all finite energy solutions of (\ref{1.3}) are in $W_{loc}^{k+2, p}(\R^N)$ for any 
$ p \in [1, \infty )$ (see, for instance, Proposition 2.2 (ii) in \cite{M8}). 
If $F$ is analytic, it can be proved that finite energy traveling waves 
 are also analytic. 
In the case of the Gross-Pitaevskii equation, this has been done in \cite{BGS}.

\medskip

Our next result concerns the symmetry of those solutions of (\ref{1.3}) that  minimize 
 $E_c$ in $\Co$.

\begin{Proposition}
\label{P7.1} 
Assume that $N\geq 3$ and the  conditions (A1), (A2)  in the introduction hold. 
Let $ u \in \Co $ be a minimizer of $ E_c$ in  $ \Co$. 
Then, after a translation in the variables $(x_2, \dots, x_N)$, 
$u$ is axially symmetric with respect to $Ox_1$. 
\end{Proposition}

{\it Proof. } 
Let $T_c$ be as in Lemma \ref{L4.7}.
We know that any minimizer $u$ of $E_c$ in $ \Co $ satisfies $ A( u) = \frac{N-1}{2} T_c > 0$.
Using Lemma \ref{L4.8} (i), it is easy to prove that a function $ u \in \Xo $ is a minimizer 
of $E_c$ in $ \Co $ if and only if 
\beq
\label{7.1}
u \mbox{ minimizes the functional  }  P_c \mbox{ in the set } 
\{ v \in \Xo \; | \; A(v) = \frac{N-1}{2} T_c \} . 
\eeq

The minimization problem (\ref{7.1}) is of the type studied in \cite{M7}. 
All we have to do is to verify that the assumptions ({\bf A1$_c$}) and ({\bf A2$_c$})  in \cite{M7}, p. 329 
are satisfied, then to apply the general theory developed there. 

Let $ \Pi $ be an affine hyperplane in $\R^N$ parallel to $Ox_1$. 
We denote by $s_{\Pi}$ the symmetry of $ \R^N$ with respect to $ \Pi $ and by 
$\Pi ^+$, $\Pi ^-$ the two half-spaces determined by $ \Pi $. 
Given a function $ v \in \Xo$, we denote
$$
v_{\Pi ^+}(x) = \left\{ 
\begin{array}{ll}
v(x) & \mbox{if } x \in \Pi ^+ \cup \Pi, 
\\
v(s_{\Pi }(x)) & \mbox{if } x \in \Pi ^ -, 
\end{array}
\right. 
\quad \mbox{ and } \quad 
v_{\Pi ^-}(x) = \left\{ 
\begin{array}{ll}
v(x) & \mbox{if } x \in \Pi ^- \cup \Pi, 
\\
v(s_{\Pi }(x)) & \mbox{if } x \in \Pi ^ + . 
\end{array}
\right.  
$$
It is easy to see that $ v_{\Pi ^+} , v_{\Pi ^-} \in \Xo $.
Moreover, for any $ v \in \Xo $ we have
$$
A( v_{\Pi ^+}) + A( v_{\Pi ^-}) = 2 A(v) 
\quad \mbox{ and } \quad 
P_c ( v_{\Pi ^+}) + P_c ( v_{\Pi ^-}) = 2 P_c (v) .
$$
This implies that assumption ({\bf A1$_c$}) in \cite{M7} is satisfied. 

By Propositions \ref{P5.6} and \ref{P6.4} and Lemma \ref{L5.5}  we know that any minimizer of (\ref{7.1}) 
is $C^1$ on $\R^N$, hence  assumption ({\bf A2$_c$}) in \cite{M7}  holds. 
Then the axial symmetry of solutions of (\ref{7.1})  follows directly from Theorem 2' p. 329 in 
\cite{M7}.
\hfill
$\Box $

\bigskip

\noindent
{\it Acknowledgements. } I am indebted to Jean-Claude Saut, 
who brought this subject to my attention 
almost ten  years ago, for many  useful and stimulating  conversations. 
I am also very grateful to Radu Ignat, Cl\'ement Gallo, 
Alberto Farina,  David Chiron  and Mihai Bostan for interesting  discussions 
and to the anonymous referees who suggested  improvements in the exposition. 

This paper was written when I was a member of the  Department of Mathematics at 
Universit\'e de Franche-Comt\'e in Besan\c con, France. 
I would like to thank people at the Department of Mathematics in Besan\c con for their hospitality.

{}

\end{document}